\documentclass[twoside,a4paper]{amsart}
\usepackage[centertags]{amsmath}
\usepackage{amsfonts}
\usepackage{dsfont}
\usepackage{amssymb}
\usepackage{amsthm}
\usepackage[ansinew]{inputenc}
\usepackage[cmtip,all,poly]{xy}
\usepackage{graphicx}
\newtheorem{thm}[equation]{Theorem}
\newtheorem{cor}[equation]{Corollary}
\newtheorem{lem}[equation]{Lemma}
\newtheorem{prop}[equation]{Proposition}
\theoremstyle{definition}
\newtheorem{defn}[equation]{Definition}
\theoremstyle{remark}
\newtheorem{rem}[equation]{Remark}
\newtheorem{exm}[equation]{Example}
\numberwithin{equation}{section}

\newcommand{\abs}[1]{\left\vert#1\right\vert}
\newcommand{\set}[1]{\left\{#1\right\}}
\newcommand{\Real}{\mathbb R}

\newcommand{\To}{\longrightarrow}

\newcommand{\C}[1]{\mathbf{#1}} 

\def\op{{op}}
\def\nill{{nil}}
\def\abb{{ab}}
\def\add{{add}}
\def\adc{{adc}}

\def\r{\rightarrow} 
\def\rr{\Rightarrow} 

\def\hom{\operatorname{Hom}}

\newcommand{\sym}[1]{\operatorname{Sym}(#1)}
\newcommand{\symm}{\operatorname{Sym}}

\newcommand{\symt}[1]{\operatorname{Sym}_\vc(#1)}
\newcommand{\symtt}{\operatorname{Sym}_\vc}

\def\sign{\operatorname{sign}}

\def\aut{\operatorname{Aut}}

\def\endo{\operatorname{End}}
\def\fc{\mathrm{fc}}

\def\der{\mathrm{Der}}
\def\inn{\mathrm{Inn}}

\def\wc{\ul{\wedge}}
\def\wa{\wedge}

\def\st{\stackrel} 

\def\ul{\underline} 
\def\ol{\overline} 

\newcommand{\hopf}{\mathit{Hopf}}

\def\colim{\mathop{\operatorname{colim}}}

\def\coker{\operatorname{Coker}}
\renewcommand{\ker}{\operatorname{Ker}}

\def\Z{\mathbb{Z}}

\def\N{\mathds{N}}

\def\S{\Sigma}

\def\L{\Omega}

\newcommand{\grupo}[1]{\langle #1\rangle}

\newcommand{\vc}{\Box}    
\newcommand{\vi}{\boxminus}  

\begin{document}

\title{The algebra of secondary homotopy operations in ring spectra}%
\author{Hans-Joachim Baues and Fernando Muro}%

\address{Max-Planck-Institut f\"ur Mathematik, Vivatsgasse 7, 53111 Bonn, Germany}%
\email{baues@mpim-bonn.mpg.de}%
\address{Universitat de Barcelona, Departament d'$\grave{\text{A}}$lgebra i Geometria, Gran Via de les Corts Catalanes, 585, 08007 Barcelona, Spain}
\email{fmuro@ub.edu}

\thanks{The second author was partially supported
by the project MTM2004-01865 and the MEC postdoctoral fellowship EX2004-0616.}%

\date{\today}
\begin{abstract}
The primary algebraic model of a ring spectrum $R$ is the ring $\pi_*R$ of homotopy groups. We introduce the
secondary model $\pi_{*,*}R$ which has the structure of a secondary analogue of a ring. The homology of
$\pi_{*,*}R$ is $\pi_*R$ and triple Massey products in $\pi_{*,*}R$ coincide with Toda brackets in $\pi_*R$. We
also describe the secondary model of a commutative ring spectrum $Q$ from which we derive the cup-one square
operation in $\pi_*Q$. As an application we obtain for each ring spectrum $R$ new derivations of the ring
$\pi_*R$.
\end{abstract}
\maketitle
\begin{footnotesize}
\tableofcontents
\end{footnotesize}

\section*{Introduction}

A ring spectrum $R$ is a topological analogue of a ring. The primary algebraic model of $R$ is the ring $\pi_*R$
of its homotopy groups. We study in this paper secondary homotopy operations in $\pi_*R$ which lead to the
secondary model $\pi_{*,*}R$ extending the ring $\pi_*R$. The model $\pi_{*,*}R$ has the structure of a
\emph{secondary algebra} which generalizes the notion of algebra in the same sort of way as a $2$-group, or
crossed module in the sense of Whitehead \cite{chII}, generalizes the classical notion of group by adding a
$2$-dimensional part. The secondary algebras needed are pair algebras (resp. $E_\infty$-pair algebras) and 
quadratic pair algebras (resp. $E_\infty$-quadratic pair algebras). They are in a natural way secondary models of
ring spectra (resp. commutative ring spectra) and describe new algebraic porperties enriching the commonly used
paradigm of ring.

A \emph{pair algebra} $B$ is a ring $B_0$ together with a $B_0$-bimodule morphism
\begin{equation*}
\partial\colon B_1\To B_0
\end{equation*}
satisfying $\partial (a)b=a\partial(b)$ for $a,b\in B_1$. Here $B_0$, $B_1$ and $\partial$ are $\N$-graded. In the ungraded case a
pair algebra is also termed a \emph{crossed bimodule}, see \cite[E.1.5.1]{ch}. Pair algebras or crossed bimodules
are known to represent elements in Hochschild and Shukla cohomology, see \cite{ch, chsml}. Moreover, quadratic
pair algebras represent elements in Mac Lane cohomology \cite{3mlc}. The cohomology class 
$$\grupo{\pi_{*,*} R}\in HML^3(\pi_*R,\Sigma^{-1}\pi_*R)$$
represented by $\pi_{*,*}R$ is the universal Toda bracket (\cite{ccglg,steffen}) which represents the homotopy
category of free $R$-modules as a linear extension of categories.

It is well known that homotopy groups $\pi_*$ form a lax symmetric monoidal functor carrying the symmetric monoidal
category of spectra to the symmetric monoidal category of abelian groups. This, in fact, implies that homotopy
groups of a ring spectrum $R$ form a ring. In a similar way we show that there if a lax symmetric monoidal
functor $\pi_{*,*}$ on $\L$-spectra, given by secondary homotopy groups, which yields the structure of a
secondary algebra on $\pi_{*,*}R$. We do not know of any other lax symmetric monoidal functor in the literature
leading to such algebraic models of ring spectra extending the ring $\pi_*R$. 

A functor in the opposite direction is
studied by Shipley in \cite{hzas} yielding a ring spectrum $\mathbb{H}D_*$ for a differential graded algebra
$D_*$. In this case $\pi_{*,*}\mathbb{H}D_*$ can be described by the secondary algebra associated to $D_*$, see
Remark \ref{exmpa}. Moreover, ungraded quadratic pair algebras yield examples of rings and commutative rings in
permutative categories, see Remarks \ref{per2} and \ref{per3}, and Elmendorf and Mandell define in \cite{rmailst} a
functor carrying such objects to ring spectra.

Algebraic Massey products in the secondary algebra $\pi_{*,*}R$ coincide with Toda brackets in $\pi_*R$ (Theorems
\ref{a1} and \ref{a2}). Moreover, if $R$ is commutative the algebraic cup-one square in $\pi_{*,*}R$ coincides
with the topologycal cup-one square in $\pi_*R$ (Theorems \ref{ca1} and \ref{ca2}).

The new invariant $\pi_{*,*}R$ is strictly stronger than the collection of secondary operations in $\pi_*R$ as we
illustrate in Example \ref{3stem} by the $3$-local sphere spectrum. This shows examplary computations in the
algebraic model $\pi_{*,*}R$.

In case $R_p$ is the endomorphism ring spectrum of the Eilenberg-Mac Lane spectrum $H\Z/p$ the ring
$\pi_*R_p=\mathcal{A}$ is the mod $p$ Steenrod algebra. In \cite{asco} the pair algebra of secondary cohomology
operations $\mathcal{B}$ is computed which corresponds to $\pi_{*,*}R_p$. It is the purpose of this paper to
achieve analogous pair algebras associated to arbitrary ring spectra. 
As an example the $0$-dimensional part $\pi_{0,*}K\C{W}$ of Waldhausen $K$-theory is computed in \cite{1tk}.

The sphere spectrum $S$ (which is a connective commutative ring spectrum) yields the
$E_\infty$-quadratic pair algebra $\pi_{*,*}S$ which enriches the structure of the commutative algebra
$\pi_*S$ considerably. This $E_\infty$-quadratic pair algebra acts on $\pi_{*,*}R$ for any other ring spectrum
$R$. The laws of this action (see Definition \ref{meqpa}) allow the construction of the derivations $\theta(b)$
as follows. Let $b\in\pi_nS$ be an element which maps to zero in $\pi_nR$ by the unit map $u\colon S\r R$ of the
ring spectrum $R$. Then a derivation of degree $n$
$$\theta(b)\colon\pi_*R\To\S^{-1}\pi_*R$$
is constructued which is well defined up to inner derivation. Moreover $\theta(b)$ is 
$\pi_*S$-linear, so it determines a Hochschild cohomology
class 
\begin{eqnarray*}
\theta(b)&\in&HH^1_{\pi_*S}(\pi_*R,\S^{-1}\pi_*R).
\end{eqnarray*}
This class vanishes in case $R$ is commutative. The properties of $\theta$ are described in Theorem
\ref{app}.

Special cases of the derivation $\theta(b)$ are known in the literature. For example for the endomorphism
spectrum of $H\Z/p$ one has $u_*(p\cdot 1)=0$ for $1\in\pi_0S$ and
$$\theta(p\cdot 1)=\kappa\colon\mathcal{A}\To\S\mathcal{A}$$
is the Kristensen derivation (\cite{osco}) for $p=2$ which carries $Sq^n$ to $Sq^{n-1}$. For odd primes $\kappa$
is computed in \cite{asco} by $\kappa(P^n)=0$ and $\kappa(\beta)=1$. Moreover, if $R$ is the endomorphism
spectrum of a $\Z/p$-space in the sense of Toda then $\theta(p\cdot 1)$ for $p$ odd coincides
with the derivation constructed by Toda in \cite{deritoda}.

In order to introduce the reader smoothly to the new theory we begin this paper with a section recalling 
classical secondary homotopy operations  and defining the derivations
$\theta(b)$ in a topological language.
Afterwards we give linear versions of our main theorems for the case of ring spectra neglecting the Hopf map.
This will help the reader to understand the more general quadratic versions
which are needed to deal with arbitrary spectra.

We are mainly concerned in this paper with connective spectra in order to avoid further technicalities.

\section{Ring spectra and module spectra}\label{1}

In this paper the framework for stable homotopy theory will be the stable model category of symmetric spectra of
compactly generated topological spaces defined in \cite[9]{mcds}. The smash product of symmetric spectra $X\wedge
Y$ defines a symmetric monoidal structure in this category. 
The unit of this monoidal structure is the \emph{sphere spectrum} $S$.
Monoids in the category of symmetric spectra are called ring spectra. 
Ring spectra, modules over a ring spectrum, and algebras over a commutative ring spectrum also form model categories, see \cite[12]{mcds}. Fibrant objects in
all these categories coincide with the objects which are $\L$-spectra. This crucial fact allows the
construction of all the new algebraic invariants presented in this paper. 

The homotopy groups of a ring spectrum $\pi_*R$ carry Toda bracket operations which enrich the
ring structure of $\pi_*R$. Toda brackets have been considered, for instance, in \cite{toda} for the sphere spectrum $S$ and
in \cite{cmp}, under the name of Massey products, for various cobordism spectra. 

The homotopy group $\pi_nR$ coincides with the group of morphisms $\S^nR\r R$ in the stable homotopy category of right $R$-modules.
Therefore 
\emph{Toda brackets} in $\pi_*R$ are defined
following \cite{shc} by using the triangulated structure in the homotopy category of right $R$-modules, see also \cite{steffen}. More
precisely, given elements $a,b,c\in\pi_*R$ of degree $p,q,r$, respectively, with $ab=0$ and $bc=0$ a generic element in the Toda
bracket $\grupo{a,b,c}\subset\pi_{p+q+r+1}R$ is a morphism $g$ in the homotopy category of right $R$-modules fitting into
the commutative diagram
\begin{equation}\label{toda1}
\xymatrix{\S^{p+q+r+1}R\ar[d]_g&C\ar[l]\ar[d]&\S^{p+q}R\ar[l]\ar@{=}[d]&\S^{p+q+r}R\ar[l]_{\S^{p+q}c}\ar@{=}[d]\\
R&\S^pR\ar[l]^a&\S^{p+q}R\ar[l]^{\S^pb}&\S^{p+q+r}R\ar[l]^{\S^{p+q}c}}
\end{equation}
where the upper row is an exact triangle. The existence of such a commutative diagram follows from the axioms of a triangulated
category.

One can also define Toda brackets in $\pi_*R$ following \cite{ccglg} by using \emph{tracks} in the category of
right $R$-modules (i.e. homotopy classes of homotopies between maps). Let us sketch this alternative construction.
Since the homotopy groups of any fibrant replacement of $R$ are isomorphic to $\pi_*R$ we can suppose without
loss of generality that $R$ is a fibrant ring spectrum. In that case the lower row of diagram (\ref{toda1}) can
be realized by a diagram in the category of right $R$-modules
\begin{equation*}
\xymatrix{R&\S^pR\ar[l]_{\bar{a}}&\S^{p+q}R\ar[l]_{\S^p \bar{b}}&\S^{p+q+r}R.\ar[l]_{\S^{p+q}\bar{c}}}
\end{equation*}
The vanishing hypotheses $ab=0$ and $bc=0$ imply the existence of null-homotopies
\begin{equation}\label{toda2}
\xymatrix{R&\S^pR\ar[l]_{\bar{a}}&\S^{p+q}R\ar[l]|{\S^p \bar{b}}^<(.85){\;}="c"\ar@/^30pt/[ll]^0_{\;}="d"&
\S^{p+q+r}R\ar[l]^{\S^{p+q}\bar{c}}_<(.88){\;}="a"\ar@/_30pt/[ll]_0^{\;}="b"\ar@{=>}"a";"b"^f\ar@{=>}"c";"d"^e}.
\end{equation}
The pasting of this diagram is a self-track of the trivial map $0\colon \S^{p+q+r}R\r R$.
Such a self-track is the same as a homotopy class
$$g\colon\S^{p+q+r+1}R\To R,$$
which is again a generic element of the Toda bracket $\grupo{a,b,c}$.
This is the more convenient approach from the
perspective of this paper. More general Toda brackets for a right $R$-module $M$,
$\grupo{a,b,c}\subset\pi_{p+q+r+1}M$, are defined simply by replacing
$R$  by $M$ on the lower left corner of (\ref{toda1}) or
on the left hand side of diagram (\ref{toda2}), so $a\in \pi_*M$ and $b,c\in\pi_*R$.

The stable homotopy groups of a commutative ring spectrum $Q$ carry
an additional operation, the \emph{cup-one square}, defined as follows. 
Let $LQ$ be a fibrant replacement of $Q$ in the category of
all ring spectra. The ring spectrum $LQ$ is no longer commutative, but it remains commutative up to a coherent 
track $\alpha_1$ satisfying the idempotence and the hexagon axioms for symmetric monoidal
categories, compare Lemma \ref{haxa}. Given $a\in\pi_{2n}Q$, $n\geq0$, we take a representative $\bar{a}\colon
S^{2n}\r LQ$ where the spectrum $S^m$ is the $m$-fold suspension of the sphere spectrum $S$, $S^m=\S^mS$. 
The symmetry isomorphism for the smash square of
an even-dimensional sphere $\tau_\wedge\colon S^{2n}\wedge S^{2n}\cong S^{2n}\wedge S^{2n}$ is homotopic to the
identity. 
We can choose a track $\hat{\tau}\colon \tau_{\wedge}\rr 1_{S^{2n}\wedge S^{2n}}$, there are two such choices. 
Consider the
following diagram where $\mu$ is the product in $LQ$.
\begin{equation}\label{tcup1}
\xymatrix{S^{2n}\wedge S^{2n}\ar[dd]_{\;}="b"^{\tau_\wedge}\ar@/_30pt/[dd]_1^{\;}="a"
\ar[r]^-{\bar{a}\wedge\bar{a}}&LQ\wedge
LQ\ar[dd]_{\tau_\wedge}^{\;}="c"\ar[rd]_{\;}="d"^\mu&\\
&&LQ\\
S^{2n}\wedge S^{2n}\ar[r]_-{\bar{a}\wedge\bar{a}}&LQ\wedge LQ\ar[ru]_\mu&
\ar@{<=}"a";"b"^{\hat{\tau}}\ar@{=>}"c";"d"_{\alpha_1}}
\end{equation}
The pasting of this diagram is a self-track of $\mu(\bar{a}\wedge\bar{a})$. The classical Barcus-Barratt-Rutter 
isomorphism allows us to indetify this self-track 
with a homotopy class
$$Sq_1(a)\colon S^{4n+1}=\S(S^{2n}\wedge S^{2n})\To Q$$
measuring the difference between the pasting of (\ref{tcup1}) and the identity self-track on $\mu(\bar{a}\wedge\bar{a})$.
This element $Sq_1(a)\in\pi_{4n+1}Q$ is the cup-one square of $a$. One can check that $Sq_1(a)$ does not depend
on the representative $\bar{a}$. However in general it does depend on the choice of $\hat{\tau}$. 

Assume now that $R$ is a $Q$-algebra. The main example is $Q=S$ since all ring spectra are algebras over the
sphere spectrum. Let $u\colon Q\r R$ be the unit, let $\tau_{\wedge_Q}$ be the symmetry isomorphism for the smash product of $Q$-modules $\wedge_Q$, 
and let $\bar{L}$ be a fibrant replacement
functor in the category of $Q$-algebras. The diagram
\begin{equation*}
\xymatrix@C=40pt{Q\wedge_Q R\ar[r]^{u\wedge_Q 1}\ar[dd]_{\tau_{\wedge_Q }}&R\wedge_Q R\ar[rd]^{\text{mult.}}&\\
&&R\\
R\wedge_Q Q\ar[r]_{1\wedge_Q u}&R\wedge_Q R\ar[ru]_{\text{mult.}}&}
\end{equation*}
commutes but
\begin{equation*}
\xymatrix@C=40pt{\bar{L}Q\wedge_Q \bar{L}R\ar[r]^{\bar{L}u\wedge_Q 1}_{\;}="b"\ar[dd]_{\tau_{\wedge_Q }}&\bar{L}R\wedge_Q \bar{L}R\ar[rd]^{\text{mult.}}&\\
&&\bar{L}R\\
\bar{L}R\wedge_Q \bar{L}Q\ar[r]_{1\wedge_Q \bar{L}u}^{\;}="a"&\bar{L}R\wedge_Q \bar{L}R\ar[ru]_{\text{mult.}}&\ar@{=>}"a";"b"_{\alpha_1}}
\end{equation*}
is only commutative up to a certain track $\alpha_1$, compare Remark \ref{hint}.

If $b\in\pi_nQ$ is in the kernel of the homomorphism $\pi_*Q\r\pi_*R$ induced by 
the unit map $u\colon Q\r R$ we can choose a representative $x\colon S^n\r \bar{L}Q$ and a track $y\colon
(\bar{L}u)x\rr 0$. Given $a\in\pi_mR$ represented by $\bar{a}\colon S^m\r \bar{L}R$ the pasting of the diagram
\begin{equation*}
\xymatrix@C=40pt{&&&\\
S^n\wedge S^m\ar[r]|{x\wedge \bar{a}}^<(.56){\;}="c"\ar@/^30pt/[rr]_<(.28){\;}="d"^0\ar[dd]_{\tau_\wedge }&\bar{L}Q\wedge_Q \bar{L}R\ar[r]|{\bar{L}u\wedge_Q 1}_{\;}="b"\ar[dd]_{\tau_{\wedge_Q }}&\bar{L}R\wedge_Q \bar{L}R\ar[rd]^{\text{mult.}}&\\
&&&\bar{L}R\\
S^m\wedge S^n\ar[r]|{\bar{a}\wedge
x}_<(.56){\;}="e"\ar@/_30pt/[rr]^<(.28){\;}="f"_0&\bar{L}R\wedge_Q \bar{L}Q\ar[r]|{1\wedge_Q \bar{L}u}^{\;}="a"&\bar{L}R\wedge_Q \bar{L}R\ar[ru]_{\text{mult.}}&
\ar@{=>}"a";"b"_{\alpha_1}\ar@{=>}"c";"d"_{y\wedge\bar{a}}\ar@{=>}"e";"f"^{\bar{a}\wedge y}}
\end{equation*}
is a selft track of $0\colon S^{n+m}=S^n\wedge S^m\r \bar{L}R$, or equivalently an element
\begin{equation}\label{teta}
\theta_{(x,y)}(a)\in \pi_{n+m+1}R.
\end{equation}
One can easily check that this element does not depend on the choice of $\bar{a}$. In Theorem \ref{ae2} we
identify $\theta_{(x,y)}(a)$ in a purely algebraic way for connective spectra. This allows to deduce that
$$\theta_{(x,y)}\colon\pi_*R\To\S^{-1}\pi_*R$$
is a degree $n$ derivation of the graded ring $\pi_*R$ with coefficients in the desuspended bimodule
$\S^{-1}\pi_*R$. Moreover, Proposition \ref{loh} implies the following theorem.

\begin{thm}\label{app}
Let $Q$ be a connective commutative ring spectrum, let $R$ be a connective $Q$-algebra, and
let $I_Q(R)$ be the kernel of the ring
homomorphism $\pi_*u\colon \pi_*Q\r \pi_*R$. Then there is a $\pi_*Q$-module homomorphism
$$\theta\colon I_Q(R)/I_Q(R)^2\To HH^1_{\pi_*Q}(\pi_*R,\S^{-1}\pi_*R),\;\;b+I_Q(R)^2\mapsto \set{\theta_{(x,y)}}.$$
This homomorphism is natural in $R$ and
in $Q$ in the obvious way. This implies that $\theta$ vanishes when $R$ is commutative.
\end{thm}

All ring spectra are $S$-algebras so the theorem applies to all connective ring spectra for $Q=S$ the sphere
spectrum. The example of $\theta(p\cdot 1)$ of Toda shows that $\theta$ is non-trivial, see \cite{deritoda}.

For a brief remainder of $1$-dimensional Hochschild cohomology see the paragraph preceding Definition
\ref{deride1}.

In Section \ref{sss} we review more technical details about spectra which are needed for the proofs of
the main results in this paper.

\section{Pair algebras associated to ring spectra}\label{2}

Let $\C{pm}$ be the category of chain complexes of abelian groups concentrated in degree $0$ and $1$. Such chain
complexes are termed \emph{pair modules}, see \cite{asco}, and are denoted by
$$M=\left(\partial\colon M_1\r M_0\right).$$
The \emph{tensor product} of two pair modules
$M\ol{\otimes}N$ is obtained from the tensor product of chain complexes $M\otimes N$ by quotienting out the subcomplex
generated by elements in
dimension $2$, compare \cite{cmoo}. Hence the category $\C{pm}$ is symmetric monoidal. The unit object is $\Z=(0\r\Z)$. 
\emph{Quasi-isomorphisms} in $\C{pm}$ are
morphisms inducing isomorphisms on the ``homology'' functors,
\begin{equation}\label{hi}
h_0,h_1\colon\C{pm}\To\C{Ab},
\end{equation}
defined by
\begin{eqnarray*}
h_0 M&=&\coker\left(\partial\colon M_1\r M_0\right),\\
h_1 M&=&\ker\left(\partial\colon M_1\r M_0\right).
\end{eqnarray*}
These homology functors come equipped with a natural isomorphism
\begin{eqnarray}\label{lsm1}
h_0(M)\otimes h_0(N)&\cong&h_0(M\ol{\otimes}N),
\end{eqnarray}
and with natural pairings
\begin{equation}\label{lsm2}
\begin{array}{c}
h_0(M)\otimes h_1(N)\To h_1(M\ol{\otimes}N),\\
h_1(M)\otimes h_0(N)\To h_1(M\ol{\otimes}N).
\end{array}
\end{equation}
The isomorphism (\ref{lsm1}) and the pairings (\ref{lsm2}) satisfy the obvious 
associativity and commutativity properties, i.e. they define a lax symmetric
monoidal structure on the functor $(h_0,h_1,0,\dots)$ from $\C{pm}$ to the category of $\N$-graded abelian groups
where $\N=\set{0,1,2,\dots}$.

We mostly deal with $\N$-graded objects. An \emph{ungraded} object is an $\N$-graded object concentrated in
degree $0$.

The category $\C{pm}^\N$ of $\N$-graded pair modules $M=(M_{n,*},n\in\N)$ is also symmetric
monoidal with the tensor product
\begin{eqnarray*}
(M\ol{\otimes}N)_{n,*}&=&\bigoplus_{i+j=n}M_{i,*}\ol{\otimes}N_{j,*}.
\end{eqnarray*}
The symmetry isomorphism $\tau_{\overline{\otimes}}\colon M\ol{\otimes} N\cong N\ol{\otimes}M$ in
$\C{pm}^\N$ is defined by
$a_{p,q}\otimes b_{r,s}\mapsto (-1)^{pr}b_{r,s}\otimes a_{p,q}$. 
Here the first subscript of $a_{p,q}\in M_{p,q}$ will always denote the $\N$-grading and the second one the $\set{0,1}$-grading of the
corresponding pair module $M_{p,*}$. We point out that $a_{p,1}\otimes b_{r,1}$ represents $0$ in $M\ol{\otimes}N$.

A \emph{pair algebra} is a monoid in $\C{pm}^\N$. It is given by
an $\N$-graded ring $B_{*,0}$, an $\N$-graded $B_{*,0}$-bimodule $B_{*,1}$, and a
$B_{*,0}$-bimodule homomorphism $\partial\colon B_{*,1}\r B_{*,0}$ such that 
\begin{eqnarray}\label{cca}
a\cdot\partial(b)&=&\partial(a)\cdot b,\;\; a,b\in B_{*,1}.
\end{eqnarray}
If $B$ is a pair algebra then by using (\ref{lsm1}) and (\ref{lsm2}) $h_0B$ is an $\N$-graded ring and $h_1B$ is an
$h_0B$-bimodule in a natural way. 

Let $A$ be a ring. The most basic examples of pair algebras are the inclusion of a two-sided ideal $I\subset A$
and the zero morphism $0\colon M\r A$ where $M$ is an $A$-bimodule. An ungraded pair algebra is also termed a
``crossed bimodule'', see \cite[E.1.5.1]{ch}.

\begin{rem}\label{sh}
Given any ring $A$ and any $A$-bimodule $M$ consider the category consisting of pair
algebras $B$ together with a ring isomorphism $h_0B\cong A$ and a bimodule isomorphism $h_1B\cong M$. Morphisms
are pair algebra morphisms over $A$ and under
$M$. Then the set of connected components of this category is in natural bijection with 
$3$-dimensional Shukla cohomology 
$$SH^3(A,M),$$
see \cite{chsml}. Shukla cohomology is derived Hochschild cohomology. If the inclusions $h_1B\subset B_1$ and 
$\partial(B_1)\subset B_0$ split additively then the Shukla cohomology class associated to a pair algebra $B$ 
is in the image of the natural homomorphism from Hochschild cohomology
$$HH^3(A,M)\To SH^3(A,M),$$
see \cite{chsml} and \cite[Exercise E.1.5.1]{ch}.
\end{rem}

There is a notion of Massey product for pair algebras which is defined as follows.

\begin{defn}\label{mp1}
Let $B$ be a pair algebra. Given elements $a,b,c\in h_0 B$ of degree $p,q,r\in\N$ 
with $ab=0$ and $bc=0$ the
\emph{Massey product}  is the subset
$$\grupo{a,b,c}\subset h_1B_{p+q+r}$$
which is a coset of the subgroup $$(h_1B_{p+q})c+a(h_1B_{q+r})$$
defined as follows. Given $\bar{a}\in B_{p,0}$, $\bar{b}\in
B_{q,0}$, $\bar{c}\in B_{r,0}$ representing $a$, $b$, $c$, there exist $e\in B_{p+q,1}$, $f\in B_{q+r,1}$ such
that $\partial(e)=\bar{a}\cdot\bar{b}$,
$\partial(f)=\bar{b}\cdot\bar{c}$. Then one can easily check that 
$$-e\cdot\bar{c}+\bar{a}\cdot f\in h_1B_{p+q+r}\subset B_{p+q+r,1}.$$
The coset $\grupo{a,b,c}\subset h_1B_{p+q+r}$ 
coincides with the set of elements obtained in this way for all
different choices of $\bar{a}$, $\bar{b}$, $\bar{c}$, $e$ and $f$.
\end{defn}

We say that a ring spectrum $R$ \emph{neglects the Hopf map} if the stable Hopf map $\eta\in\pi_1S$ is in the
kernel of the ring homomorphism $\pi_*S\r\pi_*R$ induced on stable homotopy groups by the unit $S\r R$.

Recall that the \emph{desuspension} $\S^{-1}A$ of an $\N$-graded ring
$A$ is $(\S^{-1}A)_n=A_{n+1}$, $n\geq 0$, with the $\N$-graded $A$-bimodule structure defined by the formula
\begin{eqnarray*}
a\cdot(\S^{-1}b)\cdot c&=&(-1)^{\abs{a}}\S^{-1}(a\cdot b\cdot c).
\end{eqnarray*}
Here $\abs{a}$ is the degree of $a\in A$ and for any $x\in A_{n+1}$, $n\geq 0$, we denote by 
$\S^{-1}x\in (\S^{-1}A)_n$ the corresponding element in the desuspension. One of our main results is the following. 

\begin{thm}\label{a1}
There is a functor
$$\pi_{*,*}^\add\colon \left(\begin{array}{c}
\text{\emph{connective ring}}\\
\text{\emph{spectra neglecting} }\eta
\end{array}\right)\To\left(\text{\emph{pair algebras}}\right)$$
together with natural isomorphisms
\begin{eqnarray*}
h_0\pi_{*,*}^\add R&\cong&\pi_*R, \text{ of rings},\\
h_1\pi_{*,*}^\add R&\cong&\S^{-1}\pi_*R, \text{ of bimodules},
\end{eqnarray*}
such that the Massey products in $\pi_{*,*}^\add R$ coincide with the Toda brackets in $\pi_*R$.
\end{thm}

\begin{proof}
The functor $\pi_{*,*}^\add$ is defined by the formula $\pi_{*,*}^\add R=(\pi_{*,*}R)^\add$ where $\pi_{*,*}$ is the
functor in Theorem \ref{a2} and $(-)^\add$ is the additivization functor in (\ref{ref2}). Now the theorem
follows from Theorem \ref{a2} and Proposition \ref{til}.
\end{proof}

\begin{exm}\label{exmpa}
Common examples of connective ring spectra neglecting the Hopf map are the spectra associated to multiplicative
cohomology theories concentrated in even non-negative dimensions, such as connective complex $K$-theory $ku$, complex cobordism
$MU$, Brown-Peterson theory $BP$ \dots\quad If $R$ denotes any of these ring spectra then the pair algebra
$\pi_{*,*}^{add}R$ is quasi-isomorphic to 
$$0\colon \S^{-1}\pi_*R\To\pi_*R.$$
In fact if $B$ is any pair algebra such that $h_0B$ is concentrated in even degrees and $h_1B$ is concentrated in odd
degrees then we have a diagram of quasi-isomorphisms
$$\xymatrix{h_1B\ar[d]_0&
B_{\text{even},1}\oplus h_1B\ar@{->>}[l]_-{(0,1)}\ar@{^{(}->}[r]\ar[d]^{(\partial_{\text{even}},0)}&
B_{*,1}\ar[d]^\partial\\
h_0B&B_{\text{even},0}^{}\ar@{->>}[l]\ar@{^{(}->}[r]&B_{*,0}}$$

Algebras over the Eilenberg-Mac Lane spectrum $H\Z$ also neglect the Hopf map. Any connective $H\Z$-algebra is
weakly equivalent to the $H\Z$-algebra $\mathbb{H}D_*$ of a differential graded algebra $D_*$ concentrated in
non-negative dimensions, see \cite{hzas}. The theory developed in this paper allows one to compute a small model
of $\pi_{*,*}^\add\mathbb{H}D_*$. Indeed this pair algebra is quasi-isomorphic to
$$\bar{d}\colon\S^{-1}(D_*/d(D_*))\To Z_*.$$
Here $Z_*\subset D_*$ is the subring of cycles and $\bar{d}$ is induced by the differential in $D_*$, compare
\cite[3.6]{ceahh}.

The algebra of secondary mod $p$ cohomology operations $\mathcal{B}$, for $p$ a fixed prime, computed in
\cite{asco}, is a pair algebra. This pair algebra has
proved to be useful for computations of $d_2$ differentials in the classical Adams spectral sequence, see
\cite{ce3}. The pair algebra $\mathcal{B}$ corresponds to the coconnective version of $\pi_{*,*}R$ for
$R=\endo(H\Z/p)$ the endomorphism spectrum of the mod $p$ Eilenberg-Mac Lane spectrum. The
coconnective theory, however, is not considered in this paper.
\end{exm}

\begin{rem}\label{sh2}

Given an ungraded ring $A$ and an ungraded $A$-bimodule $M$, the Shukla cohomology $SH^3(A,M)$ is naturally included in Mac Lane cohomology, see \cite{slecfgfm,chsml},
which is isomorphic to the
topological Hochschild cohomology of the corresponding Eilenberg-Mac Lane ring spectrum, compare \cite{mlhthh},
\begin{equation}\label{momok}
SH^3(A,M)\hookrightarrow HML^3(A,M)\cong THH^3(HA,HM).
\end{equation}
If $R$ is a connective ring spectrum neglecting $\eta$ we can consider the 
ungraded pair algebra $\pi_{0,*}^\add R$ in the bottom degree of $\pi_{*,*}^\add R$ in Theorem \ref{a1} yielding
the associated Shukla cohomology class
\begin{equation*}
\grupo{\pi_{0,*}^\add R}\in SH^3(\pi_0R,\pi_1R).
\end{equation*}
On the other hand Lazarev \cite{htars} introduced the ``first $k$-invariant''
\begin{equation*}
k^1_R\in THH^3(H(\pi_0R), H(\pi_1R))
\end{equation*}
of the ring spectrum $R$.
We claim that the image under (\ref{momok}) of the class $\grupo{\pi_{0,*}^\add R}$ in topological Hochschild
cohomology can be identified with the first Postnikov invariant $k^1_R$ of the ring spectrum
$R$. 
\end{rem}

\section{$E_\infty$-pair algebras associated to commutative ring spectra}\label{3}

A \emph{commutative pair algebra} is a commutative monoid in $\C{pm}^\N$. 
We need a weaker notion given by a pair algebra which is 
commutative up to a coherent homotopy as follows. 

\begin{defn}\label{eipa}
An \emph{$E_\infty$-pair algebra} $C=(C,\smile_1)$ is a pair algebra $C=C_{*,*}$ together with a homomorphism
$$\smile_1\colon C_{*,0}\otimes C_{*,0}\To C_{*,1}$$
also called \emph{cup-one product}
such that given $x_i\in C_{n_i,0}$ and
$s\in C_{m,1}$,
\begin{eqnarray}
\nonumber (-1)^{n_1\cdot n_2}x_2\cdot x_1+\partial(x_1\smile_1x_2)&=&x_1\cdot x_2,\\
\nonumber (-1)^{m\cdot n_2}x_2\cdot s+\partial(s)\smile_1x_2&=&s\cdot x_2,\\
\nonumber (-1)^{n_1\cdot n_2}x_2\smile_1x_1+x_1\smile_1x_2&=&0,\\
\nonumber (x_1\cdot x_2)\smile_1x_3&=&(-1)^{n_2\cdot n_3}(x_1\smile_1 x_3)\cdot x_2+x_1\cdot (x_2\smile_1x_3).
\end{eqnarray}
The ring $h_0C$ is commutative and $h_1C$ is an $h_0C$-module in a natural way.
\end{defn}

\begin{rem}\label{hex}
If $\mu\colon C\ol{\otimes} C\r C$ denotes the multiplication on the $E_\infty$-pair algebra $C$ then the cup-one
product consists exactly of a chain homotopy
$\smile_1\colon\mu\tau_{\ol{\otimes}}\rr\mu$ which is idempotent
$$\xymatrix{\mu\ar@{=>}[r]^{\smile_1\tau_{\ol{\otimes}}}\ar@{=}[rd]&
\mu\tau_{\ol{\otimes}}\ar@{=>}[d]^{\smile_1}\\
&\mu}$$
and such that the hexagon
$$\xymatrix{
&\mu(1\otimes(\mu\tau_{\ol{\otimes}}))(\tau_{\ol{\otimes}}\otimes1)\ar@{=}[ld]\ar@{=>}[rd]^-{\;\;\;\;\;\;\;\;\;\;\mu(1\otimes\smile_1)(\tau_{\ol{\otimes}}\otimes1)}&\\
\mu\tau_{\ol{\otimes}}(1\otimes\mu)\ar@{=>}[d]_{\smile_1(1\otimes\mu)}&&
\mu(1\otimes\mu)(\tau_{\ol{\otimes}}\otimes1)\ar@{=}[d]\\
\mu(1\otimes\mu)\ar@{=}[rd]&&\mu((\mu\tau_{\ol{\otimes}})\otimes1)\ar@{=>}[ld]^-{\;\;\mu(\smile_1\otimes1)}\\
&\mu(\mu\otimes1)&}$$
commutes. These are the usual axioms in a symmetric monoidal category, see \cite{borceux2}. Indeed there is an
additional diagram, called the pentagon, which should commute. In this case it commutes
automatically since the product in $C$ is strictly associative.
\end{rem}

Pair algebras carry a notion of Massey product as we saw in Definition \ref{mp1}. $E_\infty$-pair algebras carry in
addition a notion of cup-one square.

\begin{defn}\label{cup1a}
Let $C$ be an $E_\infty$-pair algebra. 
Given an element $a\in h_0C_{2n,*}$ we define the \emph{cup-one square} $Sq_1(a)$ of $a$
in the following way. 
Choose a representative $\bar{a}\in C_{2n,0}$ of $a$.
Then we have
\begin{eqnarray*}
\partial(\bar{a}\smile_1\bar{a})&=&0,
\end{eqnarray*}
so we can set
$$Sq_1(a)=\bar{a}\smile_1\bar{a}\in h_1C_{4n,*}.$$
One readily checks that the cup-one square construction $Sq_1(a)$ does not depend on the choice of a
representative $\bar{a}$ made for its definition.
\end{defn}

\begin{thm}\label{ca1}
There is a diagram of functors 
$$\xymatrix{{\left(\begin{array}{c}
\text{\emph{connective commutative}}\\
\text{\emph{ring spectra neglecting} }\eta
\end{array}\right)}
\ar[d]_{\text{inclusion}}\ar[r]^-{\pi^\adc_{*,*}}&
{\left(E_\infty\text{\emph{-pair algebras}}\right)}
\ar[d]^{\text{forget }\smile_1}\\
{\left(\begin{array}{c}
\text{\emph{connective ring}}\\
\text{\emph{spectra neglecting} }\eta
\end{array}\right)}
\ar[r]^-{\pi^\add_{*,*}}&
{\left(\text{\emph{pair algebras}}\right)}}$$
Here the lower horizontal arrow is the functor in Theorem \ref{a1}. 
This diagram commutes up to natural quasi-isomorphisms as in (\ref{zz}) below.
Moreover, given a connective commutative ring
spectrum $Q$ the algebraic cup-one squares in 
$\pi_{*,*}^\adc Q$ coincide with the topologically-defined cup-one squares in $\pi_*Q$.
\end{thm}

We prove this theorem in Section \ref{pt2}. It follows from Theorems \ref{ca1} and \ref{a1} that Massey products
in $\pi_{*,*}^\adc R$ coincide with Toda brackets in $R$.

\begin{exm}\label{exmepa}
Let $D_*$ be a differential graded algebra concentrated in non-negative dimensions endowed with a $\smile_1$
operation consisting of a degree $1$ homomorphism
\begin{equation}\label{ah}
\smile_1\colon D_*\otimes D_*\To D_*
\end{equation}
satisfying the Hirsch formulas
\begin{eqnarray*}
d(x_1\smile_1x_2)&=&x_1\cdot x_2-(-1)^{n_1\cdot n_2}x_2\cdot x_1-d(x_1)\smile_1 x_2-(-1)^{n_1}x_1\smile_1d(x_2),\\
(x_1\cdot x_2)\smile_1x_3&=&(-1)^{n_1}x_1\cdot(x_2\smile_1x_3)+(-1)^{n_2\cdot n_3}(x_1\smile_1x_3)\cdot x_2,
\end{eqnarray*}
compare \cite[page 267]{ugss}. Then the pair algebra
$$\bar{d}\colon\S^{-1}(D_*/d(D_*))\To Z_*$$
in Example \ref{exmpa} is an $E_\infty$-pair algebra with cup-one product induced by $\smile_1$ in (\ref{ah}).
\end{exm}

\begin{exm}\label{3stem}
Let $S_{(3)}$ be the $3$-local sphere spectrum, which is a commutative ring spectrum, and let $\Z_{(3)}$ be the
$3$-local integers. The computations in this example are based on classical results which can be found for
instance in \cite{toda} and \cite{ccshgs}. Up to dimension $13$ the abelian group $\pi_{n}S_{(3)}$ is given by
$$\begin{array}{|c||c|c|c|c|c|c|c|c|c|c|c|c|c|c|}
n&0&1&2&3&4&5&6&7&8&9&10&11&12&13\\
\hline \pi_{n}S_{(3)}&\Z_{(3)}&0&0&\Z/3&0&0&0&\Z/3&0&0&\Z/3&\Z/9&0&\Z/3\\
\hline \text{generators} &1&-&-&\alpha_1&-&-&-&\alpha_2&-&-&\beta_1&\alpha_3'&0&\alpha_1\cdot\beta_1
\end{array}$$
The stable homotopy ring $\pi_*S_{(3)}$ is a commutative $\Z_{(3)}$-algebra. 
The only product in degrees $\leq 13$  which is not indicated above is
\begin{eqnarray*}
\alpha_1\cdot\alpha_2&=&0.
\end{eqnarray*} 
Moreover, the non-vanishing secondary homotopy operations in this range
are the Toda brackets
\begin{eqnarray*}
\grupo{\alpha_1,3,\alpha_1}&=&\alpha_2,\\
\grupo{\alpha_1,\alpha_1,\alpha_1}&=&\beta_1,\\
\grupo{\alpha_1,3,\alpha_2}&=&3\cdot\alpha_3'.
\end{eqnarray*} 
We are going to show in this example 
that the knowledge of the secondary homotopy operations gives us some information about the $E_\infty$-pair algebra
$\pi_{*,*}^{adc}S_{(3)}$ in degrees $\leq 12$, but it does not determine the whole structure.
This proves that the algebraic object $\pi_{*,*}^{adc}S_{(3)}$ contains more information than just the collection 
of classical secondary homotopy operations.

The natural morphism $\pi_{*,*}^{adc}S_{(3)}\r\pi_{*,*}^{adc}S_{(3)}\otimes\Z_{(3)}$ is a
quasi-isomorphism since $\pi_*S_{(3)}$ is $3$-local. Let us choose $a_1,a_2,b_1,a'_3\in \pi_{*,0}S_{(3)}\otimes\Z_{(3)}$
representing $\alpha_1,\alpha_2,\beta_1$, $\alpha'_3$, respectively. Taking pull-backs, as in Section
\ref{13},
we can define a quasi-isomorphic $E_\infty$-pair algebra $C\r\pi_{*,*}^{adc}S_{(3)}\otimes\Z_{(3)}$ such that
$C_{*,0}$ coincides in
dimensions $\leq 12$ with the free $\Z_{(3)}$-algebra generated by $a_1,a_2,b_1,a_3'$. 
Any submodule of a free $\Z_{(3)}$-module is free since $\Z_{(3)}$ is a principal ideal
domain, therefore, since $h_0C\cong\pi_*S_{(3)}$ and $h_1C\cong\S^{-1}\pi_*S_{(3)}$, 
the pair modules $C_{n,*}$ are given for $n\leq 12$ by
$$\begin{array}{|c||c|c||c|c|}
n&C_{n,0}&\text{generators}&C_{n,1}&\text{generators}\\
\hline 0&\Z_{(3)}&1&0&-\\
\hline 1&0&-&0&-\\
\hline 2&0&-&\Z/3&\S^{-1}\alpha_1\\
\hline 3&\Z_{(3)}&a_1&\Z_{(3)}&\bar{a}_1\\
\hline 4&0&-&0&-\\
\hline 5&0&-&0&-\\
\hline 6&\Z_{(3)}&a_1^2&\Z_{(3)}\oplus\Z/3&\bar{a}_{1}^{(2)},\S^{-1}\alpha_2\\
\hline 7&\Z_{(3)}&a_2&\Z_{(3)}&\bar{a}_{2}\\
\hline 8&0&-&0&-\\
\hline 9&\Z_{(3)}&a_1^3&\Z_{(3)}\oplus\Z/3&\bar{a}_{1}^{(3)},\S^{-1}\beta_1\\
\hline 10&\Z_{(3)}\oplus\Z_{(3)}\oplus\Z_{(3)}&a_1\cdot a_2, a_2\cdot
a_1,b_1&\Z_{(3)}\oplus\Z_{(3)}\oplus\Z_{(3)}\oplus\Z/9&\bar{a}_{1,2},\bar{a}_{2,1},\bar{b}_1,\S^{-1}\alpha'_3\\
\hline 11&\Z_{(3)}&a_3'&\Z_{(3)}&\bar{a}_{3}'\\
\hline 12&\Z_{(3)}&a_1^4&\Z_{(3)}\oplus\Z/3&\bar{a}_{1}^{(4)},\S^{-1}(\alpha_1\cdot\beta_1)
\end{array}$$
where $\bar{a}_1,\bar{a}_{1}^{(2)},\bar{a}_2,\bar{a}_{1}^{(3)},\bar{a}_{1,2},\bar{a}_{2,1},\bar{b}_1,
\bar{a}_{3}',\bar{a}_{1}^{(4)}$ are arbitrarily chosen elements satisfying 
\begin{eqnarray*}
\partial(\bar{a}_i)&=&3\cdot a_i,\;\; i=1,2,\\
\partial(\bar{a}_{1}^{(i)})&=&a_1^i,\;\;i=2,3, 4,\\
\partial(\bar{a}_{i,j})&=&a_i\cdot a_j,\;\;\set{i,j}=\set{1,2},\\
\partial(\bar{b}_1)&=&3\cdot b_1,\\
\partial(\bar{a}_3')&=&9\cdot a_3'.
\end{eqnarray*}
The desuspended elements are in $h_1C=\ker\partial$. Since $\partial$ is a $C_{*,0}$-bimodule homomorphism
\begin{eqnarray*}
a_1\cdot \bar{a}_1&=&3\cdot \bar{a}_{1}^{(2)}+m\cdot (\S^{-1}\alpha_2),\\
a_1\cdot \bar{a}_2&=&3\cdot \bar{a}_{1,2}+m'\cdot (\S^{-1}\alpha_3'),\\
\bar{a}_1\cdot a_2&=&3\cdot \bar{a}_{1,2}+m''\cdot (\S^{-1}\alpha_3'),
\end{eqnarray*}
for some $m\in\Z/3$ and $m',m''\in\Z/9$, 
and we can choose
$\bar{a}_{1}^{(3)}$ and $\bar{a}_{1}^{(4)}$ so that
\begin{eqnarray*}
a_1^i\cdot \bar{a}_{1}^{(2)}&=&\bar{a}^{(i+2)}_1,\;\;i=1,2.
\end{eqnarray*}
Moreover,
\begin{eqnarray*}
\bar{a}_1^{(2)}\cdot a_1^2&=&\bar{a}_1^{(2)}\cdot\partial(\bar{a}_1^{(2)})\\
&=&\partial(\bar{a}_1^{(2)})\cdot\bar{a}_1^{(2)}\\
&=&a_1^2\cdot\bar{a}_1^{(2)}\\
&=&\bar{a}_1^{(4)}.
\end{eqnarray*}

By the laws of an $E_\infty$-pair algebra we see that we can choose $a_1^{(2)},\bar{a}_{1,2},\bar{a}_{2,1}$ in such a way that
\begin{eqnarray*}
a_1\smile_1a_1&=&2\cdot \bar{a}_{1}^{(2)},\\
a_1\smile_1a_2&=&\bar{a}_{2,1}+\bar{a}_{1,2},
\end{eqnarray*}
and moreover
\begin{eqnarray*}
\bar{a}_1\cdot a_1&=&-a_1\cdot \bar{a}_1+a_1\smile_1(3\cdot a_1)\\
&=&-3\cdot \bar{a}_{1}^{(2)}-m\cdot (\S^{-1}\alpha_2)+6\cdot \bar{a}_{1}^{(2)}\\
&=&3\cdot \bar{a}_{1}^{(2)}-m\cdot (\S^{-1}\alpha_2),  \\
\bar{a}_2\cdot a_1&=&-a_1\cdot \bar{a}_2+a_1\smile_1(3\cdot a_2)\\
&=&-3\cdot \bar{a}_{1,2}-m'\cdot (\S^{-1}\alpha_3')+3\cdot\bar{a}_{2,1}+3\cdot\bar{a}_{1,2}\\
&=&3\cdot \bar{a}_{2,1}-m'\cdot (\S^{-1}\alpha_3'),  \\
a_2\cdot \bar{a}_1&=&-\bar{a}_1\cdot a_2+a_2\smile_1(3\cdot a_1)\\
&=&-3\cdot \bar{a}_{1,2}-m''\cdot (\S^{-1}\alpha_3')+3\cdot\bar{a}_{2,1}+3\cdot\bar{a}_{1,2}\\
&=&3\cdot \bar{a}_{2,1}-m''\cdot (\S^{-1}\alpha_3').  
\end{eqnarray*}

The action of $a_1$ on $\S^{-1}\alpha_1,\S^{-1}\alpha_2,$ yields $0$ from the left and from the right, as well
as the action of $a_2$ on $\S^{-1}\alpha_1$, and 
\begin{eqnarray*}
a_1\cdot(\S^{-1}\beta_1)&=&-\S^{-1}(\alpha_1\cdot\beta_1)\\
&=&-\S^{-1}(\beta_1\cdot\alpha_1)\\
&=&-(\S^{-1}\beta_1)\cdot a_1.
\end{eqnarray*}

Finally the known Toda brackets yield the equalities
\begin{eqnarray*}
\S^{-1}\alpha_2&=&-\bar{a}_1\cdot a_1+a_1\cdot \bar{a}_1\\
&=&-m\cdot(\S^{-1}\alpha_2),\\
\S^{-1}\beta_1&=&
-\bar{a}_{1}^{(2)}\cdot a_1+ a_1\cdot \bar{a}_{1}^{(2)}\\
&=&-\bar{a}_{1}^{(2)}\cdot a_1+\bar{a}_{1}^{(3)},\\
3\cdot(\S^{-1}\alpha_3')&=&-\bar{a}_1\cdot a_2+a_1\cdot\bar{a}_2\\
&=&(m'-m'')\cdot (\S^{-1}\alpha_3').
\end{eqnarray*}
In particular 
\begin{eqnarray*}
m&=&-1\in\Z/3,\\
m''&=&m'-3\in\Z/9. 
\end{eqnarray*}

The equations above determine, up to dimension $9$, the structure of $C$,
and hence the structure of $\pi_{*,*}^{adc}S_{(3)}$ up to quasi-isomorphism. In dimensions $\leq 12$ the structure is
determined up to the unknown constant $m'\in\Z/9$, which can not be deduced from the secondary homotopy
operations in $\pi_*S_{(3)}$. The reader who wants to familiarize himself with the new algebraic invariants introduced in this
paper can try to compute more about $\pi_{*,*}^{adc}S_{(3)}$ with the help of the classical calculations
available in the literature. One can also try with primes $p>3$.

We remark that the elements $\S^{-1}\alpha_1,
a_1,\bar{a}_1,a_2,\bar{a}_2,b_1,\bar{a}_{1,2},\bar{b}_1,\S^{-1}\alpha_3',a_3',\bar{a}_3'$, generate $C$ as an 
$E_\infty$-pair algebra in dimensions $\leq 12$ since the rest of elements $\bar{a}_{1}^{(i)}$,
$\S^{-1}\alpha_2$, $\bar{a}_{2,1},$ $\S^{-1}\beta_1,
\S^{-1}(\alpha_1\cdot\beta_1)$, $i=2,3,4$, can be written in terms of these generators, the product and the
cup-one product in $C$.

\end{exm}

\begin{rem}\label{BE}
The coherence conditions in Remark \ref{hex} are related to the Barratt-Eccles operad introduced in \cite{be1}.
Indeed, one can check that $E_\infty$-pair algebras are algebras over the following operad in $\C{pm}$. Let $\mathcal{BE}$ be the
simplicial Barratt-Eccles operad and let $C_*\mathcal{BE}$ be the operad of chain
complexes obtained by taking normalized cochains on $\mathcal{BE}$. 
If we truncate this operad by dividing out the subcomplexes generated by the elements in dimensions $\geq 2$ 
we obtain an operad $t_{\leq
2}C_*\mathcal{BE}$ in $\C{pm}$ whose associated algebras are exactly the $E_\infty$-pair algebras. 

Now let $A$ be an ungraded commutative ring and let $M$ be an ungraded $A$-module. We consider the
category consisting of ungraded $E_\infty$-pair
algebras $C_*$ with specified isomorphisms $h_0C_*\cong A$ and $h_1C_*\cong M$ and $E_\infty$-pair algebra morphisms over $A$ and under
$M$. Then the set of components of this category is in natural bijection with the operadic cohomology group
$$H^3_{t_{\leq 2}C_*\mathcal{BE}}(A,M)$$
defined in \cite{cmoo}, which in dimension $3$ is the ``commutative analogue'' of Shukla cohomology, see Remark \ref{sh}. 
This cohomology is likely to be related to the $\Gamma$-cohomology of Robinson and Whitehouse,
compare \cite[7]{cmoo}. Moreover, this operadic cohomology group should map to the third topological
Andre-Quillen cohomology group, where the first Postnikov invariant of a
connective commutative ring spectrum $R$ lives, see \cite{aqccsa}. In fact, if $R$ neglects the Hopf
map the class of the ungraded $E_\infty$-pair algebra $\pi_{0,*}^\adc R$ should be mapped this way to the first Postnikov
invariant of $R$ as in the non-commutative case, see Remark \ref{sh2}.
\end{rem}

\section{Algebras over commutative ring spectra neglecting the Hopf map}\label{31}

In stable homotopy theory we have the notion of a connective algebra over a connective 
commutative ring spectrum $Q$, which
is a connective ring spectrum $R$ with a two-sided action of $Q$. If $Q$ neglects the Hopf map then so $R$
does and by Theorems \ref{a1} and \ref{ca1} the ring spectra $R$ and $Q$ give rise to a
pair algebra $\pi_{*,*}^{add}R$ and to an $E_\infty$-pair
algebra $\pi_{*,*}^{adc} Q$. In this section we study the algebraic properties of the induced
action of $\pi_{*,*}^{adc} Q$ on  $\pi_{*,*}^{add}R$, obtaining Theorem \ref{app} for spectra neglecting the
Hopf map as an application. 

\begin{defn}\label{meqpa1}
Let $C$ be an $E_\infty$-pair algebra. A \emph{$C$-algebra} is a
pair algebra $B$ 
together with a pair algebra morphism $u\colon C\r B$, called the \emph{unit}, and 
two external cup-one product operations,
\begin{equation*}
\begin{array}{c}
\smile_1\colon B_{*,0}\otimes C_{*,0}\To C_{*,1},\\
\smile_1\colon C_{*,0}\otimes B_{*,0}\To C_{*,1},
\end{array}
\end{equation*}
such that 
given $x_i\in B_{n_i,0}$, $s\in B_{m,1}$,
$\tilde{x}_i\in C_{n_i,0}$, $\tilde{s}\in C_{m,1}$, the following equations hold, compare Definition \ref{eipa}.
\begin{eqnarray}
\label{compaa} u_1(\tilde{x}_1\smile_1 \tilde{x}_2)&=&u_0(\tilde{x}_1)\smile_1\tilde{x}_2\\
\nonumber &=&\tilde{x}_1\smile_1u_0(\tilde{x}_2),\\
\label{loes} (-1)^{n_1\cdot n_2}x_2\cdot u_0(\tilde{x}_1)+\partial(\tilde{x}_1\smile_1x_2)&=&u_0(\tilde{x}_1)\cdot x_2,\\
\label{biende} (-1)^{m\cdot n_2}u_0(\tilde{x}_2)\cdot s+\partial(s)\smile_1\tilde{x}_2&=&s\cdot u_0(\tilde{x}_2),\\
\label{node} (-1)^{m\cdot n_2}x_2\cdot u_1(\tilde{s})+\partial(\tilde{s})\smile_1x_2&=&u_1(\tilde{s})\cdot x_2,\\
\label{comi} (-1)^{n_1\cdot n_2}x_2\smile_1\tilde{x}_1+\tilde{x}_1\smile_1x_2&=&0,\\
\label{otrolin} \qquad (-1)^{n_2\cdot n_3}(\tilde{x}_1\smile_1 x_3)\cdot
u_0(\tilde{x}_2)+u_0(\tilde{x}_1)\cdot (\tilde{x}_2\smile_1x_3)&=&(\tilde{x}_1\cdot \tilde{x}_2)\smile_1x_3,\\
\label{isder} (-1)^{n_2\cdot n_3}(x_1\smile_1 \tilde{x}_3)\cdot
x_2+x_1\cdot (x_2\smile_1\tilde{x}_3)&=&(x_1\cdot x_2)\smile_1\tilde{x}_3.
\end{eqnarray}

The ring $h_0B$ is an $h_0C$-algebra in such a way that the $h_0B$-bimodule $h_1B$ restricts to an $h_0C$-module.
\end{defn}

\begin{rem}
Any $E_\infty$-pair algebra is an algebra over itself. Moreover, given a morphism of $E_\infty$-pair algebras $f\colon \bar{C}\r C$ and a $C$-algebra $B$ then
$B$ has a $\bar{C}$-algebra structure, denoted by $f^*B$, with unit $uf\colon\bar{C}\r B$ and 
external cup-one products obtained by precomposing with $f$.
\end{rem}

The main new secondary homotopy operation on algebras over $E_\infty$-pair algebras is introduced in the
next definition. We need to recall first some basics facts on the $1$-dimensional Hochschild cohomology of
graded algebras.

Recall that given a commutative algebra $K$, a $K$-algebra $A$,
an $A$-bimodule $M$ which is also a $K$-module, and $n\in\Z$, a \emph{degree $n$ derivation} is a
degree $n$ homomorphism of $K$-modules $d\colon A\r M$
such that $d(a\cdot b)=d(a)\cdot b+(-1)^{n\abs{a}}a\cdot d(b)$. Any $m\in M$ gives rise to a degree $\abs{m}$
derivation $d_m(a)=m\cdot a-(-1)^{\abs{m}\abs{a}}a\cdot m$.
These derivations are called \emph{inner derivations} and they form a submodule $\inn_K(A,M)\subset \der_K(A,M)$ of
the $K$-module of all derivations.
The \emph{1-dimensional Hochschild cohomology} group of $A$ with coefficients in $M$ coincides with the quotient
$K$-module
\begin{eqnarray*}
HH^1_K(A,M)&=&\der_K(A,M)/ \inn_K(A,M).
\end{eqnarray*}

\begin{defn}\label{deride1}
Let $C$ be an $E_\infty$-pair algebra and let $B$ be a $C$-algebra. We associate to any pair $(x,y)$
with $x\in C_{0,n}$ and $y\in B_{n,1}$ such that $\partial(y)=u_0(x)$
the degree $n$ derivation $$\theta_{(x,y)}\colon h_0B\To h_1B$$
defined by
\begin{eqnarray*}
\theta_{(x,y)}(a)&=&-y\cdot\bar{a}+(-1)^{n\abs{a}}\bar{a}\cdot y+x\smile_1\bar{a}.
\end{eqnarray*}
Here $\bar{a}\in B_{*,0}$ is a representative of $a\in h_0B$. This is indeed an element in $h_1B$ by (\ref{loes}). Moreover, 
$\theta_{(x,y)}=0$ provided $B$ is an $E_\infty$-pair algebra and $u$ is a morphism of
$E_\infty$-pair algebras. It follows from equations (\ref{biende}) and (\ref{comi}) that $\theta_{(x,y)}(a)$ does not depend
on the choice of $\bar{a}$. In order to check the derivation property for $\theta_{(x,y)}$ one uses equations
(\ref{isder}) and (\ref{comi}). It is $h_0C$-linear since $\theta_{(x,y)}$ vanishes in the image of $h_0u$. This follows from
equations (\ref{compaa}) and (\ref{node}).
\end{defn}

\begin{prop}\label{loh1}
Given an $E_\infty$-pair algebra $C$ and a $C$-algebra $B$ let $I_C(B)$ be the kernel of the ring
homomorphism $h_0u\colon h_0C\r h_0B$. There is an $h_0C$-module homomorphism
$$\theta\colon I_C(B)/I_C(B)^2\To HH^1_{h_0C}(h_0B,h_1B)$$
which sends an element in $I_C(B)/I_C(B)^2$
represented by $x\in C_{n,0}$ to the element in Hochschild cohomology represented by a derivation
$\theta_{(x,y)}$ as in Definition \ref{deride} for any $y\in B_{n,1}$ with $\partial(y)=u_0(x)$. The homomorphism
$\theta$ is natural in $B$ and in $C$ in the obvious way. In particular $\theta=0$ when $B$ is an $E_\infty$-pair algebra and $u$ is a morphism of $E_\infty$-pair algebras.
\end{prop}

\begin{proof}
Using equation (\ref{node}) one can check that the cohomology class of the derivation $\theta_{(x,y)}$ only depends on the
projection of $x$ to $I_C(B)$. In this way one obtains a homomorphism $\theta$ from $I_C(B)$. 
Given $z\in C_{m,0}$ representing $c\in h_0C$
\begin{eqnarray*}
\theta_{(z\cdot x,u_0(z)\cdot y)}(a)&=&-u_0(z)\cdot y\cdot \bar{a}+(-1)^{(m+n)\abs{a}}\bar{a}\cdot u_0(z)\cdot y+(z\cdot x)\smile_1\bar{a}\\
\mbox{\scriptsize(\ref{loes})}
\qquad&=&-u_0(z)\cdot y\cdot \bar{a}+(-1)^{n\abs{a}}u_0(z)\cdot \bar{a}\cdot y\\
&&-(-1)^{n\abs{a}}\partial(z\smile_1\bar{a})\cdot y+(z\cdot x)\smile_1\bar{a}\\
&=&-u_0(z)\cdot y\cdot \bar{a}+(-1)^{n\abs{a}}u_0(z)\cdot \bar{a}\cdot y\\
&&-(-1)^{n\abs{a}}(z\smile_1\bar{a})\cdot u_0(x)+(z\cdot x)\smile_1\bar{a}\\
\mbox{\scriptsize(\ref{otrolin})}\qquad&=&-u_0(z)\cdot y\cdot \bar{a}+(-1)^{n\abs{a}}u_0(z)\cdot \bar{a}\cdot y+u_0(z)\cdot (x\smile_1\bar{a})\\
&=&u_0(z)\cdot(-y\cdot \bar{a}+(-1)^{n\abs{a}} \bar{a}\cdot y+x\smile_1\bar{a})\\
&=&(h_0u)(c)\cdot\theta_{(x,y)}(a).
\end{eqnarray*}
This formula shows that $\theta$ is $h_0C$-linear and also that $\theta(I_C(B)^2)=0$.
\end{proof}

Algebras over commutative ring spectra neglecting the Hopf map are a source of algebras over $E_\infty$-pair
algebras as the following theorem shows.

\begin{thm}\label{ae1}
Let $Q$ be a connective commutative ring spectrum neglecting the Hopf map 
and let $\pi_{*,*}^{adc}Q$ be the $E_\infty$-pair algebra
given by Theorem \ref{ca1}. There is a diagram of functors
$$\xymatrix{{\left({\emph{connective $Q$-algebras}}\right)}
\ar[d]_{\text{forget}}\ar[r]^-{\pi^{ada}_{*,*}}&
{\left(\pi_{*,*}^{adc}Q\text{\emph{-algebras}}\right)}
\ar[d]^{\text{forget}}\\
{\left(\begin{array}{c}
\text{\emph{connective ring}}\\
\text{\emph{spectra neglecting} }\eta
\end{array}\right)}
\ar[r]^-{\pi_{*,*}^{add}}&
{\left(\text{\emph{pair algebras}}\right)}}$$
which commutes up to natural quasi-isomorphisms.
Here the lower arrow is the functor in Theorem \ref{a1}.
Moreover, if $f\colon Q'\r Q$ is a morphism of connective commutative ring spectra neglecting the Hopf map then the square
$$\xymatrix{{\left({\emph{connective $Q$-algebras}}\right)}
\ar[d]_{f^*}\ar[r]^-{\pi_{*,*}^{ada}}&
{\left(\pi_{*,*}^{adc}Q\text{\emph{-algebras}}\right)}
\ar[d]^{(\pi_{*,*}^{adc}f)^*}\\
{\left({\emph{connective $Q'$-algebras}}\right)}
\ar[r]^-{\pi_{*,*}^{ada}}&
{\left(\pi_{*,*}^{adc}Q'\text{\emph{-algebras}}\right)}}$$
commutes up to natural quasi-isomorphisms. Furthermore, if $R$ is a connective $Q$-algebra, $a\in\pi_mR$, and $b\in\pi_nQ$ is in
$I_Q(R)=\ker[\pi_*Q\r\pi_*R]\cong I_{\pi_{*,*}^{adc}Q}(\pi_{*,*}^{ada}R)$, then  the element
$\theta_{(x,y)}(a)\in\pi_{n+m+1}R\cong(h_1\pi_{*,*}^{ada}R)_{n+m}$ defined in (\ref{teta}) can be identified with
$\theta_{(x,y)}(a)$ in the sense of Definition \ref{deride1}.
\end{thm}

This theorem can be derived from Theorem \ref{ae2} in the same way as we derive Theorem \ref{ca1} from Theorem \ref{ca2} in Section \ref{16}, see Remark \ref{kiki} for details.

\section{Square groups and quadratic pair modules}\label{sqpm}\label{4}

We now describe quadratic generalizations of pair modules, pair algebras and $E_\infty$-pair algebras. The
quadratic concepts are needed to achieve results as above on general ring spectra not neglecting the Hopf map. In
fact, the Hopf map requires the use of quadratic structure as developed in \cite{ecg,qaI}.
In this section we quickly recall the basic quadratic algebra which is needed in this paper. We do not recall, for
example, the definition of the quadratic tensor product, but we later give explicit definitions of monoids and
modules in the monoidal category of quadratic pair modules,
see Definitions \ref{alg1} and \ref{mod}. We relate the quadratic concepts to the additive (non-quadratic)
situation which was used in the previous sections.

A \emph{square group} $X$ is a diagram
$$X=(X_e\mathop{\leftrightarrows}\limits^P_HX_{ee})$$
where $X_e$ is a group with an additively written group law, $X_{ee}$
is an abelian group, $P$ is a homomorphism, $H$ is a \emph{quadratic map}, i.e. a function such
that the cross effect
$$(a|b)_H=H(a+b)-H(b)-H(a),\;\;a,b\in X_e,$$
is bilinear, and the following relations are
satisfied, $x,y\in X_{ee}$,
\begin{eqnarray*}
(Px|b)_H&=&(a|Py)\;=\;0,\\
P(a|b)_H&=&-a-b+a+b,\\
PHP(x)&=&P(x)+P(x).
\end{eqnarray*}
In particular $X_e$ is a group of nilpotency class $2$ and $P$ maps $X_{ee}$ to the center of $X_e$, so its image is a normal subgroup.
Let $\C{SG}$ 
be the category of square groups where a morphism of square groups $f\colon X\r Y$ is given by
homomorphisms $f_e\colon X_e\r Y_e$,
$f_{ee}\colon X_{ee}\r Y_{ee},$
commuting with $P$ and $H$. An abelian group is identified with a square group $X$ with $X_{ee}=0$. In this way
we have a full inclusion
$\C{Ab}\subset\C{SG}$ of the category of abelian groups into a category of square groups.

For any square group $X$ the function $$T=HP-1\colon X_{ee}\To X_{ee}$$ is an involution, i.e. a homomorphism with $T^2=1$. Moreover,
$$\Delta\colon X_e\To X_{ee}\colon x\mapsto (x|x)_H-H(x)+TH(x)$$
is a homomorphism which satisfies $T\Delta=-\Delta$. The cross effect induces a homomorphism
$$(-|-)_H\colon\otimes^2\coker P\To X_{ee},$$
where $\otimes^2A=A\otimes A$ is the tensor square of an abelian group. We say that $X$ is \emph{good} if this homomorphism is an isomorphism.

As an example of square group we can consider
$$\Z_\nill=(\Z\mathop{\leftrightarrows}\limits^P_H\Z)$$
with $P=0$ and $H(n)=\binom{n}{2}=\frac{n(n-1)}{2}$. 

For a pointed set $E$ with base-point $*\in E$ let 
$\Z[E]$ be the free abelian group with basis $E-\set{*}$, and let $\grupo{E}_\nill$ be the
quotient of the free group $\grupo{E}$ with basis $E-\set{*}$ by triple commutators. Both $\Z[-]$ and
$\grupo{-}_\nill$ are functors on the category of pointed sets and we have the natural abelianization
homomorphism $\grupo{E}_\nill\twoheadrightarrow\Z[E]$ which is supressed from notation. 

There is a square group 
\begin{equation}\label{znil}
\Z_\nill[E]=(\grupo{E}_\nill \mathop{\leftrightarrows}\limits^P_H \otimes^2\Z[E]),
\end{equation}
where $P(a\otimes b)=-b-a+b+a$, $H(e)=0$ for any $e\in E$ and
$(s|t)_H=t\otimes s$ so that $\Z_\nill[S^0]=\Z_\nill$. 
These are canonical examples of good square groups. 

A \emph{quadratic pair module} $C$ is a morphism $\partial\colon
C_{(1)}\r C_{(0)}$ between square groups
\begin{eqnarray*}
C_{(0)}&=&(C_{0}\mathop{\leftrightarrows}^{P_0}_HC_{ee}),\\
C_{(1)}&=&(C_{1}\mathop{\leftrightarrows}^P_{H_1}C_{ee}),
\end{eqnarray*}
such that $\partial_{ee}=1\colon C_{ee}\r C_{ee}$ is the identity
homomorphism. Hence $C$ is completely determined by
the diagram
\begin{equation}\label{ya}
\xymatrix{&C_{ee}\ar[ld]_P&\\C_{1}\ar[rr]_\partial&&C_{0}\ar[lu]_H}
\end{equation}	
where $\partial=\partial_e$, $H_1=H\partial$ and $P_0=\partial P$. We say that $C$ is \emph{$0$-good} if
$C_{(0)}$ is a good square group. Moreover, $C$ is \emph{$0$-free} if $C_{(0)}$ is of the form (\ref{znil}).

Morphisms of quadratic pair modules $f\colon C\r D$ are 
given by group homomorphisms $f_{0}\colon C_{0}\r D_{0}$,
$f_{1}\colon C_{1}\r D_{1}$, $f_{ee}\colon C_{ee}\r D_{ee}$,
commuting with $H$, $P$ and $\partial$ in (\ref{ya}). They form a
category denoted by $\C{qpm}$. A quadratic pair module with $C_{ee}=0$ is the same as a pair module. 
In this way the category $\C{pm}$ of pair modules is regarded as a full subcategory
$\C{pm}\subset\C{qpm}$.  
We write $\C{qpm}_H$ for the full subcategory of $0$-good quadratic pair modules.

We have functors
$$h_0,h_1\colon\C{qpm}\To\C{Ab}$$ defined by
$h_0C=\coker\partial$ and  
$h_1C=\ker\partial$ as in (\ref{hi}).  A \emph{quasi-isomorphism} of quadratic pair modules is a morphism inducing isomorphisms on
$h_0$ and $h_1$. 

\begin{rem}
As shown in
\cite{qaI}, the category $\C{SG}$ is a symmetric monoidal category with the tensor product $X\odot Y$ of square
groups and with unit object $\Z_\nill$. Moreover, there is a symmetric monoidal structure $\odot$ on $\C{qpm}$
defined by using the
tensor product of square groups, see \cite{2hg3}. The unit object for the monoidal structure in $\C{qpm}$ is 
$$\ol{\Z}_\nill\;=\;\left(\begin{array}{c}\xymatrix{&\Z\ar@{->>}[ld]_P&\\\Z/2\ar[rr]_0&&\Z\ar[lu]_{H}}\end{array}\right)$$
where $P$ is the non-trivial homomorphism and $H(n)=\binom{n}{2}$ as in (\ref{znil}).
\end{rem}

The functors $h_0$ and $h_1$ also admit pairings as in (\ref{lsm1}) and (\ref{lsm2}), so they assemble to a lax symmetric monoidal functor from $\C{qpm}$ to the category of $\N$-graded
abelian groups.
This implies, since $h_1\overline{\Z}_\nill=\Z/2$, that the unit isomorphism $C\odot\overline{\Z}_\nill\cong C$
induces a natural homomorphism
\begin{equation}\label{-w}
k\colon h_0C\otimes\Z/2\To h_1C
\end{equation}
which, in fact, is given by the formula $k(x)=P(x|x)_H$. This homomorphism, termed the \emph{$k$-invariant} of $C$, corresponds to the action of the Hopf map in homotopy groups,
see Propositions \ref{qlo} and \ref{qlo2} below. 

The full inclusion $\C{Ab}\subset\C{SG}$ above
admits a reflection given by
\begin{equation}\label{ref1}
(-)^\add\colon\C{SG}\To\C{Ab},\;\; X^\add=\coker P.
\end{equation}

The full inclusion $\C{pm}\subset\C{qpm}$ also admits a reflection functor
\begin{equation}\label{ref2}
(-)^\add\colon \C{qpm}\To\C{pm}\colon C\mapsto C^\add
\end{equation}
which sends $\partial\colon C_{(1)}\r C_{(0)}$ to the abelian group homomorphism
$$\partial^\add\colon C_{(1)}^\add=\coker P\To\coker \partial P=C_{(0)}^\add$$
induced by $\partial\colon C_1\r C_0$.

The natural projection $C\r C^\add$, which is the unit of the reflection (\ref{ref2}), induces an isomorphism on
$h_1$ but not
necessarily on $h_0$. For $h_0$ we have the following result.

\begin{prop}\label{til}
Let $C$ be a $0$-free quadratic pair module. The natural projection $C\r C^\add$ is a quasi-isomorphism 
if and only if the $k$-invariant of $C$, as defined in (\ref{-w}), is trivial.
\end{prop}

\begin{proof}
The $k$-invariant is invariant under quasi-isomorphisms, and pair modules have obviously trivial $k$-invariant, therefore
the ``only if'' part follows.

Suppose now that $C$ is a $0$-free quadratic pair module with trivial $k$-invariant and $C_{(0)}=\Z_\nill[E]$. 
Then the homomorphism
$$P\colon \otimes^2\Z[E]\To C_1$$
factors through the exterior square by $$\tilde{P}\colon\wedge^2\Z[E]\To C_1.$$
The composite $\partial\tilde{P}\colon\wedge^2\Z[E]\r\grupo{E}_\nill$ is known to be an injection with cokernel $\Z[E]$, therefore $\tilde{P}$ is
also injective. This implies that we have a commutative diagram with exact rows
$$\xymatrix{\wedge^2\Z[E]\ar@{^{(}->}[r]^-{\tilde{P}}\ar@{=}[d]&C_1\ar[d]^\partial\ar@{->>}[r]&
C_1^\add\ar[d]^{\partial^\add}\\
\wedge^2\Z[E]\ar@{^{(}->}[r]&C_0\ar@{->>}[r]&
C_0^\add}$$
hence the proposition follows from the ``snake lemma''.
\end{proof}

\begin{rem}
The inclusion $\C{Ab}\subset\C{SG}$ 
strictly preserves the tensor product.  Moreover,  $\C{Ab}\subset\C{SG}$ is a lax symmetric monoidal functor with 
lax symmetric monoidal structure given on units by the square group morphism $\Z_\nill\r\Z$
which is the identity on the $e$-level. The reflection (\ref{ref1}) is strict
symmetric monoidal
since it preserves the tensor product and the unit
$\Z_\nill^\add=\Z$.

Similarly the inclusion $\C{pm}\subset\C{qpm}$ is lax symmetric monoidal. It is strictly compatible with the
tensor product, and the structure morphism $\overline{\Z}_\nill\r \Z$ between the units is the identity on the
$0$-level. As in the case of square groups the reflection functor (\ref{ref2}) is strict symmetric monoidal.
\end{rem}

\begin{rem}\label{per1}
A \emph{permutative category} is a small category $\C{P}$ with a symmetric monoidal structure $\oplus$ with unit object $0$
which is strictly associative
\begin{eqnarray*}
(A\oplus B)\oplus C&=&A\oplus(B\oplus C),
\end{eqnarray*}
and strictly unital 
\begin{eqnarray*}
A\oplus 0&=&A\\
&=&0\oplus A,
\end{eqnarray*}
but not necessarily strictly commutative, i.e. the isomorphism
\begin{eqnarray*}
\tau_\oplus\colon A\oplus B&\cong &B\oplus A.
\end{eqnarray*}
needs not be the indentity.

A quadratic pair module $C$ yields a permutative category $\C{P}C$ with object set $C_0$. The morphism set is
the semidirect product $C_0\ltimes C_1$, where $C_0$ acts on the right of $C_1$ by the formula
$c_1^{c_0}=c_1+P(\partial(c_1)|c_0)_H$. An element $(c_0,c_1)\in C_0\ltimes C_1$ is a morphism $(c_0,c_1)\colon
c_0\r c_0+\partial(c_1)$.The composition in $\C{P}C$ is determined by the sum in $C_0\ltimes C_1$. Moreover, the
monoidal structure $\oplus$ is also defined by the sum in $C_0$ and $C_0\ltimes C_1$, that is
$A\oplus B=A+B$. The symmetry isomorphism is
\begin{eqnarray*}
\tau_\oplus=(A+B,P(B|A)_H)\colon A\oplus B&\cong &B\oplus A
\end{eqnarray*}
Notice that $\C{P}C$ is indeed a groupoid.
\end{rem}

\section{Quadratic pair algebras associated to ring spectra}\label{5}

In this section we give an explicit description of quadratic pair algebras. We also define Massey
products for quadratic pair algebras and show that they correspond to Toda brackets. 

\begin{defn}\label{alg1}
A \emph{quadratic pair algebra} $B$ is an $\N$-graded quadratic pair
module $\{B_{n,*}, n\in \N\}$, together with multiplications, $n,m\in \N$,
$$\begin{array}{c}
B_{n,0}\times B_{m,0}\st{\cdot}\To B_{n+m,0},\\
B_{n,0}\times B_{m,1}\st{\cdot}\To B_{n+m,1},\\
B_{n,1}\times B_{m,0}\st{\cdot}\To B_{n+m,1},\\
B_{n,ee}\times B_{m,ee}\st{\cdot}\To B_{n+m,ee},
\end{array}$$
and an element $1\in B_{0,0}$ with $H(1)=0$ which is a (two-sided) unit for the
first three multiplications and such that $(1|1)_H\in B_{0,ee}$ is a (two-sided) unit for the fourth
multiplication. These multiplications are associative in all possible ways.
Moreover, the following lists
of equations are satisfied for $x,x_i\in B_{*,0}$, $s,s_i\in
B_{*,1}$ and $a_i\in B_{*,ee}$.
The multiplications $\cdot$ are
always right linear
\begin{eqnarray*}
x_1\cdot(x_2+x_3)&=&x_1\cdot x_2 + x_1\cdot x_3,\\
x\cdot(s_1+s_2)&=&x\cdot s_1 + x\cdot s_2,\\
s\cdot(x_1+x_2)&=&s\cdot x_1 + s\cdot x_2,\\
a_1\cdot(a_2+a_3)&=&a_1\cdot a_2 + a_1\cdot a_3.
\end{eqnarray*}
The multiplications $\cdot$ satisfy the following left
distributivity laws
\begin{eqnarray*}
(x_1+x_2)\cdot x_3&=&x_1\cdot x_3+x_2\cdot x_3 + \partial P((x_2|x_1)_H\cdot H(x_3)),\\
(x_1+x_2)\cdot s&=&x_1\cdot s+x_2\cdot s + P((x_2|x_1)_H\cdot H\partial(s)),\\
(s_1+s_2)\cdot x&=&s_1\cdot x+s_2\cdot x + P((\partial(s_2)|\partial(s_1))_H\cdot H(x)),\\
(a_1+a_2)\cdot a_3&=&a_1\cdot a_3+a_2\cdot a_3.
\end{eqnarray*}
The homomorphisms $\partial$ are compatible with the multiplications
$\cdot$ in the following sense
\begin{eqnarray*}
\partial(x\cdot s)&=&x\cdot\partial(s),\\
\partial(s\cdot x)&=&\partial(s)\cdot x,\\
\partial(s_1)\cdot s_2&=&s_1\cdot\partial(s_2).
\end{eqnarray*}
And finally, we have compatibility conditions for the
multiplications $\cdot$ and the maps $P$, $H$, $\Delta$, and
$(-|-)_H$,
\begin{eqnarray*}
P((x|x)_H\cdot a)&=&x\cdot P(a),\\
P(a\cdot\Delta(x))&=&P(a)\cdot x,\\
H(x_1\cdot x_2)&=& (x_1|x_1)_H\cdot H(x_2)+H(x_1)\cdot\Delta(x_2),\\
H\partial P(a_1\cdot a_2)&=&H\partial P(a_1)\cdot a_2+a_1\cdot H\partial P(a_2)-H\partial P(a_1)\cdot
H\partial P(a_2),\\
(x_1\cdot x_2|x_3\cdot x_4)_H&=&(x_1|x_3)_H\cdot (x_2|x_4)_H.
\end{eqnarray*}
\end{defn}

\begin{rem}\label{al}
In \cite{2hg3} we define the monoidal category $\C{qpm}^\N_H$ of $\N$-graded $0$-good quadratic pair modules.
A quadratic pair algebra $B$ which is $0$-good is the same as a monoid in $\C{qpm}_H^\N$. The product $a_i\cdot
b_j$ of two elements in levels $i,j\in\set{0,1}$ with $i+j\leq 1$ is given by the image of the right linear
generator $a_i\ul{\circledcirc}b_j$ in the tensor product of square groups under the structure
morphism $B\odot B\r B$ of the monoid $B$. If $i=j=ee$ then $a_i\cdot b_j$ is the image of $a_i\otimes b_j \in 
(B\odot B)_{ee}=B_{*,ee}\otimes B_{*,ee}$.
To see this requires a technical computation which only uses the formulas for the tensor products of square groups and
quadratic pair modules in \cite{3mlc,2hg3}. We leave it to the reader.
\end{rem}

If $B$ is a quadratic pair algebra then $h_0B$ is an $\N$-graded ring and $h_1B$ is an
$h_0B$-bimodule in a natural way. Massey products are defined for quadratic pair algebras as in Definition \ref{mp1} above.

\begin{defn}\label{mp}
Let $B$ be a quadratic pair algebra. Given elements $a,b,c\in h_0 B$,
of degree $p,q,r\in\N$ with $ab=0$ and $bc=0$ the
\emph{Massey product}  is the subset
$$\grupo{a,b,c}\subset h_1B_{p+q+r},$$
which is a coset of the subgroup $$(h_1B_{p+q})c+a(h_1B_{q+r}),$$
defined as follows. Given $\bar{a}\in B_{p,0}$, $\bar{b}\in
B_{q,0}$, $\bar{c}\in B_{r,0}$ representing $a$, $b$, $c$, there exist $e\in B_{p+q,1}$, $f\in B_{q+r,1}$ such
that $\partial(e)=\bar{a}\cdot\bar{b}$,
$\partial(f)=\bar{b}\cdot\bar{c}$ and one can easily check that 
$$-e\cdot\bar{c}+\bar{a}\cdot f\in h_1B_{p+q+r}\subset B_{p+q+r,1}.$$
The coset $\grupo{a,b,c}\subset h_1B_{p+q+r}$ 
coincides with the set of elements obtained in this way for all
different choices of $\bar{a}$, $\bar{b}$, $\bar{c}$, $e$ and $f$.
\end{defn}

The next result identifies Toda brackets for ring spectra with Massey products for quadratic pair algebras.

\begin{thm}\label{a2}
There is a functor
$$\pi_{*,*}\colon \left(
\text{\emph{connective ring spectra}}\right)\To\left(\text{\emph{quadratic pair algebras}}\right)$$
together with natural isomorphisms
\begin{eqnarray*}
h_0\pi_{*,*} R&\cong&\pi_*R, \text{ of rings},\\
h_1\pi_{*,*} R&\cong&\S^{-1}\pi_*R, \text{ of bimodules},
\end{eqnarray*}
such that the Massey products in $\pi_{*,*} R$ coincide with the Toda brackets in $\pi_*R$. Moreover, using the
isomorphisms as identifications the algebraically-defined
$k$-invariant of the quadratic pair algebra $\pi_{*,*} R$ in (\ref{-w}),
$$k\colon \pi_*R\otimes\Z/2\To\S^{-1}\pi_*R,$$
coincides with the multiplication by $\eta$ where $\eta$ is the stable Hopf map.
\end{thm}

We prove this theorem in Section \ref{pt1}.

\begin{exm}\label{exmqpa}
Waldhausen defined in \cite{akttsI} the $K$-theory spectrum $K\C{W}$ of a category $\C{W}$ with cofibrations and weak
equivalences. This spectrum is a symmetric spectrum, see \cite{tchs}. Moreover, if $\C{W}$ is a strict
monoidal category with biexact tensor product then $K\C{W}$ is a ring spectrum. 
In \cite{1tk} we give a small algebraic
model  for the quadratic pair algebra $\pi_{0,*}K\C{W}$ which is generated just by the objects, weak
equivalences, and cofiber sequences in $\C{W}$.
\end{exm}

The next result is an illustrating example for computations in a quadratic pair algebra.

\begin{prop}\label{2a2}
Let $B$ be a quadratic pair algebra. Given $\alpha\in h_0B_{n,*}$ with $2\alpha=0$ the Massey product
$\grupo{2,\alpha,2}$ is defined and the
$k$-invariant satisfies
$$k(\alpha)=P(\alpha|\alpha)_H\in\grupo{2,\alpha,2}.$$
\end{prop}

\begin{proof}
Let $a\in B_{n,0}$ be a representative of $\alpha$. We choose $1+1\in B_{0,0}$ as a representative of $2\in
h_0B_{0,*}$. In order to construct an element in the Massey product
$\grupo{2,\alpha,2}$ we choose $b\in B_{n,1}$ such that $\partial(b)=a(1+1)=a+a$. Notice that
\begin{eqnarray*}
\partial(b+PH(a))&=&a+a+\partial PH(a)\\
&=&1\cdot a+1\cdot a+\partial P((1|1)_HH(a))\\
&=&(1+1)a.
\end{eqnarray*}
Here we use the axioms of a quadratic pair algebra.
These equations show that the following element is in the Massey product $\grupo{2,\alpha,2}$
\begin{eqnarray*}
-(b+PH(a))(1+1)+(1+1)b&=&-(b+PH(a)+b+PH(a))\\&&+b+b+PH(a+a)\\
&=&-PH(a)-b-PH(a)-b\\&&+b+b+PH(a+a)\\
&=&-PH(a)-PH(a)+PH(a+a)\\
&=&P(a|a)_H.
\end{eqnarray*}
Here we use the laws of a quadratic pair algebra together with the fact that $P$ is central in any
quadratic pair module.
\end{proof}

It is a result of Toda that for the sphere spectrum $S$ and $\alpha\in\pi_*S$ with $2\alpha=0$ one has
$\alpha\eta\in\grupo{2,\alpha,2}$ for $\eta$ the stable Hopf map. 
The same result for some cobordism spectra can be found in \cite{cmp}. In fact Proposition \ref{2a2} together with
Theorem \ref{a2} generalizes this result for all connective ring spectra and shows: 

\begin{cor}
For any connective ring spectrum and $\alpha\in\pi_*R$ with $2\alpha=0$ we have
$\alpha\eta\in\grupo{2,\alpha,2}$.
\end{cor}

\begin{rem} 
Given an ungraded ring $A$ and an ungraded $A$-bimodule $M$ consider the
category consisting of ungraded quadratic
pair algebras $B$ with specified isomorphisms $h_0B\cong A$ and $h_1B\cong M$ such that any element in $A$ is the
image of some $x\in B_0$ with $H(x)=0$. Morphisms are quadratic pair algebra morphisms over $A$ and under
$M$. It is proved in \cite{3mlc} that the set of connected components of this category 
is in natural bijection with the $3$-dimensional Mac Lane cohomology group
$$HML^3(A,M).$$
As we mention in Remark \ref{sh2} this cohomology group coincides with topological Hochschild cohomology of the corresponding Eilenberg-Mac Lane
ring spectrum 
$$THH^3(HA,HM).$$ 
For $R$ a connective ring spectrum the Mac Lane cohomology class of the ungraded quadratic pair algebra $\pi_{0,*}R$
$$\grupo{\pi_{0,*}R}\in HML^3(\pi_0R,\pi_1R)$$
corresponds in $THH$ to the first Postnikov invariant of the ring spectrum $R$, as studied in \cite{htars}. Compare Remarks
\ref{sh} and \ref{sh2}.
\end{rem}

\begin{rem}\label{per2}
Given an ungraded quadratic pair algebra $B$ the permutative category $\C{P}B$ is actually a \emph{ring
category}, i.e. it has an additional monoidal structure $\otimes$ satisfying the axioms of \cite[3.3]{rmailst}.
The additional monoidal structure $\otimes$ is given by the product $\cdot$ on the set of objects, which is $B_0$. On
morphisms it is given by
\begin{eqnarray*}
(x_0,s_0)\otimes(x_1,s_1)&=&(x_0\cdot x_1,s_0\cdot\partial(s_1))\\
&=&(x_0\cdot x_1,\partial(s_0)\cdot s_1).
\end{eqnarray*}
The monoidal structure $\otimes$ is strictly right distributive
\begin{eqnarray*}
(x_1\oplus x_2)\otimes (x_1\oplus x_3)&=&x_1\otimes(x_2\oplus x_3),
\end{eqnarray*}
and there is a natural left distributivity isomorphism
\begin{equation*}
(x_1\cdot x_3+x_2\cdot x_3,P((x_2|x_1)_H\cdot H(x_3)))\colon (x_1\otimes x_3)\oplus (x_2\otimes x_3) \To
(x_1\oplus x_2)\otimes x_3.
\end{equation*}
\end{rem}

\section{Modules over quadratic pair algebras and module spectra}\label{6}

In this section we explicitly describe modules over quadratic pair algebras and we show their connection with modules
over ring spectra.

\begin{defn}\label{mod}
Let $B$ be a quadratic pair algebra. A \emph{right $B$-module} $M$ 
is an $\N$-graded quadratic pair
module together with multiplications, $n,m\geq 0$,
$$\begin{array}{c}
M_{n,0}\times B_{m,0}\st{\cdot}\To M_{n+m,0},\\
M_{n,0}\times B_{m,1}\st{\cdot}\To M_{n+m,1},\\
M_{n,1}\times B_{m,0}\st{\cdot}\To M_{n+m,1},\\
M_{n,ee}\times B_{m,ee}\st{\cdot}\To M_{n+m,ee}.
\end{array}$$
These multiplications are associative with respect to the multiplications in $B$.
Moreover, $1\in B_{0,0}$ acts trivially on $M_{*,0}$ and $M_{*,1}$, and  $(1|1)_H\in B_{0,ee}$ acts trivially on
$M_{*,ee}$. Furthermore, all equations in Definition \ref{alg1} hold if we replace the elements on the left of
any multiplication $\cdot$ by elements in $M$.

Notice that $h_0M$ and $h_1M$ are naturally right $h_0B$-modules, and there is a natural right $h_0B$-module morphism 
\begin{equation}\label{mm}
\cdot\,\colon h_0M\otimes_{h_0B}h_1B\r h_1M.
\end{equation}
This morphism extends the $k$-invariant of $M$ since $k(x)=x\cdot P(1|1)_H$ for $x\in h_0M$.
Moreover, given $a\in h_0M$ and $b,c\in h_0B$
of degree $p,q,r,$ such that $ab=0$ and $bc=0$ there is defined as in Definition \ref{mp} a \emph{Massey product}
$$\grupo{a,b,c}\subset h_1M_{p+q+r}$$
which is a coset of
$$(h_1M_{p+q})c+a(h_1B_{q+r}).$$
\end{defn}

\begin{rem}
In \cite{2hg3} we define the monoidal category $\C{qpm}^\N_H$ of $\N$-graded $0$-good quadratic pair modules.
For $B$ a $0$-good quadratic pair algebra the $0$-good right $B$-modules are exactly the right modules over the monoid $B$ in
$\C{qpm}_H^\N$, compare Remark \ref{al}.
\end{rem}

Recall that for an $\N$-graded ring $A$ the \emph{desuspension} $\S^{-1}M$ of a right $A$-module
$M$ is $(\S^{-1}M)_n=M_{n+1}$, $n\geq 0$, with the right $A$-module structure defined by the formula
$(\S^{-1}m)\cdot a=\S^{-1}(m\cdot a)$.

In the following theorem we consider modules over quadratic pair algebras associated to right modules over a
ring spectrum. 

\begin{thm}\label{m2}
Let $R$ be a connective ring spectrum. There is a functor
$$\pi_{*,*}\colon\left(\text{\emph{connective right $R$-modules}}\right)\To
\left(\text{\emph{right $(\pi_{*,*}R)$-modules}}\right).$$
Here the quadratic pair algebra $\pi_{*,*}R$ is obtained by using the functor in Theorem \ref{a2}.
Moreover, if we use the first isomorphism in Theorem \ref{a2} as an identification 
then for any right $R$-module $M$ there are natural isomorphisms of right $\pi_*R$-modules
\begin{eqnarray*}
h_0\pi_{*,*}M&\cong&\pi_*M,\\
h_1\pi_{*,*}M&\cong&\S^{-1}\pi_*M.
\end{eqnarray*} 
Using these
isomorphisms and the isomorphisms in Theorem \ref{a2} as identifications the algebraically-defined
homomorphism (\ref{mm}) associated to $\pi_{*,*} M$
$$\cdot\;\colon \pi_*M\otimes_{\pi_*R} \S^{-1}\pi_*R\To \S^{-1}\pi_*M$$
is defined by the right right $\pi_*R$-module structure of $\pi_*M$ according to the formula
\begin{eqnarray*}
m\cdot (\S^{-1}a)&=&(-1)^{\abs{m}}\S^{-1}(ma).
\end{eqnarray*}
In particular the $k$-invariant of $\pi_{*,*}M$ coincides with the multiplication by the stable Hopf map $\eta$.
Furthermore, Massey products in $\pi_{*,*}M$
coincide with Toda brackets in $\pi_*M$.
\end{thm}

We prove this theorem in Section \ref{pt1}.



\section{Sign group actions on quadratic pair modules}\label{7}

For the definition of an $E_\infty$-quadratic pair algebra below we need an action of the symmetric track groups
$\symt{n}$ on a quadratic pair module. For this we recall the notion of sign group and action of a sign group
from \cite{2hg2,2hg3}.

\begin{defn}\label{sigg}
A \emph{sign group} is a diagram in the category of groups
\begin{equation*}
\set{\pm1}\st{\imath}\hookrightarrow G_\vc\st{\delta}\twoheadrightarrow G\st{\varepsilon}\To \set{\pm1}
\end{equation*}
where the first two morphisms form an extension. By abuse of notation we denote this sign group just by $G_\vc$.
The group law of the groups defining a sign group is denoted multiplicatively.
\end{defn}

Sign groups as above
were introduced in \cite{2hg2}. The main examples are the \emph{symmetric track groups} $\symt{n}$ associated
to the sign homomorphism of the symmetric groups $$\sign\colon\sym{n}\r\set{\pm1}.$$ The group $\symt{n}$ has a presentation with
generators $\omega$, $t_i$, $1\leq i\leq n-1$,
and relations
\begin{eqnarray*}
t_1^2&=&1\text{ for }1\leq i\leq n-1,\\
(t_it_{i+1})^3&=&1\text{ for }1\leq i\leq n-2,\\
\omega^2&=&1,\\
t_i\omega&=&\omega t_i\text{ for }1\leq i\leq n-1,\\
t_it_j&=&\omega t_jt_i\text{ for }1\leq i<j-1\leq n-1.
\end{eqnarray*}
Moreover, the structure of sign groups is given by $\imath(-1)=\omega$, $\delta(\omega)=0$, and
$\delta(t_i)=(i\; i+1)$, the permutation exchanging
$i$ and $i+1$ in $\set{1,\dots,n}$.

\begin{defn}\label{laac}
A sign group $G_\vc$ \emph{acts on the right} of a quadratic pair module
$C$ if $G$ acts on the right of $C$ by morphisms $g^*\colon C\r
C$, $g\in G$, in $\C{qpm}$, and there is a bracket
$$\grupo{-,-}=\grupo{-,-}_{G}\colon C_0\times G_\vc\To C_1$$
satisfying the following properties, $x,y\in C_0$, $z\in C_1$,
$s,t\in G_\vc$.
\begin{enumerate}
\item
$\grupo{x+y,t}=\grupo{x,t}+\grupo{y,t}+P(-{\delta(t)^*(x)}+\varepsilon\delta(t)x|{\delta(t)^*(y)})_H$,
\item $\varepsilon\delta(t)x+\binom{\varepsilon\delta(t)}{2}\partial PH(x)=\delta(t)^*(x)+\partial\grupo{x,t}$,
\item $\varepsilon\delta(t)z+\binom{\varepsilon\delta(t)}{2}PH\partial(x)=\delta(t)^*(z)+\grupo{\partial(z),t}$,
\item $\grupo{x,s\cdot
t}=\grupo{\delta(s)^*(x),t}+\grupo{\varepsilon\delta(t)x+\binom{\varepsilon\delta(t)}{2}\partial PH(x),s}$,
\item for the element $\omega=\imath(-1)\in G_\vc$ we have the
\emph{$\omega$-formula}:
$$\grupo{x,\omega}=P(x|x)_H.$$
Notice that the $\omega$-formula corresponds to the $k$-invariant in (\ref{-w}).
\end{enumerate}
Here $\varepsilon\delta(t)$ can take the value $-1$. In this case $(-1)x$ means $-x\in C_0$, and $(-1)z=-z\in C_1$.
\end{defn}

Given a sign group $G_\vc$ the \emph{``group ring''} $A(G_\vc)$ is the ungraded 
quadratic pair algebra with
generators
\begin{itemize}
\item $[g]$ for any $g\in G$ on the $0$-level,
\item $[t]$ for any $t\in G_\vc$ on the $1$-level,
\item no generators on the $ee$-level,
\end{itemize}
and relations
\begin{itemize}
\item $H[g]=0$ for $g\in G$,
\item $[1]=1$ the unit element,
\item $[gh]=[g]\cdot[h]$ for $g,h\in G$,
\item $\partial[t]=-[\delta(t)]+\varepsilon\delta(t)$,
\item $[st]=[\delta(s)]\cdot[t]+\varepsilon\delta(t)\cdot[s]
+\binom{\varepsilon\delta(s)}{2}\binom{\varepsilon\delta(t)}{2}P(1|1)_H$ for $s,t\in G_\vc$,
\item $[\omega]=P(1|1)_H$ where $\omega=\imath(-1)$.
\end{itemize}
The relations above show that $A_{(0)}(G_\vc)=\Z_\nill[G_+]$ where $G_+$ is the group $G$ 
together with an outer base point, so $A(G_\vc)$ is $0$-free and hence $0$-good.

If $C$ is a $0$-good quadratic pair module then an action of $G_\vc$ on $C$ is the same as a right
ungraded $A(G_\vc)$-module structure on $C$. The correspondence is given by the following equations, see \cite{2hg3}. Given
$g\in G$, $t\in G_\vc$,
$x\in C_0$, $y\in C_1$, $z\in C_{ee}$,
\begin{eqnarray*}
g^*x&=&x\cdot[g],\\
g^*y&=&y\cdot[g],\\
g^*z&=& z\cdot \Delta[g],\\
\grupo{x,t}&=&x\cdot[t]+\binom{\varepsilon\delta(t)}{2}PH(x).
\end{eqnarray*}

\begin{defn}
A morphism $f_\vc\colon G_\vc\r K_\vc$ in the category $\C{Gr}_\pm$ of sign groups is given by two group
homomorphisms $f \colon G \r K $, $f_\vc\colon G_\vc\r K_\vc$ with $f_\vc\imath = \imath$, $f\delta=\delta
f_\vc$, and $\varepsilon =\varepsilon f$. In this situation we say that the homomorphism $f_\vc$ \emph{covers}
$f$.
A \emph{twisted bilinear morphism} of sign groups
$$(f_\vc,g_\vc)\colon G_\vc\times L_\vc\To K_\vc$$
is given by a pair of sign group morphisms
$f_\vc\colon G_\vc\r K_\vc$, $g_\vc\colon L_\vc\r K_\vc$,
such that given $a\in G$ and $b\in L$ the equality
$$f(a)g(b)=g(b)f(a)$$ holds in $K$, and therefore the group homomorphism
$$(f,g)\colon G\times L\To K\colon (a,b)\mapsto f(a)g(b)$$
is defined, and given $x\in G_\vc$, $y\in L_\vc$ the following
equality is satisfied in $K_\vc$
$$f_\vc(x)g_\vc(y)=g_\vc(y)f_\vc(x)\imath\left((-1)^{\binom{\varepsilon\delta(x)}{2}\binom{\varepsilon\delta(y)}{2}}\right).$$
\end{defn}

There is a universal twisted bilinear morphism
\begin{equation*}
(i_{G_\vc},i_{L_\vc})\colon G_\vc\times L_\vc\To
G_\vc\tilde{\times}L_\vc.
\end{equation*}
Here $G_\vc\tilde{\times}L_\vc$ is the \emph{twisted product} of sign groups defined in \cite{2hg3} which is a
sign group
$$\set{\pm1}\st{\imath}\hookrightarrow G_\vc\tilde{\times}L_\vc\twoheadrightarrow G\times
L\st{\varepsilon}\twoheadrightarrow\set{\pm1}.$$
The twisted product defines a non-symmetric monoidal structure in $\C{Gr}_\pm$.
We show in \cite{2hg3} that the ``group ring'' $A(G_\vc)$ is a strict monoidal functor, i.e. the ``group ring''
of a twisted product
of sign groups is the tensor product of the ``group rings'' of the factors and the ``group ring'' of the trivial
sign group is $\ol{\Z}_\nill$.

\begin{rem}
The isomorphism class of a sign group $G_\vc$ is determined by two group cohomology classes
$$\varepsilon_{G_\vc}\in H^1(G,\Z/2),\;\;\set{G_\vc}\in H^2(G,\Z/2),$$
corresponding to the homomorphism $\varepsilon$ and the extension in Definition \ref{sigg}, respectively.
For the symmetric track groups $\symt{n}$ these cohomology classes are the first two Stiefel-Whitney classes
associated to the usual inclusion of the symmetric group into the orthogonal group $\sym{n}\subset O(n)$.
The cohomological
invariants of $G_\vc\tilde{\times}L_\vc$ are given by
\begin{eqnarray*}
\varepsilon_{G_\vc\tilde{\times} L_\vc}&=&p_1^*\varepsilon_{G_\vc}+
p_2^*\varepsilon_{L_\vc}\in H^1(G\times L,\Z/2),\\
\set{G_\vc\tilde{\times} L_\vc}&=&p_1^*\set{G_\vc}
+p_1^*\varepsilon_{G_\vc}\smile p_2^*\varepsilon_{L_\vc}+p_2^*\set{L_\vc}\in H^2(G\times L,\Z/2).
\end{eqnarray*}
Here $p_1$ and $p_2$ are the projections of the factors of the product $G\times K$ and $\smile$ denotes the
cup-product.
\end{rem}

We define in \cite{2hg3} sign group morphisms
\begin{equation}\label{sgx}
\times_\vc\colon\symt{n}\tilde{\times}\symt{m}\To\symt{n+m}
\end{equation}
covering the usual cross product homomorphisms
\begin{equation}\label{cp}
\sym{n}{\times}\sym{m}\To\sym{n+m}\colon (\sigma,\tau)\mapsto \sigma\times\tau.
\end{equation}
Recall that the block permutation $\sigma\times\tau$ permutes the first $n$ elements $\set{1,\dots,n}$ of 
$\set{1,\dots,n,n+1,\dots,n+m}$ according to $\sigma$ and the last $m$ elements of this set $\set{n+1,\dots,n+m}$ according to
$\tau$.
The sign group morphism (\ref{sgx}) corresponds to the twisted bilinear morphism
\begin{equation*}
(-\wedge S^m,S^n\wedge-)\colon\symt{n}\times\symt{m}\To\symt{n+m}
\end{equation*}
defined on generators by
\begin{eqnarray*}
t_i\wedge S^m&=&t_i,\;\;1\leq i\leq n,\\
\omega\wedge S^m&=&\omega,\\
S^n\wedge t_i&=&t_{n+i},\;\;1\leq i\leq m,\\
S^n\wedge \omega&=&\omega.
\end{eqnarray*}
The sign group morphism (\ref{sgx}) was defined in \cite{2hg3} in a geometric way. We obtained this description in terms
of generators by using the positive pin representations in \cite{2hg2} and the computations in
\cite[Section 17]{2hg3}, see \cite{2hg5} for further details.

\section{$E_\infty$-quadratic pair algebras associated to commutative ring spectra}\label{seqpa}\label{8}

In order to deal with secondary operations in the homotopy groups of a connective commutative ring spectrum 
in general we need the following algebraic structure.

\begin{defn}\label{eqpa}
An \emph{$E_\infty$-quadratic pair algebra} is a 
quadratic pair algebra $C$ as in Definition \ref{alg1} 
together with a cup-one product operation
$$\smile_1\colon C_{n,0}\times C_{m,0}\To C_{n+m,1},\;\;n,m\geq 0,$$
such that the quadratic pair module $C_{n,*}$ is endowed with a right action of the symmetric track group $\symt{n}$
and the following compatibility conditions hold.
Let $x_i\in C_{n_i,0}$, $s_i\in C_{n_i,1}$, $a_i\in C_{n_i,ee}$, 
$g_i,g'_i\in\sym{n_i}$, and $r_i\in\symt{n_i}$. The product in the quadratic pair algebra $C$ is compatible with the sign
group actions and the sign group morphisms $\times_\vc$ in (\ref{sgx}),
\begin{eqnarray}\label{equi}
(x_1\cdot[g_1])\cdot(x_2\cdot[g_2])&=&(x_1\cdot x_2)\cdot [g_1\times g_2],\\
\nonumber (s_1\cdot[g_1])\cdot(x_2\cdot[g_2])&=&(s_1\cdot x_2)\cdot [g_1\times g_2],\\
\nonumber (a_1\cdot([g_1]|[g_1'])_H)\cdot(a_2\cdot([g_2]|[g_2'])_H)&=&
\nonumber (a_1\cdot a_2)\cdot ([g_1\times g_2]| [g_1'\times g_2'])_H,\\
\nonumber x_1\cdot (x_2\cdot[r_2])&=&(x_1\cdot x_2)\cdot[S^{n_1}\wedge r_2].
\end{eqnarray}
The cup-one product measures the lack of commutativity, i.e. if $\tau_{p,q}\in\sym{p+q}$ denotes the permutation
exchanging the blocks $\set{1,\dots,p}$ and $\set{p+1,\dots,p+q}$, $p,q\geq 0$, then
\begin{eqnarray}\label{commh} 
&&\\
\nonumber (x_2\cdot x_1)\cdot[\tau_{n_1,n_2}]+\partial(x_1\smile_1x_2)&=&x_1\cdot x_2+\partial P(H(x_2)\cdot
TH(x_1))\cdot[\tau_{n_1,n_2}],\\
\nonumber (x_2\cdot s_1)\cdot[\tau_{n_1,n_2}]+\partial(s_1)\smile_1x_2&=&
s_1\cdot x_2+P(H(x_2)\cdot TH\partial(s_1))\cdot[\tau_{n_1,n_2}].
\end{eqnarray}
The cup-one product is itself commutative in the following sense
\begin{eqnarray}\label{comcup1}
(x_2\smile_1x_1)\cdot[\tau_{n_1,n_2}]+x_1\smile_1x_2&=&-P(TH(x_1)\cdot H(x_2))\\
\nonumber &&+P(H(x_2)\cdot
TH(x_1))\cdot[\tau_{n_1,n_2}].
\end{eqnarray}
Let $1_n\in\sym{n}$ be the unit element.
The cup-one product also satisfies the following rules with respect to addition
\begin{eqnarray}\label{cup1+}
x_1\smile_1(x_2+x_3)&=&x_1\smile_1x_2+x_1\smile_1x_3\\
\nonumber &&+P(\partial(x_1\smile_1 x_2)|(x_3\cdot x_1)\cdot[\tau_{n_1,n_3}])_H,
\end{eqnarray}
multiplication
\begin{eqnarray}\label{cup1prod}
&&\\
\nonumber (x_1\cdot x_2)\smile_1x_3&=&((x_1\smile_1 x_3)\cdot x_2)\cdot[1_{n_1}\times\tau_{n_2,n_3}]+x_1\cdot (x_2\smile_1x_3)\\
\nonumber &&+P((\partial(x_1\smile_1 x_3)|(x_3\cdot x_1)\cdot[\tau_{n_1,n_3}])_H\cdot
H(x_2))\cdot[1_{n_1}\times\tau_{n_2,n_3}]\\
\nonumber &&+P(H(x_3)\cdot(x_1|x_1)_H\cdot TH(x_2))\cdot[\tau_{n_1+n_2,n_3}]\\
\nonumber &&-P((x_1|x_1)_H\cdot H(x_3)\cdot TH(x_2))\cdot [1_{n_1}\times\tau_{n_2,n_3}], 
\end{eqnarray}
and symmetric group action
\begin{eqnarray}\label{cup1sga}
(x_1\cdot[g_1])\smile_1(x_2\cdot[g_2])&=&(x_1\smile_1x_2)\cdot[g_1\times g_2].
\end{eqnarray}
\end{defn}

\begin{rem}\label{qia}
One can check that an $E_\infty$-pair algebra $C$ as defined in Definition \ref{eipa} is the same as an $E_\infty$-quadratic pair algebra such that $C_{*,ee}=0$ and for any $x\in
C_{n,0}$ and $t\in\symt{n}$ the equation $x\cdot[t]=0$ holds. This defines a full inclusion of categories
\begin{eqnarray*}
\left(E_\infty\text{-pair algebras}\right)&\subset&\left(E_\infty\text{-quadratic pair algebras}\right).
\end{eqnarray*}
\end{rem}

\begin{lem}
Given an $E_\infty$-quadratic pair algebra $C$ the ring $h_0C$ is commutative and $h_1C$ is an $h_0C$-module in a natural way.
\end{lem}

\begin{proof}
Let $x_i\in C_{n_i,0}$. By the laws of an $E_\infty$-quadratic pair algebra the class of $x_1\cdot x_2$ in $h_0C$
coincides with the class of $(x_2\cdot x_1)\cdot[\tau_{n_1,n_2}]$. It is easy to see that $\sign \tau_{n_1,n_2}=(-1)^{n_1n_2}$. Let
$\hat{\tau}\in\symt{n_1+n_2}$ be an element with $\delta(\hat{\tau})=\tau_{n_1,n_2}$. By the laws of a 
module over a quadratic pair algebra and the definition of the ``group ring'' of a sign group
\begin{eqnarray*}
\partial((x_2\cdot x_1)\cdot[\hat{\tau}])&=&(x_2\cdot x_1)\cdot\partial[\hat{\tau}]\\
&=&(x_2\cdot x_1)\cdot(-[\tau_{n_1,n_2}]+(-1)^{n_1n_2})\\
&=&-(x_2\cdot x_1)\cdot[\tau_{n_1,n_2}]+(-1)^{n_1n_2}(x_2\cdot x_1),
\end{eqnarray*}
therefore $(x_2\cdot x_1)\cdot[\tau_{n_1,n_2}]$ represents the same class as $(-1)^{n_1n_2}(x_2\cdot x_1)$ in
$h_0C$, so $h_0C$ is commutative.

Let now $s_1\in C_{n_1,1}$ be an element with $\partial(s_1)=0$. By the laws of an $E_\infty$-quadratic pair algebra
$(x_2\cdot s_1)\cdot[\tau_{n_1,n_2}]=s_1\cdot x_2$. Moreover, by the laws of quadratic pair algebras,  
modules over a quadratic pair algebra, and the definition of the ``group ring'' of a sign group
\begin{eqnarray*}
0&=&
x_2\cdot\partial(s_1)\cdot[\hat{\tau}]\\
&=&x_2\cdot s_1\cdot\partial[\hat{\tau}]\\
&=&x_2\cdot s_1\cdot(-[\tau_{n_1,n_2}]+(-1)^{n_1n_2})\\
&=&-(x_2\cdot s_1)\cdot[\tau_{n_1,n_2}]+(-1)^{n_1n_2}(x_2\cdot s_1),
\end{eqnarray*}
hence we are done.
\end{proof}

As in the non-quadratic case $E_\infty$-quadratic pair algebras carry a notion of cup-one square.

\begin{defn}\label{cup1}
Let $C$ be an $E_\infty$-quadratic pair algebra. Given an element $a\in h_{0}C_{2n,*}$ we define the \emph{cup-one square} of $a$
$$Sq_1(a)\in h_1C_{4n,*}$$
in the following way. 
Choose a representative $\bar{a}\in C_{2n,0}$ of $a$ and an element in the symmetric track group $\hat{\tau}\in \symt{4n}$ whose
boundary is the shuffle permutation
$\delta(\hat{\tau})=\tau_{2n,2n}$. 
Then we have
\begin{eqnarray*}
\partial(-(\bar{a}\cdot \bar{a})\cdot[\hat{\tau}]+\bar{a}\smile_1\bar{a}-P(H(\bar{a})\cdot TH(\bar{a}))\cdot[\tau_{2n,2n}])&=&0.
\end{eqnarray*}
so we can set
$$Sq_1(a)=-(\bar{a}\cdot \bar{a})\cdot[\hat{\tau}]+\bar{a}\smile_1\bar{a}-P(H(\bar{a})\cdot
TH(\bar{a}))\cdot[\tau_{2n,2n}]\in h_1C_{4n,*}.$$
We leave it to the reader to check that this element does not depend on the choice of $\bar{a}$. However it does depend
on the choice of $\hat{\tau}$. There are two possible choices, namely $\hat{\tau}$ and
$\omega\hat{\tau}$. The next result determines the difference between the two possible cup-one squares.
\end{defn}

\begin{lem}
Let $Sq_1$ be cup-one square in an $E_\infty$-quadratic pair algebra $C$ associated to the lift $\hat{\tau}$ of the shuffle
permutation and let $Sq_1^\omega$ be the cup-one square associated to $\omega\hat{\tau}$. Then given $a\in h_0C_{2n,*}$
\begin{eqnarray*}
Sq^\omega_1(a)&=&Sq_1(a)+k(a\cdot a).
\end{eqnarray*}
Here $k$ is the $k$-invariant defined in (\ref{-w}).
\end{lem}

\begin{proof}
Let $\bar{a}\in C_{2n,0}$ be a representative of $a$. The laws of a quadratic pair algebra, of a module over it, and the definition of the ``group ring'' of a sign group yield equations
\begin{eqnarray*}
(\bar{a}\cdot\bar{a})\cdot [\omega\hat{\tau}]&=&(\bar{a}\cdot\bar{a})\cdot([\hat{\tau}]+[\omega])\\
&=&(\bar{a}\cdot\bar{a})\cdot[\hat{\tau}]+(\bar{a}\cdot\bar{a})\cdot[\omega]\\
&=&(\bar{a}\cdot\bar{a})\cdot[\hat{\tau}]+(\bar{a}\cdot\bar{a})\cdot P(1|1)_H\\
&=&(\bar{a}\cdot\bar{a})\cdot[\hat{\tau}]+P(\bar{a}\cdot\bar{a}|\bar{a}\cdot\bar{a})_H,
\end{eqnarray*}
hence the lemma follows from the very definition of the cup-one product and the $k$-invariant of a quadratic pair
module.
\end{proof}

One of the main results in this paper is the following theorem.

\begin{thm}\label{ca2}
There is a commutative diagram of functors
$$\xymatrix{{\left(\begin{array}{c}
\text{\emph{connective commutative}}\\
\text{\emph{ring spectra}}
\end{array}\right)}
\ar[d]_{\text{inclusion}}\ar[r]^-{\pi_{*,*}}&
{\left(E_\infty\text{\emph{-quadratic pair algebras}}\right)}
\ar[d]^{\text{forget}}\\
{\left(\text{\emph{connective ring spectra}}\right)}
\ar[r]^-{\pi_{*,*}}&
{\left(\text{\emph{quadratic pair algebras}}\right)}}$$
Here the lower arrow is the functor in Theorem \ref{a2}. Moreover, for a commutative ring spectrum $Q$
the algebraic cup-one squares in $\pi_{*,*} Q$ correspond to the topologically-defined cup-one
squares in $\pi_*Q$.
\end{thm}

We prove this theorem in Section \ref{pt1}.

\begin{rem}\label{caen}
$E_\infty$-quadratic pair algebras can also be described as coherently homotopy commutative monoids as in the
additive case, see Remark \ref{hex}. For this one uses the symmetric monoidal category $\C{qpm}_H^{\symtt}$ of enhanced symmetric
sequences of $0$-good quadratic pair modules introduced in \cite{2hg3}. In fact, this is the way in which we
obtain the equations in Definition \ref{eqpa}. The cup-one product $x_1\smile_1x_2$ of two level $0$ elements
corresponds to the image of the right linear generator $x_1 \ul{\circledcirc} x_2$ in the tensor product of
square groups $C_{*,(0)}\odot C_{*,(0)}=(C\odot C)_{*,(0)}$ under the track $\smile_1$. The definition of tracks
for quadratic pair modules is recalled in Definition \ref{altra} below.

As we recall in Remark \ref{BE} the coherence conditions in Remark \ref{hex} (i.e. the hexagon and the idempotence) are encoded by the
Barratt-Eccles operad in the simplicial setting. Indeed, $E_\infty$-quadratic pair algebras are also the algebras 
over an operad in $\C{qpm}_H$ which can be
obtained from the simplicial Barratt-Eccles operad, i.e. the additive situation described in Remark \ref{BE}
extends to the quadratic world.
\end{rem}

\begin{rem}\label{per3}
If $C$ is an ungraded $E_\infty$-quadratic pair algebra then the permutative category $\C{P}C$ in Remark
\ref{per1} is actually a \emph{bipermutative category} in the sense of \cite[3.6]{rmailst}. The monoidal structures $\oplus$ and
$\otimes$ are given by Remarks \ref{per1} and \ref{per2}, and the symmetry isomorphism for $\otimes$ is
\begin{equation*}
(x_1\cdot x_2,P(H(x_2)\cdot TH(x_1))-x_1\smile_1x_2)\colon x_1\otimes x_2\To x_2\otimes x_1.
\end{equation*}
\end{rem}

\section{Algebras over commutative ring spectra}\label{9}

Given a connective commutative ring spectrum $Q$ and a $Q$-algebra $R$ we study in this section the algebraic 
properties of the induced
action of the $E_\infty$-quadratic pair algebra $\pi_{*,*}Q$ defined by Theorem \ref{ca2} on the quadratic
pair algebra $\pi_{*,*}R$ given by Theorem \ref{a2}, obtaining Theorem \ref{app} as an application. This
generalizes the results of Section \ref{31} to spectra which do not neglect the Hopf map.

\begin{defn}\label{meqpa}
Let $C$ be an $E_\infty$-quadratic pair algebra. A \emph{$C$-algebra} is a
quadratic pair algebra $B$ as in Definition \ref{alg1} 
together with a quadratic pair algebra morphism $u\colon C\r B$, called the \emph{unit}, and 
two external cup-one product operations, $n,m\geq 0$,
\begin{equation}\label{extcup1}
\begin{array}{c}
\smile_1\colon B_{n,0}\times C_{m,0}\To C_{n+m,1},\\
\smile_1\colon C_{n,0}\times B_{m,0}\To C_{n+m,1},
\end{array}
\end{equation}
such that the quadratic pair module $B_{n,*}$ is endowed with a right action of the symmetric track group $\symt{n}$
and the following compatibility conditions hold. Let 
$x_i\in B_{n_i,0}$, $s_i\in B_{n_i,1}$, $a_i\in B_{n_i,ee}$, 
$\tilde{x}_i\in C_{n_i,0}$, $\tilde{s}_i\in C_{n_i,1}$, $\tilde{a}_i\in C_{n_i,ee}$, 
$g_i,g'_i\in\sym{n_i}$, and $r_i\in\symt{n_i}$.
The product in $B$ is equivariant, i.e. it satisfies
equations (\ref{equi}) and
\begin{eqnarray*}
(x_1\cdot[r_1])\cdot x_2&=&(x_1\cdot x_2)\cdot [r_1\wedge S^{n_2}].
\end{eqnarray*}
The morphism $u$ is also equivariant, 
\begin{eqnarray*}
u_0(\tilde{x}_1\cdot[g_1])&=&u_0(\tilde{x}_1)\cdot[g_1],\\
u_1(\tilde{s}_1\cdot[g_1])&=&u_1(\tilde{s}_1)\cdot[g_1],\\
u_1(\tilde{x}_1\cdot[r_1])&=&u_0(\tilde{x}_1)\cdot[r_1],
\end{eqnarray*}
and compatible with the cup-one products in $C$ and $B$
\begin{eqnarray*}
u_1(\tilde{x}_1\smile_1 \tilde{x}_2)&=&u_0(\tilde{x}_1)\smile_1\tilde{x}_2,\\
&=&\tilde{x}_1\smile_1u_0(\tilde{x}_2).
\end{eqnarray*}
The external cup-one products satisfy the following analogue of (\ref{commh}),
\begin{eqnarray*}
\nonumber (x_2\cdot u_0(\tilde{x}_1))\cdot[\tau_{n_1,n_2}]+\partial(\tilde{x}_1\smile_1x_2)&=&
u_0(\tilde{x}_1)\cdot x_2+\partial P(H(x_2)\cdot
THu_0(\tilde{x}_1))\cdot[\tau_{n_1,n_2}],\\
\nonumber (x_2\cdot u_1(\tilde{s}_1))\cdot[\tau_{n_1,n_2}]+\partial(\tilde{s}_1)\smile_1x_2&=&
u_1(\tilde{s}_1)\cdot x_2+P(H(x_2)\cdot TH\partial u_1(\tilde{s}_1))\cdot[\tau_{n_1,n_2}],\\
\nonumber (u_0(\tilde{x}_2)\cdot s_1)\cdot[\tau_{n_1,n_2}]+\partial(s_1)\smile_1 \tilde{x}_2&=&
s_1\cdot u_0(\tilde{x}_2)+P(H u_0(\tilde{x}_2)\cdot TH\partial(s_1))\cdot[\tau_{n_1,n_2}].
\end{eqnarray*}
The two external cup-one products are related by the following formula, which is similar to (\ref{comcup1}),
\begin{eqnarray*}
(\tilde{x}_2\smile_1x_1)\cdot[\tau_{n_1,n_2}]+x_1\smile_1\tilde{x}_2&=&-P(TH(x_1)\cdot Hu_0(\tilde{x}_2))\\
\nonumber &&+P(Hu_0(\tilde{x}_2)\cdot
TH(x_1))\cdot[\tau_{n_1,n_2}].
\end{eqnarray*}
in particular it is redundant to give the two external cup-one products in (\ref{extcup1}) 
since any of them determines the other, but having two makes the equations shorter.
The external cup-one products also satisfy the following rules with respect to addition
\begin{eqnarray*}
x_1\smile_1(\tilde{x}_2+\tilde{x}_3)&=&x_1\smile_1\tilde{x}_2+x_1\smile_1\tilde{x}_3+P(\partial(x_1\smile_1 \tilde{x}_2)|(u_0(\tilde{x}_3)\cdot x_1)\cdot[\tau_{n_1,n_3}])_H,\\
\tilde{x}_1\smile_1(x_2+x_3)&=&\tilde{x}_1\smile_1x_2+\tilde{x}_1\smile_1x_3+P(\partial(\tilde{x}_1\smile_1 x_2)|(x_3\cdot u_0(\tilde{x}_1))\cdot[\tau_{n_1,n_3}])_H,
\end{eqnarray*}
compare (\ref{cup1+}),
multiplication
\begin{eqnarray*}
\nonumber (\tilde{x}_1\cdot \tilde{x}_2)\smile_1x_3&=&((\tilde{x}_1\smile_1 x_3)\cdot
u_0(\tilde{x}_2))\cdot[1_{n_1}\times\tau_{n_2,n_3}]+u_0(\tilde{x}_1)\cdot (\tilde{x}_2\smile_1x_3)\\
\nonumber &&+P((\partial(\tilde{x}_1\smile_1 x_3)|(x_3\cdot u_0(\tilde{x}_1))\cdot[\tau_{n_1,n_3}])_H\cdot
Hu_0(\tilde{x}_2))\cdot[1_{n_1}\times\tau_{n_2,n_3}]\\
\nonumber &&+P(H(x_3)\cdot u_{ee}(\tilde{x}_1|\tilde{x}_1)_H\cdot THu_0(\tilde{x}_2))\cdot[\tau_{n_1+n_2,n_3}]\\
\nonumber &&-P(u_{ee}(\tilde{x}_1|\tilde{x}_1)_H\cdot H(x_3)\cdot TH u_0(\tilde{x}_2))\cdot [1_{n_1}\times\tau_{n_2,n_3}], \\
\nonumber (x_1\cdot x_2)\smile_1\tilde{x}_3&=&((x_1\smile_1 \tilde{x}_3)\cdot
x_2)\cdot[1_{n_1}\times\tau_{n_2,n_3}]+x_1\cdot (x_2\smile_1\tilde{x}_3)\\
\nonumber &&+P((\partial(x_1\smile_1 \tilde{x}_3)|(u_0(\tilde{x}_3)\cdot x_1)\cdot[\tau_{n_1,n_3}])_H\cdot
H(x_2))\cdot[1_{n_1}\times\tau_{n_2,n_3}]\\
\nonumber &&+P(Hu_0(\tilde{x}_3)\cdot(x_1|x_1)_H\cdot TH(x_2))\cdot[\tau_{n_1+n_2,n_3}]\\
\nonumber &&-P((x_1|x_1)_H\cdot Hu_0(\tilde{x}_3)\cdot TH(x_2))\cdot [1_{n_1}\times\tau_{n_2,n_3}], 
\end{eqnarray*}
compare (\ref{cup1prod}),
and symmetric group action
\begin{eqnarray*}
(x_1\cdot[g_1])\smile_1(\tilde{x}_2\cdot[g_2])&=&(x_1\smile_1\tilde{x}_2)\cdot[g_1\times g_2],\\
(\tilde{x}_1\cdot[g_1])\smile_1(x_2\cdot[g_2])&=&(\tilde{x}_1\smile_1x_2)\cdot[g_1\times g_2],
\end{eqnarray*}
compare (\ref{cup1sga}).

The ring $h_0B$ is an $h_0C$-algebra in such a way that the $h_0B$-bimodule $h_1B$ restricts to an $h_0C$-module.
\end{defn}

\begin{rem}
Any $E_\infty$-quadratic pair algebra is an algebra over itself. Moreover, given a morphism of $E_\infty$-quadratic
pair algebras $f\colon \bar{C}\r C$ and a $C$-algebra $B$ then
$B$ has a $\bar{C}$-algebra structure, denoted by $f^*B$, with unit $uf\colon\bar{C}\r B$ and 
external cup-one products obtained by precomposing with $f$.
\end{rem}

The main new secondary homotopy operation on algebras over $E_\infty$-quadratic pair algebras is introduced in the
next definition, analogous to Definition \ref{deride1}.

\begin{defn}\label{deride}
Let $C$ be an $E_\infty$-quadratic pair algebra and let $B$ be a $C$-algebra. We associate to any pair $(x,y)$
with $x\in C_{0,n}$ and $y\in B_{n,1}$ such that $\partial(y)=u_0(x)$
the degree $n$ derivation $$\theta_{(x,y)}\colon h_0B\To h_1B$$
defined by
\begin{eqnarray*}
\theta_{(x,y)}(a)&=&-y\cdot\bar{a}+(\bar{a}\cdot y)\cdot[\tau_{n,\abs{a}}]+x\smile_1\bar{a}-P(H(\bar{a})\cdot
THu_0(x)).
\end{eqnarray*}
Here $\bar{a}\in B_{*,0}$ is a representative of $a\in h_0B$. This is indeed an element in $h_1B$ by the
equations defining a $C$-algebra which are analogue to (\ref{commh}), see Definition \ref{meqpa}. Moreover, by (\ref{commh})
$\theta_{(x,y)}=0$ provided $B$ is an $E_\infty$-quadratic pair algebra and $u$ is a morphism of
$E_\infty$-quadratic pair algebras.
One readily checks that $\theta_{(x,y)}$ is indeed a derivation and that $\theta_{(x,y)}(a)$ does not depend
on the choice of $\bar{a}$.
\end{defn}

\begin{prop}\label{loh}
Given an $E_\infty$-quadratic pair algebra $C$ and a $C$-algebra $B$ let $I_C(B)$ be the kernel of the ring
homomorphism $h_0u\colon h_0C\r h_0B$. There is an $h_0C$-module homomorphism
$$\theta\colon I_C(B)/I_C(B)^2\To HH^1_{h_0C}(h_0B,h_1B)$$
which sends an element in $I_C(B)/I_C(B)^2$
represented by $x\in C_{n,0}$ to the element in Hochschild cohomology represented by a derivation
$\theta_{(x,y)}$ as in Definition \ref{deride} for any $y\in B_{n,1}$ with $\partial(y)=u_0(x)$. The homomorphism
$\theta$ is natural in $B$ and in $C$ in the obvious way. In particular $\theta=0$ when $B$ is an $E_\infty$-quadratic pair algebra and $u$ is a morphism of
$E_\infty$-quadratic pair algebras.
\end{prop}

The proof of this proposition is tedious but straightforward. 
It only uses the
equations defining algebras over $E_\infty$-quadratic pair algebras, compare the proof of the additive analogue, Proposition \ref{loh1}.

Algebras over commutative ring spectra yield examples of algebras over $E_\infty$-quadratic pair
algebras as the following theorem shows.

\begin{thm}\label{ae2}
Let $Q$ be a connective commutative ring spectrum and let $\pi_{*,*}Q$ be the $E_\infty$-quadratic pair algebra
given by Theorem \ref{ca2}. There is a diagram of functors
$$\xymatrix{{\left({\emph{connective $Q$-algebras}}\right)}
\ar[d]_{\text{forget}}\ar[r]^-{\pi_{*,*}}&
{\left(\pi_{*,*}Q\text{\emph{-algebras}}\right)}
\ar[d]^{\text{forget}}\\
{\left(\text{\emph{connective ring spectra}}\right)}
\ar[r]^-{\pi_{*,*}}&
{\left(\text{\emph{quadratic pair algebras}}\right)}}$$
which commutes up to quasi-isomorphisms.
Here the lower arrow is the functor in Theorem \ref{a2}.
Moreover, if $f\colon Q'\r Q$ is a morphism of connective commutative ring spectra then the square
$$\xymatrix{{\left({\emph{connective $Q$-algebras}}\right)}
\ar[d]_{f^*}\ar[r]^-{\pi_{*,*}}&
{\left(\pi_{*,*}Q\text{\emph{-algebras}}\right)}
\ar[d]^{(\pi_{*,*}f)^*}\\
{\left({\emph{connective $Q'$-algebras}}\right)}
\ar[r]^-{\pi_{*,*}}&
{\left(\pi_{*,*}Q'\text{\emph{-algebras}}\right)}}$$
commutes up to quasi-isomorphisms. Furthermore, if $R$ is a connective $Q$-algebra, $a\in\pi_mR$, and
$b\in\pi_nQ$ is in $I_Q(R)=\ker[\pi_*Q\r\pi_*R]\cong I_{\pi_{*,*}Q}(\pi_{*,*}R)$, then  the element
$\theta_{(x,y)}(a)\in\pi_{n+m+1}R\cong(h_1\pi_{*,*}R)_{n+m}$ defined in (\ref{teta}) can be identified with
$\theta_{(x,y)}(a)$ in the sense of Definition \ref{deride}.
\end{thm}

The proof of this theorem goes along the same lines as the proof of Theorem \ref{ca2}. 
The reader 
can find the crucial hints in the following remark.

\begin{rem}\label{hint}
Suppose that $Q$ is a connective commutative ring spectrum, $R$ is a connective $Q$-algebra, and 
$\bar{L}$ is a functorial fibrant replacement in the category of $Q$-algebras. 
We assume without loss of generality that $\bar{L}Q$ concides with the fibrant relplacement chosen for the definition
of $\pi_{*,*}Q$ in the proof of Theorem \ref{ca2} (otherwise we would have to work with a quasi-isomorphic
$E_\infty$-quadratic pair algebra).
Then $\pi_{*,*}R$ is defined by the secondary homotopy groups of the spectrum $\bar{L}R$, as defined in Section
\ref{14}.
The unit $\bar{u}\colon Q\r R$ of the $Q$-algebra $R$ and the product in $\bar{L}R$ induce maps
$$\begin{array}{c}
\mu_1\colon \bar{L}Q\wedge_Q\bar{L}R\st{\bar{L}\bar{u}\wedge_Q1}\To \bar{L}R\wedge_Q\bar{L}R\st{\text{mult.}}\To \bar{L}R,\\
\mu_2\colon \bar{L}R\wedge_Q\bar{L}Q\st{1\wedge_Q\bar{L}\bar{u}}\To \bar{L}R\wedge_Q\bar{L}R\st{\text{mult.}}\To \bar{L}R,
\end{array}$$
such that there is a track $\alpha_1\colon \mu_1\tau_{\wedge_Q}\rr\mu_2$ which can be constructed as in Lemma
\ref{haxa} below. The external cup-one product
$$\smile_1\colon\pi_{*,0}R\times\pi_{*,0}Q\To \pi_{*,1}R$$
corresponds to the track $\alpha_1$, and the other external cup-one product corresponds to
$\alpha_1\tau_{\wedge_Q}$, compare the proof of Theorem \ref{ca2}.

If $L$ is the fibrant replacement functor in the category of ring spectra used to define $\pi_{*,*}$ in Theorems
\ref{a2} and \ref{ca2} then the natural quasi-isomorphisms making commutative the first square in Theorem \ref{ae2}
are obtained by taking secondary homotopy groups on
$$\bar{L}R\rightarrowtail L\bar{L}R\leftarrow LR.$$
For the second square one can argue in a similar way.

With the notation used for the definition of $\theta_{(x,y)}(a)$ in (\ref{teta}), the
elements $x$ and $y$ lie in $\pi_{n,0}Q$ and $\pi_{n,1}R$, respectively, and satisfy $\partial(y)=u_0(x)$ as
required in  Definition \ref{deride}, therefore the last part of the statement of Theorem \ref{ae2} makes sense.
\end{rem}

\section{Universal Toda brackets of ring spectra}\label{10}

Given a ring spectrum $R$ the category $\C{mod}(R)$ of right $R$-modules is a Quillen model category, see \cite[Theorem 12.1]{mcds}. 
By \cite{ah} the full
subcategory $\C{mod}^\fc(R)$ of fibrant-cofibrant right $R$-modules is a groupoid-enriched category whose morphisms are homotopy classes
of homotopies, also termed \emph{tracks}. 

A \emph{free right $R$-module} is a right $R$-module of the form
$$R^{(n_1,\dots,n_k)}\;=\;(S^{n_1}\vee\cdots\vee S^{n_k})\wedge R,\;\;k\geq0, \;\; n_1,\dots,n_k\geq0.$$
Here $\vee$ denotes the coproduct. Free right $R$-modules are cofibrant objects in
$\C{mod}(R)$ by \cite[Theorem 12.1 (vi)]{mcds} since sphere spectra $S^m$ are cofibrant in the category of
symmetric spectra.

We denote by $\C{mod}_f(R)$ the full groupoid-enriched subcategory of $\C{mod}^\fc(R)$ consisting of fibrant replacements of free right
$R$-modules. An object $X$ in $\C{mod}_f(R)$ is a trivial
cofibration $R^{(n_1,\dots,n_k)}\st{\sim}\rightarrowtail X$ in $\C{mod}(R)$ with fibrant target. This groupoid-enriched category has a
further structure described in the following definition.

\begin{defn}\label{gttc}
Let $\C{A}$ be a category and
$D\colon\C{A}^\op\times\C{A}\r\C{Ab}$ a functor, also called \emph{$\C{A}$-bimodule}. A
\emph{linear track extension} of $\C{A}$ by $D$ 
consists of
\begin{enumerate}
\item a groupoid-enriched category $\C{B}$.

\item A functor $p\colon\C{B}\r\C{A}$ from the underlying ordinary category of $\C{B}$
such that $p$ is the identity on objects and $p(f)=p(g)$ if and only
if there exists a track $f\rr g$, i.e. $\C{A}$ is the
\emph{homotopy category} of $\C{B}$.

\item Isomorphisms
$$\sigma_f\colon D(X,Y)\cong\aut_{\C{B}}(f)$$
for all morphisms $f\colon X\r Y$ in $\C{B}$ such that given a track
$\alpha\colon f\rr g$ in $\C{B}$,
\begin{equation*}
\alpha\vc\sigma_f(x)=\sigma_g(x)\vc\alpha;
\end{equation*}
and given composable morphisms
$\bullet\st{h}\r\bullet\st{g}\r\bullet\st{f}\r\bullet$ in $\C{B}$
\begin{equation*}
f\sigma_g(x)=\sigma_{fg}(D(1,p(f))x) \text{ and }
\sigma_g(x)h=\sigma_{gh}(D(p(h),1)x).
\end{equation*}
\end{enumerate}
\end{defn}

For any category $\C{A}$ and any functor $D\colon \C{A}^\op\times\C{A}\r\C{Ab}$ the
third Baues-Wirsching 
cohomology group $H^3(\C{A},D)$ (which in this case is isomorphic to a Hochschild-Mitchell cohomology) classifies
linear track extensions of $\C{A}$ by $D$,  see \cite{ccglg}.
Indeed any such linear track extension $\C{B}$ determines a characterisic cohomology class 
$$\grupo{\C{B}}\in H^3(\C{A},D)$$
also termed \emph{universal Toda bracket}, and all cohomology classes are represented by linear track extensions.

Given an $\N$-graded ring $A$ a \emph{free right $A$-module} is a right $A$-module of the form
$$A^{(n_1,\dots,n_k)}\;=\;((\Z,{n_1})\oplus\cdots\oplus(\Z,{n_k}))\otimes A,\;\;k\geq0, \;\; n_1,\dots,n_k\geq0.$$
Here $(\Z,m)$ denotes the $\N$-graded abelian group consisting of $\Z$ concentrated in degree $m$. The category
of free right $A$-modules is denoted by $\C{mod}_f(A)$.

Given an $A$-bimodule $M$ we define the graded version of Mac Lane cohomology by 
\begin{eqnarray*}
HML^*(A,M)&=&H^*(\C{mod}_f(A),\hom_A(-,-\otimes_AM)).
\end{eqnarray*}
For an ungraded ring $A$ if $\C{mod}_f(A)$ denotes the usual category of finitely generated free ungraded
right $A$-modules this equality corresponds to the natural isomorphism established in \cite{mlhthh} in homology.

\begin{lem}
Let $R$ be a connective spectrum. 
The category $\C{mod}_f(R)$ is a linear track extension of $\C{mod}_f(\pi_*R)$ 
by $\hom_{\pi_*R}(-,-\otimes_{\pi_*R}\S^{-1}\pi_*R)$.
\end{lem}

\begin{proof}
For any right $R$-module $Y$ the right $R$-module morphisms
$$R^{(n_1,\dots,n_k)}\To Y$$
coincide with the symmetric spectra morphisms
$$S^{n_1}\vee\cdots\vee S^{n_k}\To Y$$
hence it becomes obvious that $\C{mod}_f(\pi_*R)$ is the homotopy category of $\C{mod}_f(R)$ 
and the projection functor is given by taking homotopy groups $p=\pi_*$.

The group of self-tracks of a map $X\r Y$ between fibrant-cofibrant right $R$-modules coincides with 
the abelian group of homotopy classes $[\S X,Y]_R$ in the 
category of right $R$-modules, compare \cite{ccglg}. If $X$ and $Y$ are weakly equivalent to free right $R$-modules then 
this abelian group is canonically isomorphic to 
$\hom_{\pi_*R}(\pi_*X,\pi_*Y\otimes_{\pi_*R}\S^{-1}\pi_*R)$ by the same argument as before. The axioms of a linear track extension are checked
in \cite{ccglg} in a more general setting.
\end{proof}

The lemma above yields a universal Toda bracket associated to a connective ring spectrum
$R$
$$\grupo{\C{mod}_f(R)}\in HML^*(\pi_*R,\S^{-1}\pi_*R).$$
This universal Toda bracket determines all Toda brackets in $\C{mod}_f(R)$, see \cite{ccglg}, i.e. all matric Toda brackets in $R$,
compare \cite{mmp}. It has also been studied by Sagave from a different perspective, see \cite{steffen}.
The goal of this section is the construction of a small algebraic model for this universal Toda bracket. For this we will
consider universal Toda brackets associated to quadratic pair algebras.

Given a quadratic pair algebra $B$, in order to obtain a convenient groupoid-enriched category of 
right $B$-modules we need to consider weak morphisms.

\begin{defn}
A \emph{weak morphism} $f\colon C\r D$ between quadratic pair modules consists of group homomorphisms 
$$\begin{array}{c}
f_{0}\colon C_{0}\To C_{0},\\
f_{1}\colon C_{1}\To D_{1},\\
f_{ee}\colon C_{ee}\To D_{ee},
\end{array}$$
such that 
\begin{eqnarray*}
f_{0}\partial&=&\partial f_{1},\\
f_{1} P&=&P f_{ee},\\
f_{ee}(-|-)_H&=&(f_{0}|f_{0})_H,\\
f_{ee}T&=&Tf_{ee}.
\end{eqnarray*}

Weak morphisms were introduced in \cite{2hg2}. The functors $h_0$ and $h_1$ in (\ref{hi}) obviously extend to
the category $\C{wqpm}$ of quadratic pair modules and weak morphisms.

A weak morphism $f\colon M\r N$ between
right $B$-modules is a collection of weak morphisms $f_n\colon M_{n,*}\r N_{n,*}$, $n\geq 0$, which are strictly compatible
with the action of $B$, i.e. given $m_i\in M_{p,i}$ and $b_i\in B_{q,i}$ for $i=0,1,ee$,
\begin{eqnarray*}
f_{p+q,0}(m_0\cdot b_0)&=&f_{p,0}(m_0)\cdot b_0,\\
f_{p+q,1}(m_0\cdot b_1)&=&f_{p,0}(m_0)\cdot b_1,\\
f_{p+q,1}(m_1\cdot b_0)&=&f_{p,1}(m_1)\cdot b_0,\\
f_{p+q,ee}(m_{ee}\cdot b_{ee})&=&f_{p,ee}(m_{ee})\cdot b_{ee}.
\end{eqnarray*}

The reader can easily check that right $B$-modules together
with weak morphisms form a category that we denote by $\C{wmod}(B)$. The ordinary category of right $B$-modules
$\C{mod}(B)$ is a
subcategory of $\C{wmod}(B)$ with the same objects.
\end{defn}

The category $\C{wmod}(B)$ is also a groupoid-enriched category
with the tracks (i.e. $2$-morphisms) defined as follows.

\begin{defn}\label{altra}
Let $f,g\colon C\r D$ be weak morphisms between quadratic pair modules.
A \emph{track} $\alpha\colon f\rr g$ is a function $\alpha\colon C_0\r D_1$
satisfying the equations, $x,y\in C_0$, $z\in C_1$, 
\begin{eqnarray*}
\alpha(x+y)&=&\alpha(x)+\alpha(y)+P(-f_0(x)+g_0(x)|f_0(y))_H,\\
g_0(x)&=&f_0(x)+\partial\alpha(x),\\
g_1(z)&=&f_1(z)+\alpha\partial(z).
\end{eqnarray*}

This definition of track was given in \cite{2hg2}, where we show that the category $\C{wqpm}$ of quadratic pair
modules and weak morphisms is a groupoid-enriched category. The subcategory $\C{qpm}$ is also groupoid-enriched.
Tracks in $\C{qpm}$ are just tracks in $\C{wqpm}$ between morphisms in $\C{qpm}$.

The vertical composition of tracks 
$$\bullet\st{\alpha}\Longrightarrow\bullet\st{\beta}\Longrightarrow\bullet$$ 
is defined by addition
$(\beta\vc\alpha)(x)=\alpha(x)+\beta(x)$. The horizontal
composition of tracks and maps as in 
$$\xymatrix{\bullet\ar[r]^g&\bullet\ar@/^15pt/[r]_{\;}="a"\ar@/_15pt/[r]^{\;}="b"&\bullet\ar[r]^{f}&\bullet
\ar@{=>}"a";"b"^\alpha}$$
is defined as $(f\alpha)(x)=f_1\alpha(x)$ and $(\alpha
g)(x)=\alpha g_0(x)$. The identity track $0^\vc_f\colon f\rr f$ (also called trivial track) is always
defined as $(0^\vc_f)(x)=0$.

Let $B$ be a quadratic pair algebra and let $f,g\colon M\r N$ be weak right $B$-module morphisms. 
A \emph{track} $\alpha\colon f\rr g$ in $\C{mod}(B)$ is collection of tracks $\alpha_n\colon f_n\rr g_n$ between
the quadratic pair module
morphisms $f_n, g_n\colon M_n\r N_n$ such that for any $x\in M_{n,0}$ and $a\in B_{m,0}$
\begin{eqnarray*}
\alpha_{n+m}(x\cdot a)&=&\alpha_n(x)\cdot a+P((-f_{n,0}(x)+g_{n,0}(x)|f_{n,0}(x))_H\cdot H(a)).
\end{eqnarray*}
We leave it to the reader to check that these tracks with vertical and horizontal composition defined as above endow
$\C{wmod}(B)$ with a groupoid-enriched category structure. 
\end{defn}

\begin{rem}
The groupoid-enriched category $\C{qpm}$ was studied in \cite{2hg2, 2hg3}. The enrichment extends in the obvious way to the
category of $\N$-graded objects $\C{qpm}^\N$ and restricts to the subcategory of $0$-good objects $X$ in $\C{qpm}_H^\N$
which is monoidal. In this monoidal category tracks can be interpreted as morphisms from the cylinder object
$\mathbb{I}\odot X$ defined by the
interval quadratic pair module $\mathbb{I}$, see \cite[Lemma 5.3]{2hg3}. If $B$ is a $0$-good quadratic pair algebra
the cylinder of a $0$-good right $B$-module is again a right $B$-module in a natural way, so we can define tracks as right $B$-module morphisms from the cylinder. In this way we obtained the formulas for tracks in
Definition \ref{altra} above. This also works for weak morphisms.
\end{rem}

Let $B$ be a 
quadratic pair algebra. A \emph{free right $B$-module} is a $B$ module of the form
$$B^{(n_1,\dots,n_k)}\;=\;((\overline{\Z}_\nill,n_1)\vee\cdots\vee(\overline{\Z}_\nill,n_k))\odot B,\;\;k\geq0,\;\;n_1,\dots,n_k\geq0.$$
Here $(\overline{\Z}_\nill,m)$ denotes the $\N$-graded quadratic pair module given by $\overline{\Z}_\nill$
concentrated in degree $m$ and $\vee$ is the coproduct. 
We denote by $\C{wmod}_f(B)$ the full groupoid-enriched subcategory of $\C{wmod}(B)$ consisting of free
right $B$-modules.

\begin{lem}\label{lte2}
Let $B$ be a 
quadratic pair algebra. 
Then the category $\C{wmod}_f(B)$ is a linear track extension of 
$\C{mod}_f(h_0B)$ by $\hom_{h_0B}(-,-\otimes_{h_0B}h_1B)$.
\end{lem}

The proof of this lemma will be given below.

Lemma \ref{lte2} yields a universal Toda bracket in Mac Lane cohomology for any quadratic pair algebra $B$
$$\grupo{\C{wmod}_f(B)}\in HML^*(h_0B,h_1B).$$

The following theorem is the main result of this section.

\begin{thm}\label{mtb}
Let $R$ be a connective ring spectrum. If we use the isomorphisms in Theorem \ref{a2}
as identifications then the universal Toda brackets associated to the ring spectrum $R$ and to the quadratic pair
algebra $\pi_{*,*}R$ coincide, that is:
$$\grupo{\C{mod}_f(R)}\;=\;\grupo{\C{wmod}_f(\pi_{*,*}R)}.$$
\end{thm}

We prove Theorem \ref{mtb} at the end of this section.

\begin{rem}
In \cite[Theorem 2.2.1]{3mlc} it is proved in the ungraded setting that any element in Mac Lane cohomology
$HML^3(A,M)$ is represented by a quadratic pair algebra $B$ with $h_0B\cong A$ and $h_1B\cong M$. Such a
quadratic pair algebra is also termed a ``crossed extension of $A$ by $M$ in the symmetric monoidal category of
square groups'', see the Addendum of \cite[Theorem 2.2.1]{3mlc}. Since for a connective ring spectrum $R$ the
topologically defined category $\C{mod}_f(R)$ represents an element in $HML^3(\pi_*R,\Sigma^{-1}\pi_*R)$ the
result in \cite{3mlc} makes it plausible that, in fact, the element $\grupo{\C{mod}_f(R)}$ should also be
represented by a quadratic pair algebra. By Theorem \ref{mtb} above we construct such a quadratic pair algebra, given
by $\pi_{*,*}R$, even in a functorial way.
\end{rem}

In order to prove Lemma \ref{lte2} we need the following results.

\begin{lem}\label{calco}
There are natural isomorphisms
\begin{eqnarray*}
h_0(B^{(n_1,\dots,n_k)})&\cong&(h_0B)^{(n_1,\dots,n_k)},\\
h_1(B^{(n_1,\dots,n_k)})&\cong&(h_0B)^{(n_1,\dots,n_k)}\otimes_{h_0B}h_1B.\\
\end{eqnarray*}
\end{lem}

Lemma \ref{calco} follows easily from the following one.

\begin{lem}
For any quadratic pair module $C$ there are natural isomorphisms
\begin{eqnarray*}
h_0((\ol{\Z}_\nill\vee\st{k}\cdots\vee\ol{\Z}_\nill)\odot C)&\cong&h_0C\oplus\st{k}\cdots\oplus h_0C,\\
h_1((\ol{\Z}_\nill\vee\st{k}\cdots\vee\ol{\Z}_\nill)\odot C)&\cong&h_1C\oplus\st{k}\cdots\oplus h_1C.
\end{eqnarray*}
\end{lem}

\begin{proof}
There is an exact sequence of square groups in the sense of \cite{qaI}
$$h_1C\hookrightarrow C_{(1)}\st{\partial}\To C_{(0)}\twoheadrightarrow h_0C.$$

Given a pointed set $E$ the square group $\Z_\nill[E]$ in (\ref{znil}) is the coproduct of copies of the unit element $\Z_\nill$ for the
tensor product of square groups indexed by $E$ excluding the base-point. Therefore by using the definition in
\cite{2hg3} of
the tensor product in $\C{qpm}$ we see that the quadratic pair module
$$(\ol{\Z}_\nill\vee\st{k}\cdots\vee\ol{\Z}_\nill)\odot C$$
coincides with
$$\Z_\nill[\ul{k}]\odot C_{(1)}\st{\Z_\nill[\ul{k}]\odot\partial}\To \Z_\nill[\ul{k}]\odot C_{(0)}$$
where $\ul{k}=\set{0,1,\dots,k}$ pointed at $0$. 

The square groups in (\ref{znil}) are flat for the tensor product of square groups. This follows easily from
\cite[Corollary 39]{qaI} since the unit object $\Z_\nill$ is obviously flat. Moreover, it is easy to check by using the generatos and
relations for the tensor product of square groups in \cite{qaI} that for any abelian group $A$ there is a natural isomorphism $\Z_\nill[\ul{k}]\odot
A\cong\Z[\ul{k}]\otimes A$, and hence the lemma follows.
\end{proof}

The following lemma is a simple exercise.

\begin{lem}\label{bij}
Let $C$ be any quadratic pair module. Weak morphisms $f\colon \ol{\Z}_\nill\r C$ are in bijective correspondence
with the elements of $C_0$. The correspondence sends $f$ to $f_0(1)$ for $1\in\Z=(\ol{\Z}_\nill)_0$. Moreover, 
given two morphisms $f,g\colon \ol{\Z}_\nill\r C$ the set of tracks $\alpha\colon f\rr g$ is in bijection with the subset
$\partial^{-1}(-f_0(1)+g_0(1))\subset C_1$. The bijection sends $\alpha$ to $\alpha(1)$.
\end{lem}

Now we are ready to prove Lemma \ref{lte2}.

\begin{proof}[Proof of Lemma \ref{lte2}]
Weak right $B$-module morphisms $$B^{(n_1,\dots,n_k)}\To B^{(m_1,\dots,m_l)}$$ and tracks between them are in
bijection with weak morphisms between graded quadratic pair modules
$$(\overline{\Z}_\nill,n_1)\vee\cdots\vee(\overline{\Z}_\nill,n_k)\To B^{(m_1,\dots,m_l)}$$
and tracks between them.
Hence Lemma \ref{lte2} follows from Lemmas \ref{calco} and \ref{bij}. The projection onto the homotopy category is
$h_0\colon \C{wmod}_f(B)\r\C{mod}_f(h_0B)$.
\end{proof}

We finally need the following result.

\begin{lem}\label{pepitio}
Let $R$ be a connective ring spectrum. 
The restriction of the functor
$\pi_{*,*}$ in Theorem \ref{m2} to the full subcategory of connective fibrant-cofibrant right $R$-modules can be
extended to a groupoid-enriched functor.
\end{lem}

\begin{proof}
The functor $\pi_{*,*}$ in Theorem \ref{m2} is constructued in Section \ref{pt1} by using secondary homotopy
groups of spectra. Secondary homotopy groups of spectra are defined as colimits of additive secondary homotopy groups of
spaces which are groupoid-enriched functors, see \cite{2hg1, 2hg2}. Therefore secondary homotopy groups of
spectra take homotopies to tracks between quadratic pair modules, and in
particular the functor $\pi_{*,*}$ in Theorem \ref{m2} becomes a groupoid-enriched functor on connective fibrant-cofibrant right $R$-modules.   
\end{proof}

Now we are ready to prove Theorem \ref{mtb}.

\begin{proof}[Proof of Theorem \ref{mtb}]
Recall that a \emph{pseudofunctor} between groupoid-enriched categories
$\varphi\colon\C{C}\r\C{B}$ is an assignment of objects, maps and
tracks which preserves horizontal composition and identity morphisms only up to certain
given tracks
$$\varphi_{f,g}\colon\varphi(f)\varphi(g)\rr\varphi(fg)\text{ and
} \varphi_X\colon\varphi(1_X)\rr1_{\varphi(X)}.$$ These tracks must
satisfy well-known coherence and naturality properties. Track functors are pseudofunctors where these tracks are
trivial. Similarly there is a definition of \emph{pseudonatural} transformation $\alpha\colon\varphi\rr\psi$
between pseudofunctors $\varphi,\psi\colon\C{C}\r\C{B}$ given by maps $\alpha_X\colon\varphi(X)\r\psi(X)$ such
that the usual square does not commute, but for any map $f\colon X\r Y$ there is a track $\alpha_f\colon
\psi(f)\alpha_X\rr\alpha_Y\varphi(f)$ satisfying coherence and naturality properties.

Two linear track extensions $\C{B}$, $\C{C}$ of $\C{A}$ by $D$ have the same universal Toda bracket if and only
if there is a pseudofunctor $\varphi\colon\C{C}\r\C{B}$ such that for any morphism $f\colon X\r Y$ in $\C{C}$ the
equality $p_\C{B}\varphi(f)=p_\C{C}(f)$ holds, and for any $x\in D(X,Y)$ the equation
$\sigma_{\varphi(f)}(x)=\varphi\sigma_f(x)$ is satisfied, compare \cite{ptc}.

We now prove Theorem \ref{mtb} by constructing an appropriate pseudofunctor
$$\varphi\colon\C{mod}_f(R)\To \C{wmod}_f(\pi_{*,*}R).$$
The functor sending a fibrant replacement of a free right $R$-module $X$ to $\pi_{*,*}X$, see Theorem \ref{m2}, takes values in the category of right $\pi_{*,*}R$-modules, but not on the full subcategory 
$\C{wmod}_f(\pi_{*,*}R)$. However, the fixed trivial cofibration $R^{(n_1,\dots,n_k)}\st{\sim}\rightarrowtail X$ 
together with Theorem \ref{estatb} and Remark \ref{pinSn} give rise to a right $\pi_{*,*}R$-module morphism
\begin{equation}\label{explia}
{\alpha_X}\colon \varphi(X)=(\pi_{*,*}R)^{(n_1,\dots,n_k)}\To\pi_{*,*}X,
\end{equation}
see Remark \ref{expli} for details,
which induces isomorphisms in $h_0$ and $h_1$, compare Lemma \ref{calco}. Now it is easy to check that for any map $f\colon X\r Y$ between free right $R$-modules we can choose a morphism $\varphi(f)$ and a track $\alpha_f$ as in diagram
\begin{equation*}
\xymatrix{\varphi(X)\ar[r]^{\alpha_X}_{\;}="a"\ar[d]^{\;}="b"_{\varphi(f)}&\pi_{*,*}X\ar[d]^{\pi_{*,*}f}\\
\varphi(Y)\ar[r]^{\alpha_Y}&\pi_{*,*}Y\ar@{=>}"a";"b"^{\alpha_f}}
\end{equation*}
Since $\pi_{*,*}$ is  groupoid-enriched functor by Lemma \ref{pepitio} and $\alpha_X$ and $\alpha_Y$ induce isomorphisms on $h_0$ and $h_1$ then by
categorical nonsense
there is a unique pseudofunctor $\varphi$ defined on objects and maps as above such that 
$\alpha\colon \varphi\rr\pi_{*,*}$ is a pseudonatural transformation. 
This pseudofunctor proves the theorem.
\end{proof}

\section{Symmetric spectra}\label{sss}\label{11}

We work in the symmetric monoidal model category $\C{Top}^*$ of pointed (compactly generated)
spaces, see \cite[4.2.12]{hmc}. 

A \emph{spectrum} $X$ is a sequence of pointed spaces
$X_0$, $X_1$, \dots, $X_n$, \dots together with pointed maps
$\sigma_n\colon S^1\wedge X_n\r X_{1+n}$ called \emph{structure maps}. A
morphism of spectra $f\colon X\r Y$ is a sequence of pointed maps
$f_n\colon X_n\r Y_n$ with $\sigma_n f_n=f_{1+n}\sigma_n$. The
category of spectra will be denoted by $\C{Sp}^\N$.

An \emph{$\Omega$-spectrum} $X$ is a spectrum such that the adjoints
$ad(\sigma_n)\colon X_n\r\L X_{1+n}$ of the structure maps
$\sigma_n$ are weak equivalences for all $n\geq 0$.

A \emph{symmetric spectrum} is a more structured notion of spectrum defined in the following way. An object $X$ in
the category $\C{(\C{Top}^*)}^{\symm}$ of \emph{symmetric sequences} in $\C{Top}^*$ is a sequence of spaces
$X_0$, $X_1$, \dots, $X_n$,\dots such that $X_n$ carries a left action of the symmetric group $\sym{n}$. A morphism $X\r Y$ of
symmetric sequences is a sequence of equivariant maps $X_n\r Y_n$. The category of symmetric
sequences has a symmetric monoidal
structure with the tensor product $X\otimes Y$ characterized by the existence of a natural
isomorphism
\begin{eqnarray*}
\hom_{(\C{Top}^*)^{\symm}}(X\otimes Y,Z)&\cong&\prod_{p,q\geq 0}\hom_{\times}(X_p\wedge
Y_q,Z_{p+q})\colon f\mapsto f_{p,q}.
\end{eqnarray*}
Here $\hom_{\times}$ is the set of equivariant morphisms with respect to the cross product homomorphism $\times$
in (\ref{cp}).
The symmeytry isomorphism $X\otimes Y\cong Y\otimes X$ is given by the maps
$$X_p\wedge Y_q\cong Y_q\wedge X_p\st{(1_{Y\otimes X})_{q,p}}\To (Y\wedge X)_{q+p}\st{\tau_{q,p}}\To (Y\wedge
X)_{p+q}.$$
The shuffle permutation $\tau_{q,p}$ was already considered in Definition \ref{eqpa}.

The sphere spectrum $S$ is the symmetric sequence $S$ given by the spheres $S^0$, $S^1$, \dots, $S^n$, \dots, where $\sym{n}$ acts on
$S^n=S^1\wedge\cdots \wedge S^1$ by permutation of coordinates. This symmetric sequence is a commutative monoid
with multiplication $\mu\colon S\otimes S\r S$ given by $\mu_{p,q}\colon S^p\wedge S^q=S^{p+q}$.
The category $\C{Sp}^\S$ of symmetric
spectra is the category of left $S$-modules. The structure maps $\sigma_{1,n}\colon S^1\wedge X_n\r X_{1+n}$ of 
the underlying spectrum of a symmetric spectrum $X$ are defined by the left action $\sigma\colon S\otimes X\r X$.

Since the monoid $S$ is commutative the category $\C{Sp}^\S$ is symmetric monoidal for the tensor product $\otimes_S$ also
called \emph{smash product} and
denoted by $\wedge$. The smash product of symmetric spectra $X\wedge Y$ is the coequalizer of the two 
multiplications $X\otimes S\otimes Y\rightrightarrows X\otimes Y$. The smash product $X\wedge Y$ comes equipped
with natural maps, $p,q\geq 0$,
\begin{equation}\label{jmath}
\jmath_{p,q}\colon X_p\wedge Y_q\To (X\wedge Y)_{p+q}
\end{equation}
satisfying a universal property which characterizes $X\wedge Y$ as a symmetric spectrum, see \cite{mtadsmc}.

The category $\C{Sp}^\N$ of spectra is a stable model category, see \cite[4.1 and 9]{mcds}, while $\C{Sp}^\S$ is in addition a symmetric monoidal
model category, see \cite[4.2 and 9]{mcds}. We wish to emphasize that among the possible model category
structures on $\C{Sp}^\N$ and $\C{Sp}^\S$ we have chosen to work with the \emph{stable model structure} defined in
\cite[9]{mcds} since with this structure fibrant objects coincide in both cases with the $\L$-spectra.

Symmetric sequences and symmetric spectra
defined in this way are
available over monoidal categories more general than $\C{Top}^*$, see \cite{sssgmc}. They 
can also be defined by using
right actions of symmetric groups instead of left actions. 

\section{Secondary homotopy groups of spaces}\label{12}

A groupoid-enriched category, also termed \emph{track category}, is a $2$-category where all $2$-morphisms are invertible with
respect to vertical composition. We also call \emph{tracks} to the $2$-morphisms in track
categories. Horizontal composition is
denoted by juxtaposition, and we use the symbol $\vc$ for the
vertical composition. The vertical inverse of a track $\alpha$ is $\alpha^\vi$. Identity morphisms are denoted by $1$, and the
symbol $0^\vc$ is used for identity tracks. \emph{Track functors} are $2$-functors between track categories. 
The category $\C{Top}^*$ is a track category where a track between two maps is a relative homotopy class of
homotopies between them. Similarly for the category of fibrant-cofibrant spectra or symmetric spectra, compare \cite{ah}.
For the convenience of the reader we recall from \cite{2hg1}, \cite{2hg2} the following definition.

\begin{defn}
Let $n\geq 3$. For a pointed space $X$ we define the \emph{additive
secondary homotopy group}, $\Pi_{n,*}X$ which is a $0$-free quadratic pair
module with $\Pi_{n,(0)}X=\Z_\nill[\L^nX]$ where $\L^nX$ is the discrete pointed set of maps $S^n\r X$ in $\C{Top}^*$.

We describe the group $\Pi_{n,1}X$ and the homomorphisms $P$ and
$\partial$ as follows. Given a pointed set $E$ we denote by $\vee_ES^n$ to the $n$-fold suspension of $E$ in
$\C{Top}^*$. It is a one-point union of $n$-spheres indexed by the set $E-\set{*}$. The group $\Pi_{n,1}X$  is given by the set of equivalence classes $[f,F]$ 
represented
by a map $f\colon S^1\r \vee_{\L^n X}S^1$ and a track in $\C{Top}^*$ of the form:
$$\xymatrix{S^n\ar[r]_{\S^{n-1}f}^<(.98){\;\;\;\;\;}="a"\ar@/^25pt/[rr]^0_{}="b"&\vee_{\L^nX}S^n\ar[r]_{ev}&X.\ar@{=>}"a";"b"_F}$$
Here the map $ev\colon \vee_{\L^nX}S^n\r X$
is the obvious evaluation map. 

The equivalence relation $[f,F]=[g,G]$ holds provided there is a 
track $$N\colon\Sigma^{n-1}f\rr\S^{n-1}g$$ with trivial Hopf invariant $\hopf(N)=0$ in the sense of (\ref{elhopf}) below such that 
the composite track in the following diagram is the trivial track.
\begin{equation}\label{algo}
\xymatrix@C=50pt{S^n\ar@/^40pt/[rr]^0_{\;}="a"\ar@/_40pt/[rr]_0^{\;}="f"\ar@/^15pt/[r]|{\S^{n-1}f}^<(.935){\;}="b"_{\;}="c"\ar@/_15pt/[r]|{\S^{n-1}g}^{\;}="d"_<(.93){\;}="e"
&\vee_{\L^nX}S^n\ar[r]^{ev}&X\ar@{=>}"a";"b"^{F^\vi}\ar@{=>}"c";"d"^{N}\ar@{=>}"e";"f"^G}
\end{equation}
The element $\partial[f,F]\in\grupo{\L^nX}_\nill$ is given by the image of the generator $1\in\Z\cong\pi_1S^1$
under the homomorphism
$$\pi_1f\colon\pi_1S^1\To\pi_1(\vee_{\L^nX}S^1)\cong\grupo{\L^nX}.$$

Let $I_+$ be the unit interval $I=[0,1]$ with an outer base-point and let $IX=I_+\wedge X$ be the reduced
cylinder of $X$ in $\C{Top}^*$. The \emph{Hopf invariant} $\hopf(N)$ of a track $N\colon\Sigma^{n-1}f\rr\S^{n-1}g$ with $f,g\colon
S^1\r\vee_ES^1$ is defined in \cite{2hg1} 3.3 by the
homomorphism 
\begin{equation}\label{elhopf}
H_2(IS^1,S^1\vee S^1)\st{ad(N)_*}\To H_2(\L^{n-1}(\vee_ES^n),\vee_ES^1)\cong
\hat{\otimes}^2\Z[E],
\end{equation}
which carries the generator $1\in\Z\cong H_2(IS^1,S^1\vee S^1)$ to $\hopf(N)$. Here the isomorphism is induced by
the Pontrjagin product and $ad(N)_*$ is the homomorphism
induced in homology by the adjoint of $$S^{n-1}\wedge I_+\wedge S^1\cong I_+\wedge S^n\st{N}\r \vee_ES^n.$$ The
reduced tensor square $\hat{\otimes}^2$ in (\ref{elhopf}) is defined for any abelian group $A$ as the quotient
$$\hat{\otimes}^2A=\frac{A\otimes A}{a\otimes b+b\otimes a\sim 0}$$ 
and $\bar{\sigma}\colon\otimes^2A\twoheadrightarrow\hat{\otimes}^2A$ denotes the natural projection.
If $f,g\colon
\vee_{E'}S^1\r\vee_ES^1$ are maps from a coproduct the Hopf invariant $\hopf(N)\colon\Z[E']\r\hat{\otimes}^2\Z[E]$
is a homomorphism defined on the basis $E'$ by the Hopf invariants of the restriction of $N$
to the components of $\vee_{E'}S^n$.
We refer the reader to \cite[3]{2hg1} for the elementary properties of the invariant $\hopf$.

The equivalence classes $[f,F]$ define  $\Pi_{n,1}X$ as a set. The
group structure of $\Pi_{n,1}X$ is induced by the comultiplication $\mu\colon S^1\r S^1\vee S^1$, compare
\cite{2hg1} 4.4.

We now define the homomorphism $P$. Consider the diagram
$$\xymatrix{S^n\ar[rr]_{\S^{n-1}\beta}^{\;}="a"\ar@/^30pt/[rr]^0_{\;}="b"&&S^n\vee S^n\ar@{=>}"a";"b"_{B}}$$
where $\beta\colon S^1\r S^1\vee S^1$ is given such that $(\pi_1\beta)(1)=-a-b+a+b$  
is the commutator of the canonical generators of $\pi_1(S^1\vee S^1)$, i.e. $a$ corresponds to the inclusion of
the first sphere and $b$ to the inclusion of the second one. The track $B$ 
is any track with
$\hopf(B)=-\bar{\sigma}(a\otimes b)$. 
Given $x\otimes y\in\otimes^2\Z[\L^nX]$ let
$\tilde{x},\tilde{y}\colon S^1\r\vee_{\L^nX}S^1$ be maps such that the images of
$(\pi_1\tilde{x})(1)$ and $(\pi_1\tilde{y})(1)$ in $\Z[\L^nX]$ are $x$ and
$y$, respectively. Then
the diagram
\begin{equation}\label{omedia}
\xymatrix{S^n\ar[rr]_{\S^{n-1}\beta}^{\;}="a"\ar@/^30pt/[rr]^0_{\;}="b"&&S^n\vee
S^n\ar@{=>}"a";"b"_{B}\ar[rr]_{\S^{n-1}(\tilde{y},\tilde{x})}&&S^n_X\ar[r]_{ev}&X}
\end{equation}
represents $P(x\otimes y)\in\Pi_{n,1}X$.
This completes the definition of the quadratic pair module $\Pi_{n,*}X$ for $n\geq 3$. 
\end{defn}

Additive secondary homotopy groups are defined in
\cite{2hg2} also for $n=0,1,2$. In this way we get for $n\geq 0$ a
functor
$$\Pi_{n,*}\colon\C{Top}^*\To\C{qpm}$$
which is actually a track functor. Moreover, $\Pi_{n,(0)}X=\Z_\nill[\L^nX]$ so that $\Pi_{n,*}$ takes values in the full
subcategory $\C{qpm}_H$ of $0$-free, and hence $0$-good, quadratic pair modules. 

\begin{rem}\label{PinSn}
By \cite[4.16]{2hg1} $\Pi_{n,*}S^n$ is quasi-isomorphic to $\ol{\Z}_\nill$. The quasi-isomorphism $\imath_{n,0}\colon
\ol{\Z}_\nill\r\Pi_{n,*}S^n$ is defined by $\imath_{n,0}(1)=1_{S^n}\in\grupo{\L^nS^n}_\nill$ for
$1\in\Z=(\ol{\Z}_\nill)_0$. This quasi-isomorphism is an isomorphism for $n=0$.
\end{rem}

The secondary homotopy group $\Pi_{n,*}X$ is endowed with a natural action of the symmetric track group
$\symt{n}$, see \cite{2hg2}. The underlying  action of $\sym{n}$ on
$\Pi_{n,(0)}X=\Z_\nill[\L^nX]$ is given by the right action of $\sym{n}$ on the pointed set $\L^nX$ of maps $S^n\r X$
determined by the left action of $\sym{n}$
on the $n$-sphere $S^n=S^1\wedge\cdots\wedge S^1$ by permutation of coordinates.

One of the results in \cite{2hg3} states that the functor
\begin{equation}\label{lm}
\Pi_{*,*}\colon\C{Top}^*\To\C{qpm}_H^\N
\end{equation}
to the category of $\N$-graded $0$-good quadratic pair modules is lax monoidal with multiplication given by the smash product operation for secondary homotopy groups
$$\Pi_{n,*}X\odot\Pi_{m,*}Y\st{\wedge}\To\Pi_{n+m,*}(X\wedge Y).$$ 
This morphism is equivariant with respect to the sign group morphism (\ref{sgx}). 
The functor $\Pi_{*,*}$ in (\ref{lm}) is however non-symmetric. 
In order to obtain a symmetric version we have to regard $\Pi_{*,*}$ as a functor to the category
$\C{qpm}_H^{\symtt}$ of enhanced
symmetric sequences of $0$-good quadratic pair modules defined in \cite[8.14]{2hg3}. 
The properties of the smash product for secondary homotopy groups yield a lax symmetric
monoidal functor
\begin{equation}\label{lsm}
\Pi_{*,*}\colon\C{Top}^*\To\C{qpm}_H^{\symtt}.
\end{equation}
This is the main result in \cite{2hg3}.

There are natural
isomorphisms, $n\geq 3$,
\begin{eqnarray}
\label{h0} h_0\Pi_{n,*}X&\cong&\pi_nX,\\	
\nonumber h_1\Pi_{n,*}X&\cong&\pi_{n+1}X.
\end{eqnarray}
Here we use \cite{2hg1} 5.1.
Furthermore, the following property is crucial.

\begin{prop}\label{qlo}
For all $n\geq 3$ the $k$-invariant of $\Pi_{n,*}X$ defined in (\ref{-w})  coincides via (\ref{h0}) with the
homomorphism $\eta^*\colon\pi_nX\r\pi_{n+1}X$ induced by precomposition with $\S^{n-2}\eta$
where $\eta\colon S^3\r S^2$ is the Hopf map. 
\end{prop}

This follows from \cite[8.2]{2hg1}.

\section{The pull-back construction}\label{13}

Let $C$ be any quadratic pair module. Given a pointed set $E$ and a pointed map $E\r \ker H\subset C_0$, or equivalently 
a group homomorphism $\varphi\colon \grupo{E}_\nill\r C_0$  such that $H\varphi(e)=0$ for any $e\in E$ then the $0$-free quadratic pair module $\varphi^*C$ is defined as
follows.
\begin{equation}\label{relative}
\xymatrix@C=40pt{\;\;\;\;\;\;\;\;\;\;\;\;\;\;\;\;\;(\varphi^*C)_{ee}=\otimes^2\Z[E]\ar@{-->}[d]_P\ar[r]^-{(\varphi_\abb|\varphi_\abb)_H}
\ar@/_40pt/[dd]_{P\text{ in (\ref{znil})}}&
C_{ee}\ar[d]^P\\
(\varphi^*C)_1\ar[d]_\partial\ar[r]\ar@{}[rd]|{\text{pull}}&C_1\ar[d]^\partial\\
\;\;\;\;\;\;\;\;\;\;\;(\varphi^*C)_0=\grupo{E}_\nill
\ar[r]_-\varphi&C_0}
\end{equation}
This quadratic pair module comes equipped with a morphism $\varphi_\#\colon \varphi^*C\r C$ in $\C{qpm}$ which is
given by the horizontal arrows in diagram
(\ref{relative}). This construction
can be straightforwardly extended to the graded setting degreewise.

\begin{lem}\label{h0iso}
In the conditions above the morphism $\varphi_\#\colon \varphi^*C\r C$ induces always an isomorphism on $h_1$. It also induces an isomorphism on $h_0$
provided the set of elements in
$h_0C$ comming from $E$ generate the group $h_0C$.
\end{lem}


\begin{lem}\label{actpull}
Suppose that $G_\vc$ is a sign group acting on a quadratic pair module $C$ and that the group $G$ acts on
$E$ in such a way that $\varphi$ is $G$-equivariant. 
Then there is a unique action of $G_\vc$ on $\varphi^*C$ 
such that the induced $G$-action on level $0$ is given by the action of $G$ on $E$ and the natural morphism
$\varphi_\#\colon\varphi^*C\r C$ is a right 
$A(G_\vc)$-module morphism.
\end{lem}

This follows easily from the universal property of a pull-back.

\begin{lem}
If $C$ is a quadratic pair algebra, $E$ is a pointed graded monoid, and $E\r \ker H\subset C_0$ is a graded monoid morphism then
there is a unique quadratic pair algebra structure on $\varphi^*C$ such that the monoid structure of
$(\varphi^*C)_{*,(0)}=\Z_\nill[E]$ in $\C{SG}$ is induced by the monoid structure on $E$ in the sense of
\cite[Section 12.1]{qaI} and the natural morphism
$\varphi_\#\colon\varphi^*C\r C$ is a morphism of quadratic pair algebras. 
\end{lem}

This is also a consequence of the universal property of a pull-back. In can be extended to the $E_\infty$-case.

\begin{lem}\label{kqt}
Suppose that $C$ is an $E_\infty$-quadratic pair algebra, $E$ is a pointed graded monoid 
such that the symmetric group $\sym{n}$  acts on $E_n$ in such a way that the multiplications $E_n\wedge E_m\r
E_{n,m}$ are equivariant with respect to the cross product homomorphisms in
(\ref{cp}), and $E\r \ker H$ is a degreewise equivariant 
graded monoid morphism. Then
there is a unique $E_\infty$-quadratic pair algebra structure on $\varphi^*C$ such that the monoid structure of
$(\varphi^*C)_{*,(0)}=\Z_\nill[E]$ in $\C{SG}$ is induced by the monoid structure on $E$ in the sense of
\cite[Section 12.1]{qaI}, the action of $\sym{n}$ on $(\varphi^*C)_{n,(0)}=\Z_\nill[E_n]$ is given by the
action on $E_n$, and the natural morphism
$\varphi_\#\colon\varphi^*C\r C$ is a morphism of $E_\infty$-quadratic pair algebras.
\end{lem}

\section{Secondary homotopy groups of spectra}\label{shgs}\label{14}

Secondary homotopy groups of spaces introduced in the previous section can be used to define secondary homotopy
groups of spectra which are the crucial tools in this paper.

Given a space $X$ and a pointed map $E\r \L^nX$ inducing a homomorphism $\varphi\colon
\grupo{E}_\nill\r\grupo{\L^nX}_\nill$, we write 
\begin{eqnarray*}
\Pi_{n,*}(X|E)&=&\varphi^*\Pi_{n,*}X
\end{eqnarray*}
for the pull-back construction defined in (\ref{relative}).

\begin{defn}\label{dshgs}
Let $X$ be now a spectrum. The \emph{secondary homotopy group} $\pi_{n,*}X$, $n\geq 0$, is defined as the colimit in $\C{qpm}$
$$\pi_{n,*}X=\colim_{k\geq0}\Pi_{k+n,*}(X_{k}|\L^n X_0).$$
Here the map $\L^nX_0\r \L^{k+n}X_{k}$ is obtained by taking $\L^n$ on the adjoint of the structure map $\sigma_{k,0}\colon S^k\wedge
X_0\r X_{k}$, i.e. it 
sends an $n$-loop $f\colon S^n\r X_0$ to the composite
$$S^{k+n}\st{S^k\wedge f}\To S^k\wedge X_0\st{\sigma_{k,0}}\To X_{k}.$$ The bonding morphisms of the
directed system
$$\Pi_{k+n,*}(X_{k}|\L^n X_0)\To\Pi_{1+k+n,*}(X_{1+k}|\L^n X_0)$$
are defined to be the identity on levels $0$ and $ee$, so $\pi_{n,(0)}X=\Z_\nill[\L^nX_0]$ and $\pi_{n,*}X$ is $0$-free, and hence $0$-good.
On the $1$-level the bonding morphism is defined by the composite in $\C{qpm}$
$$\xymatrix@R=20pt{\Pi_{k+n,*}(X_{k})\ar[d]^\cong\\
\overline{\Z}_\nill\odot\Pi_{k+n,*}(X_{k})\ar[d]^{\imath_{1,0}\odot1}\\
\Pi_{1,*}S^1\odot\Pi_{k+n,*}(X_{k})\ar[d]^{\wedge}\\
\Pi_{1+k+n,*}(S^1\wedge X_{k})\ar[d]^{\Pi_{1+k+n,*}\sigma_{1,k}}\\\Pi_{1+k+n,*}X_{1+k}}$$
Here the first isomorphism is the unit isomorphism for the symmetric monoidal structure $\odot$ in $\C{qpm}$. 
In the second arrow we use
the morphism $\imath_{1,0}$ in Remark \ref{PinSn}. 
The morphism $\wedge$ is the smash product for secondary homotopy groups of
spaces constructed in \cite{2hg3}.

\begin{rem}\label{pinSn}
The quadratic pair module $\pi_{n,*}S^n$ is quasi-isomorphic to $\ol{\Z}_\nill$. 
The quasi-isomorphism $\imath_{n,0}\colon
\ol{\Z}_\nill\r\pi_{n,*}S^n$ is defined by $\imath_{n,0}(1)=1_{S^n}\in\grupo{\L^nS^n}_\nill$, as in Remark
\ref{PinSn}. This quasi-isomorphism is an isomorphism for $n=0$. This remark follows for instance from Proposition \ref{qlo2} below.
\end{rem}

The symmetric track group $\symt{n}$ acts on $\Pi_{k+n,*}(X_{k}|\L^nX_0)$, see Lemma \ref{actpull}. The induced action of the symmetric
group $\sym{n}$ at level $0$ is given by the action on $\L^nX_0$ considered in the previous section. The bonding
morphisms of the directed system are right 
$A(\symt{n})$-module morphisms, therefore $\pi_{n,*}X$ carries a natural action of $\symt{n}$. 
\end{defn}



For any $\L$-spectrum $X$ we derive from (\ref{h0}) and Lemma \ref{h0iso}
the natural isomorphisms, $n\geq 0$,
\begin{eqnarray}
\label{h02} h_0\pi_{n,*}X&\cong&\pi_{n}X,\\	
\nonumber h_1\pi_{n,*}X&\cong&\pi_{n+1}X.
\end{eqnarray}
The following proposition follows from Proposition \ref{qlo}.

\begin{prop}\label{qlo2}
Let $X$ be an $\L$-spectrum. The $k$-invariant of $\pi_{n,*}X$ coincides via the isomorphisms (\ref{h02}) with the
homomorphism $\eta\colon\pi_{n}X\r\pi_{n+1}X$ induced by multiplication with the stable Hopf map $\eta\in\pi_1S$. 
\end{prop}

The projection to the colimit yields a natural morphism in $\C{qpm}$
$$\phi_n\colon\Pi_{n,*}X_0\To\pi_{n,*}X.$$
In case $X$ is an $\L$-spectrum $\phi_n$ is levelwise surjective for $n=0,1$ and an isomorphism for $n\geq 2$.
The next result is a crucial ingredient for the main theorems in this paper.

\begin{thm}\label{sur2}
Let $X, Y, Z$ be symmetric spectra which are $\L$-spectra and let $f\colon X\wedge Y\r Z$ be a morphism in $\C{Sp}^\S$.
Then for all $n,m\geq 0$ there is a unique morphism $\tilde{f}_{n,m}$ in $\C{qpm}$ for which the following diagram
commutes.
$$\xymatrix@C=50pt{\Pi_{n,*}X_0\odot\Pi_{m,*}Y_0\ar[r]^-\wedge
\ar@{->>}[d]_{\phi_n\otimes\phi_m}&
\Pi_{n+m,*}(X_0\wedge Y_0)\ar[r]^-{f_*}&\Pi_{n+m}Z_0\ar@{->>}[d]^{\phi_{n+m}}\\
\pi_{n,*}X\odot \pi_{m,*}Y\ar@{-->}[rr]_-{\tilde{f}_{n,m}}&&
\pi_{n+m,*}Z}$$
\end{thm}

The filler $\tilde{f}_{n,m}$ will be constructed in Theorem \ref{estatb}. In order to define it
we consider the morphisms
\begin{equation}\label{level}
\xymatrix{\Pi_{k+n,*}(X_{k}|\L^nX_0)\odot\Pi_{l+m,*}(Y_{l}|\L^mY_0)\ar[d]\\
\Pi_{k+l+n+m,*}((X\wedge Y)_{k+l}|\L^{n+m}(X\wedge Y)_{0})}
\end{equation}
defined on the $(0)$-level by
\begin{equation}\label{s0}
\xymatrix{\Z_\nill[\L^nX_0]\odot\Z_\nill[\L^mY_0]\ar[d]^\cong\\\Z_\nill[(\L^nX_0)\wedge(\L^mY_0)]
\ar[d]^{\Z_\nill[\wedge]}\\
\Z_\nill[\L^{n+m}(X_0\wedge Y_0)]
\ar[d]^{\Z_\nill[\L^{n+m}\jmath_{0,0}]}\\\Z_\nill[\L^{n+m}(X\wedge Y)_{0}]}
\end{equation}
Here the first arrow is the isomorphism in \cite[3,7]{2hg3} and $\wedge \colon(\L^nX_0)\wedge(\L^mY_0)\r\L^{n+m}(X_0\wedge Y_0)$ is the map 
defined by $(f\colon S^n\r X_0)\wedge (g\colon S^m\r Y_0)=f\wedge g\colon S^{n+m}\r X_0\wedge Y_0$. Moreover, the
morphism (\ref{level}) is induced on the $1$-level by
\begin{equation}\label{s1}
\xymatrix{\Pi_{k+n,*}(X_{k})\odot\Pi_{l+m,*}(Y_{l})\ar[d]^{\wedge}\\
\Pi_{k+n+l+m,*}(X_{k}\wedge Y_{l})\ar[d]^{\Pi_{k+n+l+m,*}\jmath_{k,l}}\\
\Pi_{k+n+l+m,*}(X\wedge Y)_{k+l}\ar[d]^{(1_k\times\tau_{l,n}\times1_m)^*}\\
\Pi_{k+l+n+m,*}(X\wedge Y)_{k+l}}
\end{equation}
Here the last arrow is given by the right action of the sign group $\symt{k+l+n+m}$.

\begin{lem}
The morphism (\ref{level}) is well defined.
\end{lem}

\begin{proof}
We have to check that the composite (\ref{s0}) is compatible with (\ref{s1}) so that (\ref{level})
is indeed defined. This follows from the fact that given maps $f\colon S^n\r X_0$, $g\colon S^m\r Y_0$ the
following diagram commutes
$$\xymatrix@C=45pt{S^k\wedge S^l\wedge S^n\wedge S^m\ar[r]^{1_k\times\tau_{l,n}\times1_m}
\ar[d]_{S^{k+l}\wedge f\wedge g}&
S^k\wedge S^n\wedge S^l\wedge S^m\ar[r]^{S^k\wedge f\wedge S^l\wedge g}&
S^k\wedge X_0\wedge S^l\wedge Y_0\ar[d]^{\sigma_{k,0}\wedge\sigma_{l,0}}\\
S^{k+l}\wedge X_0\wedge Y_0\ar[d]_{\jmath_{0,0}}&&X_{k}\wedge Y_{l}\ar[d]^{\jmath_{k,l}}\\
S^{k+l}\wedge (X\wedge Y)_{0}\ar[rr]_{\sigma_{k+l,0}}&&(X\wedge Y)_{k+l}}$$
\end{proof}

\begin{thm}\label{equit}
Suppose that $E$ is a pointed space, $\varphi\colon E\r\L^{n+m}X$ 
is a pointed weak equivalence and the symmetric group $\sym{n}$ acts from the left on
$X$ in such a way that the adjoint
$$\psi\colon S^n\wedge S^m\wedge E\To X$$
of $\varphi$ is equivariant with respect to the left action of $\sym{n}$ on the sphere $S^{n}$. Suppose also that $m\geq 1$. Then for any $\tau\in\sym{n}$ the
two possible composites in the diagram
$$\xymatrix@C=25pt{\Pi_{n+m,*}(X|E)\ar[r]^-{\varphi_\#}&\Pi_{n+m,*}X\ar@<.5ex>[rr]^{(\tau\times
1_m)^*}\ar@<-.5ex>[rr]_{\Pi_{n+m,*}\tau}&&\Pi_{n+m,*}X}$$
coincide.
\end{thm}

\begin{proof}
The theorem is obvious in degrees $0$ and $ee$. Let us check the degree $1$ case. Obviously the statement of Theorem \ref{equit} is
non-trivial just in case $n\geq 2$ so we can suppose that we are within this range, and hence $n+m\geq3$.

Let $[f,F]\in\Pi_{n+m,1}(X,E)$. This element is represented by a map
$$f\colon S^1\To\vee_ES^1,$$
such that the image of $(\pi_1f)(1)\in\grupo{E}$ in the quotient $\grupo{E}_\nill$ is $\partial[f,F]$, and a track
$$\xymatrix@C=30pt{S^{n+m}\ar[r]_-{\S^{n+m-1}f}\ar@/^35pt/[rrr]^0_{}="b"&
\vee_ES^{n+m}\ar[r]^{\;}="a"_-{\S^{n+m}\varphi}&\vee_{\L^{n+m}X}S^{n+m}\ar[r]_-{ev}&X.\ar@{=>}"a";"b"_F}$$
Here $\varphi$ is the underlying map of $\varphi$ above between pointed discrete sets.

On the one hand the element $(\Pi_{n+m,*}\tau)\varphi_\#[f,F]$ is represented by
\begin{equation}\label{a}
\xymatrix@C=33pt{&&&X\ar@/^10pt/[rd]^\tau&\\S^{n+m}\ar@/^15pt/[rrru]^0_{\;}="b"\ar[r]_-{\S^{n+m-1}f}&
\vee_ES^{n+m}\ar[r]_-{\S^n\varphi}&
\vee_{\L^{n+m}X}S^{n+m}\ar[ru]_{ev}^{\;}="a"
\ar[r]_-{\S^{n+m}\L^{n+m}\tau}&\vee_{\L^{n+m}X}S^{n+m}\ar[r]_-{ev}&X\ar@{=>}"a";"b"^F}
\end{equation}
Here $\L^{n+m}\tau$ is also regarded as a map between pointed discrete sets.

On the other hand $(\tau\times1_m)^*\varphi_\#[f,F]$ is given by
\begin{equation}\label{b}
\xymatrix@R=70pt@C=28pt{&&&&\\S^{n+m}\ar@/^100pt/[rrrrd]_{\;}="f"^0\ar[r]_-{\S^{n+m-1}f}&
\vee_ES^{n+m}\ar[rr]^{\;}="e"|<(.3){\S^{n+m}\varphi}_<(.128){\;}="d"&&**[l]\vee_{\L^{n+m}X}S^{n+m}\ar[rd]|-{ev}&\\
S^{n+m}\ar[r]_-{\S^{n+m-1}f}\ar[u]^{(\cdot)_{n+m}^{\sign(\tau)}}&
\vee_ES^{n+m}\ar[r]_-{\S^{n+m}\varphi}^{\;}="c"&
\vee_{\L^{n+m}X}S^{n+m}\ar[r]_<(.1){\S^{n+m}(\tau\times1)^*}
\ar@/^20pt/[ru]^<(.5){(\cdot)_{n+m}^{\sign(\tau)}\wedge\L^{n+m}X}_{\;}="b"
\ar@/_20pt/[ru]_{\tau\wedge\L^{n+m}X}^{\;}="a"
&\vee_{\L^{n+m}X}S^{n+m}\ar[r]_-{ev}&X\ar@{=>}"a";"b"|<(.4){\hat{\tau}\wedge\L^{n+m}X}
\ar@{=>}"c";"d"_Q\ar@{=>}"e";"f"_F}\
\end{equation}
Here $\L^{n+m}X$ is regarded as a discrete pointed set, $(\tau\times 1)^*\colon\L^{n+m}X\r\L^{n+m}X$ is the pointed map induced by precomposition with
$\tau\times1_m=\tau\wedge S^m\colon S^{n+m}\r S^{n+m}$, $(\cdot)_{n+m}^{\sign(\tau)}\colon S^{n+m}\r S^{n+m}$ is
the $(n+m-1)$-fold suspension of the map
$(\cdot)^{\sign(\tau)}\colon S^1\r S^1$ where we use the multiplicative 
topological abelian group structure of $S^1$, $\hat{\tau}\colon
\tau\times1_m\rr (\cdot)_{n+m}^{\sign(\tau)}$ is any track, and $Q$ is the unique track with Hopf invariant
\begin{eqnarray*}
\hopf(Q)&=&-(\hat{\otimes}^2\Z[\varphi])\binom{\sign(\tau)}{2}H\partial[f,F].
\end{eqnarray*}

Since $\psi$ is $\sym{n}$-equivariant we have the equation
\begin{equation}\label{e}
(\L^n\tau)\varphi=(\tau\times1_m)^*\varphi.
\end{equation}
We use this equation in order to simplify diagram (\ref{b}) by using Lemma \ref{kita} below.

The hypothesis $m\geq1$ is also needed in order to apply Lemma \ref{kita} below.
By Lemma \ref{kita} the track (\ref{b}) coincides with
\begin{equation}\label{c}
\xymatrix@R=60pt@C=27pt{&&&&\\
S^{n+m}\ar@/^100pt/[rrrrd]_{\;}="f"^0\ar[r]_-{\S^{n+m-1}f}_{\;}="d"&
**[r]\vee_ES^{n+m}\ar[rr]^{\;}="e"_{\S^{n+m}\varphi}&&
**[l]\vee_{\L^{n+m}X}S^{n+m}\ar[rd]|-{ev}&\\
S^{n+m}\ar[r]_-{\S^{n+m-1}f}^{\;}="c"\ar[u]_{\tau\times1_m}^{\;}="a"
\ar@/^30pt/[u]|<(.2){(\cdot)_{n+m}^{\sign(\tau)}}_{\;}="b"&
\vee_ES^{n+m}\ar[r]_-{\S^{n+m}\varphi}&
\vee_{\L^{n+m}X}S^{n+m}\ar[r]_<(.1){\S^{n+m}(\tau\times1)^*}
&\vee_{\L^{n+m}X}S^{n+m}\ar[r]_-{ev}&X\ar@{=>}"a";"b"_{\hat{\tau}}\ar@{=>}"e";"f"_F}\
\end{equation}
Any track horizontally composed with a trivial map becomes a trivial track, therefore in order to complete the
proof it is enough to check that
\begin{equation}\label{d}
\tau F= F(\tau\times 1).
\end{equation}
The track $F$ can be regarded as the relative homotopy class of a path in $\L^{n+m}X$ between two points comming from $E$.
Since $\varphi$ is a weak equivalence the whole path comes from $E$ up to homotopy. Moreover, by (\ref{e}) the
two possible composites in the diagram 
$$\xymatrix{E\ar[r]^-{\varphi}&\L^{n+m}X\ar@<.5ex>[r]^{(\tau\times
1_m)^*}\ar@<-.5ex>[r]_{\L^{n+m}\tau}&\L^{n+m}X}$$
coincide. Hence (d) follows.
\end{proof}

The following lemma is used in the proof of Theorem \ref{equit} above. The map $(\cdot)_{n+m}^{\sign(\tau)}$ and
the track $\hat{\tau}$ in the statement of the lemma can also be found in the proof of Theorem \ref{equit}.

\begin{lem}\label{kita}
Let $m\geq 1$. Given two pointed sets $A$, $B$ and a map $f\colon \vee_AS^1\r\vee_BS^1$ the Hopf invariant of the composite track $M(f)$
$$\xymatrix@C=50pt@R=50pt{\vee_AS^{n+m}\ar[r]^{\S^{n+m-1}f}
\ar[d]^{\tau\wedge S^m\wedge A}_{\;}="b"\ar@/_40pt/[d]_{(\cdot)^{\sign(\tau)}_{n+m}\wedge A}^{\;}="a"&
\vee_BS^{n+m}\ar[d]_{\tau\wedge S^m\wedge
B}^{\;}="c"\ar@/^40pt/[d]^{(\cdot)^{\sign(\tau)}_{n+m}\wedge B}_{\;}="d"\\
\vee_AS^{n+m}\ar[r]_{\S^{n+m-1}f}&\vee_BS^{n+m}\ar@{=>}"a";"b"^{\hat{\tau}^\vi\wedge A}
\ar@{=>}"c";"d"^{\hat{\tau}\wedge B}}$$
is defined by the equation, $a\in A$,
$$\hopf(M(f))(a)=\binom{\sign(\tau)}{2}H((\pi_1f)_\nill(a))\in\otimes^2\Z[B],$$
where $(\pi_1f)_\nill(a)$ denotes the image under $\pi_1f$ of the generator $a\in \grupo{A}\cong\pi_1(\vee_AS^1)$ 
projected to $\grupo{B}_\nill$.
\end{lem}

\begin{proof}
By using similar techniques to \cite[16]{2hg3} one can check that it is enough to prove the statement for
$f=(\cdot)^{-1}\colon S^1\r S^1$ the complex inversion. We
have that
$$\binom{\sign(\tau)}{2}H(-1)=\binom{\sign(\tau)}{2}\in\Z/2.$$ This coincides with $\hopf(M((\cdot)^{-1}))$ provided the following identity
holds in the positive pin group $Pin_+(n+m)$ as defined in \cite[6.6]{2hg2} 
$$e_{n+m}\hat{\tau}=(-1)^{\binom{\sign(\tau)}{2}}\hat{\tau}e_{n+m}.$$
The positive pin group is defined in \cite[6.6]{2hg2}  as a subgroup of units in the Clifford algebra
$C_+(n+m)$ with generators $e_1,\dots,e_{n+m}$ associated to the Euclidean scalar product in $\Real^{n+m}$.
See also \cite[6.11]{2hg2} to see how $\hat{\tau}$ defines an element in $Pin_+(n+m)$. Now the equation can be checked by using the defining
relations of $C_+(n+m)$, see \cite[6.6]{2hg2}, and the fact that $\hat{\tau}$ does not depend on $e_{n+m}$ as an
element in $C_+(n+m)$.
\end{proof}

We use Theorem \ref{equit} to prove the following result.

\begin{thm}\label{estatb}
For any morphism of symmetric spectra $f\colon X\wedge Y\r U$ with $U$ fibrant the composite of (\ref{level}) with the morphism
\begin{equation*}\tag{a}
\xymatrix{\Pi_{k+l+n+m,*}((X\wedge Y)_{k+l}|\L^{n+m}(X\wedge Y)_{0})\ar[d]\\
\Pi_{k+l+n+m,*}(U_{k+l}|\L^{n+m}U_{0})}
\end{equation*}
induced by $f$ defines a morphism in $\C{qpm}$ 
\begin{equation*}\tag{b}
\tilde{f}_{n,m}\colon \pi_{n,*}X\odot\pi_{m,*}Y\To\pi_{n+m,*}U
\end{equation*}
by taking colimits over $k,l\geq 0$.
\end{thm}

\begin{proof}
The map $f$ clearly induces a morphism (a) which is $\Z_\nill[\L^{n+m}f_0]$ on the $(0)$-level and which is induced
by $\Pi_{k+l+n+m,*}f$ on the $1$-level.

In order to check that the morphism (b) is defined in the colimit it is enough to check that
\begin{equation*}\tag{c}
f_*\sigma_*(1_{S^1}\wc 
(1_k\times\tau_{l,m}\times1_n)^*(\jmath_{k,l})_*(a\wc b))
\end{equation*}
is equal to both
\begin{equation*}\tag{d}
f_*(1_{1+k}\times\tau_{l,m}\times1_n)^*
(\jmath_{1+k,l})_*(\sigma_*(1_{S^1}\wc a)\wc b)),
\end{equation*}
\begin{equation*}\tag{e}
f_*(1_k\times\tau_{1+l,m}\times1_n)^*(\jmath_{k,1+l})_*(a\wc\sigma_*(1_{S^1}\wc b)),
\end{equation*}
at least for $k,l$ sufficiently large. We will see that equality $\text{(c)}=\text{(d)}$ holds for any $k$ and $l$,
while we have to require $l\geq1$ to check that $\text{(c)}=\text{(e)}$.

In order to check $\text{(c)}=\text{(d)}$ we consider the following equations 
\begin{eqnarray*}
\text{(c)}&=&f_*(1_{1+k}\times\tau_{l,m}\times1_n)^*\sigma_*(1_{S^1}\wc 
(\jmath_{k,l})_*(a\wc b))\\
&=&f_*(1_{1+k}\times\tau_{l,m}\times1_n)^*\sigma_*(S^1\wedge\jmath_{k,l})_*(1_{S^1}\wc 
(a\wc b))\\
&=&\text{(d)}.
\end{eqnarray*}
For the first equality we use the fact that the sign group action on secondary homotopy groups
of spaces is natural and that the smash product operation is equivariant, see \cite{2hg3} 7.1 and 8.13. 
For the second one we use again the naturality of the smash product operation. For the last one we use the fact that the following diagram commutes
$$\xymatrix@C=80pt{S^1\wedge X_{k}\wedge Y_{l}\ar[r]^{S^1\wedge\jmath_{k,l}}\ar[d]_{\sigma\wedge Y_{l}}&
S^1\wedge(X\wedge Y)_{k+l}\ar[d]^{\sigma}\\
X_{1+k}\wedge Y_{l}\ar[r]_{\jmath_{1+k,l}}&(X\wedge Y)_{1+k+l}}$$

Let us now check $\text{(c)}=\text{(e)}$ under the assumption $l\geq1$.
Here we need to assume that $U$ is fibrant in order to apply Theorem \ref{equit}. 
Let $\tau_\wedge\colon X_{k}\wedge S^1\cong S^1\wedge X_{k}$ be the symmetry isomorphism for the smash
product of pointed spaces. The commutativity rule for the smash product operation on secondary homotopy groups,
see \cite{2hg3} 7.3., and the fact that $H(1_{S^1})=0$ show that 
\begin{eqnarray*}
(\tau_\wedge)_*(a\wc 1_{S^1})&=&\tau_{1,k}^*(1_{S^1}\wa a)\\
&=&\tau_{1,k}^*(1_{S^1}\wc a).
\end{eqnarray*}
Using this equality together with the commutativity of
$$\xymatrix@C=3pt{X_{k}\wedge S^1\wedge Y_{l}\ar[r]^-{\tau_\wedge}\ar[d]_-{X_{k}\wedge\sigma}& 
**[r] S^1\wedge X_{k}\wedge Y_{l} 
\ar[rr]^-{S^1\wedge\jmath_{k,l}}&&
S^1\wedge(X\wedge Y)_{k+l}\ar[d]^-{\sigma}\\
X_{k}\wedge Y_{1+l}\ar[rr]_-{\jmath_{k,1+l}}&&**[l](X\wedge Y)_{k+1+l}&(X\wedge Y)_{1+k+l}
\ar[l]^-{\tau_{1,k}\times1_{l}}}$$
and the arguments used to check $\text{(c)}=\text{(d)}$
one can obtain the first equality in the following chain of equations.
\begin{eqnarray*}
\text{(e)}&=&f_*(1_k\times\tau_{1+l,n}\times1_m)^*\\
&&(\tau_{k+n,1}\times1_{l+m})^*(\tau_{1,k}\times1_{l})_*\sigma_*(1_{S^1}\wc (\jmath_{k,l})_*(a\wc b))\\
&=&(1_k\times\tau_{1+l,n}\times1_m)^*\\
&&(\tau_{k+n,1}\times1_{l+m})^*(\tau_{1,k}\times1_{l})_*\\
&&(1_{1+k}\times\tau_{n,l}\times1_m)^*\\
&&f_*\sigma_*(1_{S^1}\wc 
(1_k\times\tau_{l,n}\times1_m)^*(\jmath_{k,l})_*(a\wc b))\\
&=&(1_k\times\tau_{1+l,n}\times1_m)^*\\
&&(\tau_{k+n,1}\times1_{l+m})^*(\tau_{1,k}\times1_{l})_*(1_{1+k}\times\tau_{n,l}\times1_m)^*\\
&&(\tau_{k,1}\times1_{l})_*(\tau_{1,k}\times1_{l})_*\\
&&f_*\sigma_*(1_{S^1}\wc 
(1_k\times\tau_{l,n}\times1_m)^*(\jmath_{k,l})_*(a\wc b))\\
&\st{\text{(f)}}=&(1_k\times\tau_{1+l,n}\times1_m)^*\\
&&(\tau_{k+n,1}\times1_{l+m})^*(\tau_{1,k}\times1_{l})_*\\
&&(1_{1+k}\times\tau_{n,l}\times1_m)^*\\
&&(\tau_{k,1}\times1_{l})_*(\tau_{1,k}\times1_{l+n+m})^*\\
&&f_*\sigma_*(1_{S^1}\wc 
(1_k\times\tau_{l,n}\times1_m)^*(\jmath_{k,l})_*(a\wc b))\\
&=&\text{(c)}
\end{eqnarray*}
In the second and the third steps we introduce permutations which cancel each other. For the last equality we use
the identities
\begin{eqnarray*}
(\tau_{1,k}\times1_{l+n+m})(1_{1+k}\times\tau_{n,l}\times1_m)(\tau_{k+n,1}\times1_{l+m})(1_k\times\tau_{1+l,n}\times1_m)&=&1,\\
(\tau_{1,k}\times1_{l})(\tau_{k,1}\times1_{l})&=&1.
\end{eqnarray*}
For (f) we use the hypothesis $l\geq 1$ and Theorem \ref{equit} applied to the map
$$\varphi\colon\L^{n+m}U_{0}\r\L^{1+k+l+n+m}U_{1+k+l}$$ obtained by taking $\L^{n+m}$ on the adjoint of
$\sigma^{1+k+l}\colon S^{1+k+l}\wedge U_{0}\r U_{1+k+l}$. 
The map $\varphi$ is a weak equivalence since $U$ is
fibrant. Moreover, the adjoint of $\varphi$ is the composite
$$\psi\colon S^{1+k}\wedge S^l\wedge S^{n+m}\wedge\L^{n+m}U_{0}\To 
S^{1+k}\wedge S^l\wedge U_{0}
\st{\sigma_{1+k+l}}\To U_{1+k+l},$$ 
where the first arrow is obtained by taking $S^{1+k}\wedge S^l\wedge-$ on the adjoint 
$S^{n+m}\wedge\L^{n+m}U_{0}\r U_{0}$ of the identity map on $\L^{n+m}U_{0}$, therefore $\psi$ is
$\sym{1+k}$-equivariant with respect to, the action of $\sym{1+k}$ on $S^{1+k}$ by permutation of corrdinates, and
the pull-back action of $\sym{1+k+l}$ on $U_{1+k+l}$ along the inclusion
$-\times1_{l}\colon\sym{1+k}\r\sym{1+k+l}$.
\end{proof}

One now readily checks that the homomorphism $\tilde{f}_{n,m}$ given by Theorem \ref{estatb} proves Theorem 
\ref{sur2}.

\section{Proof of Theorems \ref{a2}, \ref{m2}, and \ref{ca2}}\label{pt1}\label{15}

Since the notation $\pi_{*,*}$ is also used in the statement of Theorems \ref{a2}, \ref{m2}, and \ref{ca2} we denote in this section the secondary homotopy groups of
spectra defined in the previous section by $\bar{\pi}_{*,*}$.

\begin{proof}[Proof of Theorem \ref{a2}]
Let $L$ be a fibrant replacement functor in the model category of ring spectra. We define the functor $\pi_{*,*}$
in the statement of Theorem \ref{a2} as the composite $\bar{\pi}_{*,*}L$. Given a ring spectrum $R$ the ring
multiplication $\mu\colon LR\wedge LR\r LR$ and
Theorem \ref{sur2} induce a quadratic pair algebra structure on $\bar{\pi}_{*,*}LR$ with multiplication
$$\mu_{*,*}\colon\bar{\pi}_{*,*}LR\odot\bar{\pi}_{*,*}LR\To\bar{\pi}_{*,*}LR.$$
Moreover, $\bar{\pi}_{0,*}S^0\cong\ol{\Z}_\nill$, and the unit $\nu\colon S^0\r LR$ of the ring spectrum
$LR$ induces the unit of the quadratic pair algebra $\bar{\pi}_{*,*}LR$, which is the image of $1\in
\Z=(\ol{\Z}_\nill)_0$ by the morphism 
$$\ol{\Z}_\nill\cong\bar{\pi}_{0,*}S^0\hookrightarrow\bar{\pi}_{*,*}S^0\st{\bar{\pi}_{*,*}\nu}\To\bar{\pi}_{*,*}LR.$$
The axioms of a quadratic pair algebra are satisfied since $\bar{\pi}_{*,*}LR$ is a monoid in the symmetric
monoidal category $\C{qpm}_H^\N$ of $\N$-graded $0$-good quadratic pair algebras. This follows from Theorem
\ref{a2} and the fact that secondary homotopy groups of pointed spaces form a lax monoidal functor to
$\C{qpm}_H^\N$, as we show in \cite{2hg3}.

Let us check the statement about Toda brackets and Massey products. For the definition of Toda brackets in
(\ref{toda2}) we assumed $R$ to be fibrant. In this proof we do not make this assumption, therefore we replace
$R$ in (\ref{toda2}) by the fibrant replacement $LR$. 

The right $LR$-module maps $\bar{a}$, $\bar{b}$, $\bar{c}$
in (\ref{toda2}) are the same as maps of pointed spaces $\bar{a}\colon S^p\r LR_0$, $\bar{b}\colon S^q\r LR_0$,
$\bar{c}\colon S^r\r LR_0$. These pointed maps are basis elements in the level $0$ of $\pi_{*,*}R$
representing $a$, $b$, $c$, as in Definition \ref{mp}. Moreover, the tracks $e$ and $f$ in (\ref{toda2}) can be
regarded as tracks between pointed maps $e\colon \mu_0(\bar{a}\wedge \bar{b})\rr 0$, $f\colon \mu_0(\bar{b}\wedge
\bar{c})\rr 0$. 
These tracks in $\C{Top}^*$ represent elements in $\Pi_{p+q,1}LR_0$ and $\Pi_{q+r,1}LR_0$, respectively, which project to
level $1$ elements in the colimit $\pi_{*,*}R$ also denoted by $e$ and $f$.
These elements satisfy $\partial(e)=\mu_0(\bar{a}\wedge \bar{b})=\bar{a}\cdot \bar{b}$ and $\partial(f)=
\mu_0(\bar{b}\wedge \bar{c})=\bar{b}\cdot \bar{c}$
as required in Definition \ref{mp}. Finally we notice that the element $-e\cdot\bar{c}+\bar{a}\cdot f$ defining
the Massey product corresponds exactly to the pasting of diagram (\ref{toda2}), which defines the Toda bracket,
hence they coincide.

The statement about the $k$-invariant follows from Proposition \ref{qlo2}.
\end{proof}

\begin{proof}[Proof of Theorem \ref{m2}]
Let $L$ be a fibrant replacement functor in the model category of ring spectra and let $L'$ be a fibrant
replacement in the model category of right $LR$-modules. We define the functor $\pi_{*,*}$
in the statement of Theorem \ref{m2} as the composite $\bar{\pi}_{*,*}L'(-\wedge_RLR)$. See \cite[22.2]{mcds} for
the definition of $\wedge_R$. Given a right $R$-module
$M$ the right action of the ring spectrum $LR$ on the module spectrum $L'(M\wedge_RLR)$, $$\gamma\colon L'(M\wedge_RLR)\wedge LR\r L'(M\wedge_RLR),$$ and
Theorem \ref{sur2} induce a right $\bar{\pi}_{*,*}LR$-module structure on $\bar{\pi}_{*,*}L'(M\wedge_RLR)$ with multiplication
$$\gamma_{*,*}\colon\bar{\pi}_{*,*}L'(M\wedge_RLR)\odot\bar{\pi}_{*,*}LR\To\bar{\pi}_{*,*}LR,$$
compare the proof of Theorem \ref{a2} above.
The axioms of a right $\bar{\pi}_{*,*}LR$-module are satisfied since $\bar{\pi}_{*,*}L'(M\wedge_RLR)$ is right module
over the monoid monoid $\bar{\pi}_{*,*}LR$ in the symmetric
monoidal category $\C{qpm}_H^\N$ of $\N$-graded $0$-good quadratic pair algebras. This follows from Theorem
\ref{a2} and the fact that secondary homotopy groups of pointed spaces form a lax monoidal functor to
$\C{qpm}_H^\N$, as we show in \cite{2hg3}.

The statements about Massey products and Toda brackets follow as in the proof of Theorem
\ref{a2}. The identification of the product morphism (\ref{mm}) is an easy exercise.
\end{proof}

\begin{rem}\label{expli}
Once defined the functor $\pi_{*,*}$ in Theorem \ref{m2} we explicitly indicate how to construct the morphism
(\ref{explia}) in the proof of Theorem \ref{mtb} from a trivial cofibration $f\colon R^{(n_1,\dots,n_k)}=(S^{n_1}\vee\cdots \vee
S^{n_k})\wedge R\st{\sim}\rightarrowtail X$.
\begin{equation*}
\xymatrix@R=10pt{(\pi_{*,*}R)^{(n_1,\dots,n_k)}\ar@{=}[d]\\
((\ol{\Z}_\nill,n_1)\vee \cdots\vee (\ol{\Z}_\nill,n_k))\odot \bar{\pi}_{*,*}LR\ar[d]^{\text{using $\imath_{n_j,0}$
in Remark
\ref{pinSn}}}\\
(\bar{\pi}_{n_1,*}S^{n_1}\vee \cdots\vee \bar{\pi}_{n_k,*}S^{n_k})\odot {\bar{\pi}}_{*,*}LR\ar@{^{(}->}[d]\\
(\bar{\pi}_{*,*}S^{n_1}\vee \cdots\vee \bar{\pi}_{*,*}S^{n_k})\odot {\bar{\pi}}_{*,*}LR\ar[d]\ar[d]\\
\bar{\pi}_{*,*}(S^{n_1}\vee \cdots\vee S^{n_k})\odot {\bar{\pi}}_{*,*}LR\ar[d]
\ar[d]^{\text{induced by }(S^{n_1}\vee \cdots\vee S^{n_k})\wedge LR\st{\sim}\rightarrowtail
L'((S^{n_1}\vee \cdots\vee S^{n_k})\wedge LR)\text{ and Thm. \ref{estatb}}}\\
\bar{\pi}_{*,*}L'((S^{n_1}\vee \cdots\vee S^{n_k})\wedge LR)\ar@{=}[d]\\
\pi_{*,*}((S^{n_1}\vee \cdots\vee S^{n_k})\wedge R)\ar[d]^{\pi_{*,*}f}\\
\pi_{*,*}X
}
\end{equation*}
\end{rem}

In the statement of the following lemma we consider the interval $1$-cell complex $K^1$, and the $2$-cell
complexes with the shape of a triangle $K^2$, a hexagon $K^3$, and a disc $K^4$.
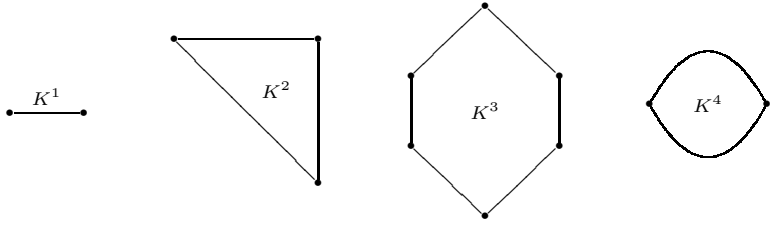
\begin{equation*}
\def\objectstyle{\scriptscriptstyle}
\begin{array}{c}
\entrymodifiers={+0}\xymatrix{\bullet\ar@{-}[r]^{K^1}&\bullet}
\end{array}
\qquad
\begin{array}{c}
\entrymodifiers={+0}\xymatrix@!=30pt{\bullet\ar@{-}[r]\ar@{-}[rd]&\bullet\ar@{-}[d]_{}="a"\\&\bullet
\ar@{}"1,1";"a"|<(.7){K^2}}
\end{array}
\qquad
\begin{array}{c}
\entrymodifiers={+0}\xymatrix{&\bullet\ar@{-}[ld]\ar@{-}[rd]&\\
\bullet\ar@{-}[d]\ar@{}[rrd]|{K^3}&&\bullet\ar@{-}[d]\\
\bullet\ar@{-}[rd]&&\bullet\ar@{-}[ld]\\&\bullet&}
\end{array}
\qquad
\begin{array}{c}
\entrymodifiers={+0}\xymatrix@!=20pt{\bullet\ar@/_20pt/@{-}[r]\ar@/^20pt/@{-}[r]\ar@{}[r]|-{K^4}&\bullet}
\end{array}
\end{equation*}
Moreover, given space $K$ we write $K_+$ for the disjoint union of $K$ with
an outer base-point. The smash product of a pointed space $K$ and a symmetric spectrum $X$ is the
symmetric spectrum defined by $(K\wedge X)_n=K\wedge X_n$. We denote by $\tau_\wedge$ to the
symmetry isomorphism for the smash product of spectra.

\begin{lem}\label{haxa}
Let $L$ be a fibrant replacement functor in the category of ring spectra. Given a commutative ring spectrum $Q$
the multiplication $\mu\colon LQ\wedge LQ\r LQ$ in $LQ$ needs not be commutative, but there is a commuting
homotopy $\alpha_1\colon K^1_+\wedge LQ\wedge LQ\r LQ$ which is coherent in the following sense: there are 
maps $\alpha_2\colon K^2_+\wedge LQ\wedge LQ\r LQ$ and
$\alpha_3\colon K^3_+\wedge LQ\wedge LQ\wedge LQ\r LQ$ defined over the boundary of $K^1$ and $K^2$ as indicated in the diagram below.
Moreover, the commuting homotopy is natural in the following sense: given a map of commutative ring spectra
$f\colon Q\r Q'$ if $\alpha_1$ and $\alpha'_1$ denote respective commuting homotopies then there is a map
$\alpha_4\colon K^4_+\wedge LQ\wedge LQ\r LQ'$ defined over the boundary of $K^4$ as indicated in the diagram below.
Furthermore, the map $\alpha_i$ can be chosen so that the restriction to $Q$ along the trivial cofibration
$Q\st{\sim}\rightarrowtail LQ$ is constant over $K_i$, $i=1,2,3,4$.
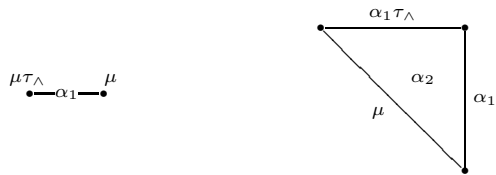
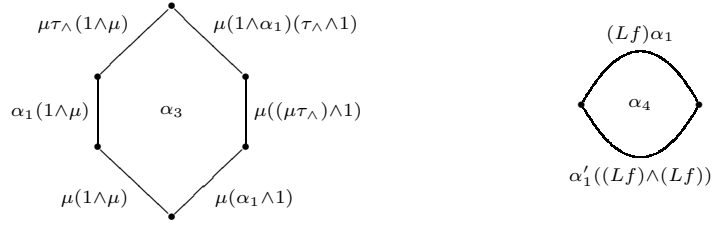
\begin{equation*}
\def\objectstyle{\scriptscriptstyle}
\begin{array}{c}
\entrymodifiers={+0}\xymatrix{\;\bullet\ar@{-}[r]|-{\alpha_1}^<(-.1){\mu\tau_\wedge}^<(1.1){\mu}&\bullet}
\end{array}
\qquad\qquad\qquad
\begin{array}{c}
\entrymodifiers={+0}\xymatrix@!=30pt{\bullet\ar@{-}[r]^-{\alpha_1\tau_\wedge}\ar@{-}[rd]_\mu&
\bullet\ar@{-}[d]^{\alpha_1}_{}="a"\\
&\bullet\ar@{}"1,1";"a"|<(.7){\alpha_2}}
\end{array}
\end{equation*}
\begin{equation*}
\def\objectstyle{\scriptscriptstyle}
\begin{array}{c}
\entrymodifiers={+0}\xymatrix{&\bullet\ar@{-}[ld]_{\mu\tau_\wedge(1\wedge\mu)}
\ar@{-}[rd]^{\mu(1\wedge\alpha_1)(\tau_\wedge\wedge1)}&\\
\bullet\ar@{-}[d]_{\alpha_1(1\wedge\mu)}\ar@{}[rrd]|{\alpha_3}&&\bullet\ar@{-}[d]^{\mu((\mu\tau_\wedge)\wedge1)}\\
\bullet\ar@{-}[rd]_{\mu(1\wedge\mu)}&&\bullet\ar@{-}[ld]^{\mu(\alpha_1\wedge1)}\\&\bullet&}
\end{array}
\qquad\qquad\qquad
\begin{array}{c}
\entrymodifiers={+0}\xymatrix@!=20pt{\bullet\ar@/_20pt/@{-}[r]_{\alpha_1'((Lf)\wedge(Lf))}
\ar@{}[r]|-{\alpha_4}\ar@/^20pt/@{-}[r]^{(Lf)\alpha_1}&\bullet}
\end{array}
\end{equation*}
Here the edges labelled with a map which only depends on $\mu$ and $\tau_\wedge$ are constant homotopies on the
indicated map.
\end{lem}


\begin{proof}
Let $P$ be the push-out 
\begin{equation*}
\xymatrix{\set{0,1}_+\wedge Q\wedge Q\ar@{>->}[r]\ar@{>->}[d]_\sim\ar@{}[rd]|{\text{push}}&
K^1_+\wedge Q\wedge Q\ar@{>->}[d]^\sim\\
\set{0,1}_+\wedge LQ\wedge LQ\ar@{>->}[r]&P}
\end{equation*}
The inclusion of the boundary $\set{0,1}\rightarrowtail K^1$ and $Q\st{\sim}\rightarrowtail LQ$ induce a trivial cofibration
$$j\colon P\st{\sim}\rightarrowtail K^1_+\wedge LQ\wedge LQ.$$

Consider the composite map 
$$a\colon K^1_+\wedge Q\wedge Q\To Q\wedge Q \st{\text{mult.}}\To Q\st{\sim}\rightarrowtail LQ$$
where the first arrow is given by the map collapsing the interval $K^1$ to a point.
Consider also the map
$$b\colon \set{0,1}_+\wedge LQ\wedge LQ\To LQ$$
which restricted to $\set{1}$ is $\mu$ and restricted to $\set{0}$ is $\mu\tau_\wedge$. The maps $a$ and $b$
induce a map
$$a\cup b\colon P\To LQ$$

The map $\alpha_1$ is defined as a lift of the following diagram
$$\xymatrix{P\ar[r]^-{a\cup b}\ar@{>->}[d]_\sim^j&LQ\ar[d]\\
K^1_+\wedge LQ\wedge LQ\ar[r]\ar@{-->}[ru]|{\alpha_1}&{*}}$$

Notice that the procedure for the construction of $\alpha_1$ consists of defining it over the restriction to $Q$ and on
the boundary of $K^1$ according to the statement. Then $\alpha_1$ is given by the
usual lifting diagram in a model category. The same procedure can be used for the definition of $\alpha_2$,
$\alpha_3$, and $\alpha_4$.
\end{proof}

\begin{proof}[Proof of Theorem \ref{ca2}]
Let $L$ be a fibrant replacement functor in the category of ring spectra. Given a commutative ring spectrum
$Q$ we have to promote the quadratic
pair algebra $\bar{\pi}_{*,*}LQ$ in the proof of Theorem \ref{a2} to an $E_\infty$-quadratic pair algebra. 
As we showed in the previous section the quadratic pair module $\bar{\pi}_{n,*}LQ$
carries an action of the symmetric
track group $\symt{n}$. The multiplication, defined by Theorem \ref{sur2}, is equivariant with respect to the sign group homomorphism
(\ref{sgx}), i.e. $\bar{\pi}_{*,*}LQ$ is a monoid in the category $\C{qpm}^{\symtt}_H$ of enhanced symmetric sequences of $0$-good
quadratic pair modules. Now consider the map $\alpha_1$ in the statement of Lemma \ref{haxa}.
By Theorem \ref{estatb}, and
using that the sphere spectrum $S$ is the unit of the smash product, the
map $\alpha_1$ yields a morphism 
\begin{equation*}
\bar{\pi}_{0,*}(K^1_+\wedge S)\odot\bar{\pi}_{*,*}LQ\odot\bar{\pi}_{*,*}LQ\To\bar{\pi}_{*,*}LQ.
\end{equation*}

The ``interval'' quadratic pair module $\mathbb{I}$ in \cite[5]{2hg3} 
is $0$-free. The basis at level $0$ is the set of vertices of $K_1$, i.e. $\set{0,1}$. 
On level $1$ it is generated by the edge $e$ of $K^1$ and $\partial(e)=-(0)+(1)$. This correspondence between the generators of
$\mathbb{I}$ and $K^1$ yields a morphism $\upsilon\colon\mathbb{I}\r \Pi_{0,*}K^1_+$, see \cite[19]{2hg3}.
Conposing $\upsilon$ with the projection to the colimit 
$\phi_0\colon \Pi_{0,*}K^1_+\r\bar{\pi}_{0,*}(K^1_+\wedge S)$ we obtain a new morphism
$$\mathbb{I}\odot\bar{\pi}_{*,*}LQ\odot\bar{\pi}_{*,*}LQ\To\bar{\pi}_{*,*}LQ.$$
Such a morphism is the same as a track to the multiplication of $\bar{\pi}_{*,*}LQ$ from the opposite
multiplication. The maps $\alpha_2$ and $\alpha_3$ in Lemma \ref{haxa} show that this track satisfies the
idempotence and the hexagon conditions, therefore $\bar{\pi}_{*,*}LQ$ is an $E_\infty$-quadratic pair algebra,
see Remark \ref{caen}. The map $\alpha_4$ in Lemma \ref{haxa} shows the naturaluty of $\smile_1$ with respect to morphisms between
commutative ring spectra.

For the statement about the smash product we notice that $\bar{a}$ in diagram (\ref{tcup1}) is the same as a
pointed map $\bar{a}\colon S^{2n}\r LQ_0$ which, regarded as a basis element in $\bar{\pi}_{2n,0}LQ$, represents
$a\in \pi_{2n}Q=\coker\partial$. Since $\bar{a}$ is in the basis $H(\bar{a})=0$. This simplifies the definition
of the algebraic cup-one square 
\begin{equation*}\tag{a}
Sq_1(a)=-(\bar{a}\cdot\bar{a})\cdot[\hat{\tau}]+\bar{a}\smile_1\bar{a},
\end{equation*} 
see Definition \ref{cup1}. Let us recall how these summands are explicitly defined.

An  element $\hat{\tau}\in\symt{n}$ as in Definition \ref{cup1} is the same 
as a track $\hat{\tau}$ as in diagram (\ref{tcup1}). Moreover, $\bar{a}\cdot\bar{a}=\mu_0(\bar{a}\wedge\bar{a})$.
The element $(\bar{a}\cdot\bar{a})\cdot[\hat{\tau}]$ is the projection to the colimit of
$\grupo{\mu_0(\bar{a}\wedge\bar{a}),\hat{\tau}}\in\Pi_{4n,1}LQ_0$ in the sense of \cite[4.5]{2hg2}.

The map $\alpha_1$ in Lemma \ref{haxa} is given on degree $0$ by a homotopy $(\alpha_1)_0$ to the product on
$LQ_0$ from the opposite product. Such a homotopy induces a track of quadratic pair modules, which is a function
$((\alpha_1)_0)_*\colon \Pi_{4n,0}LQ_0\r \Pi_{4n,1}LQ_0$. The construction of
$((\alpha_1)_0)_*(\mu_0(\bar{a}\wedge\bar{a}))$ can be found in \cite[7.3]{2hg1}. The projection to the colimit of this
element is $\bar{a}\smile_1\bar{a}$.

With this description of the summands in (a) one can check that (a) is the same as the pasting of diagram (\ref{tcup1}).
\end{proof}

\section{Proof of Theorem \ref{ca1}}\label{pt2}\label{16}

In order to derive Theorem \ref{ca1} from Theorem \ref{ca2} we need some technical results.

In Remark \ref{qia} we indicated how the category of $E_\infty$-pair algebras 
can be regarded as a full subcategory of the category of $E_\infty$-quadratic pair
algebras. This inclusion admits a reflection
$$(-)^\adc\colon\left(E_\infty\text{-quadratic pair algebras}\right)\To\left(E_\infty\text{-pair algebras}\right)$$
defined as follows. 

Let $C$ be an $E_\infty$-quadratic pair algebra. The quadratic pair module $C_{n,*}$ is a right $A(\symt{n})$-module for all $n\geq 0$, therefore
the additivization $C_{n,*}^\add$ in the sense of (\ref{ref2}) is a right $A(\symt{n})^\add$-module. The pair algebra $A(\symt{n})^\add$
is the inclusion of the two-sided ideal $I_{\sign}\subset\Z\sym{n}$ in the group-ring of the symmetric group which
is the kernel of the ring homomorphism
$\Z\sym{n}\r\Z$ induced by the sign homomorphism $\sign\colon\sym{n}\r\set{\pm1}$. The pair module $C^\adc_{n,*}$
is defined by the following diagram with exact rows
\begin{equation}\label{adc}
\xymatrix{C_{n,0}^\add\otimes_{\Z\sym{n}} I_{\sign}\ar[r]^-{\cdot}\ar@{=}[d]&C_{n,1}^\add\ar[d]^{\partial^\add}\ar@{->>}[r]&
C^\adc_{n,1}\ar[d]^{\partial^\adc}\\
C_{n,0}^\add\otimes_{\Z\sym{n}} I_{\sign}\ar[r]&C_{n,0}^\add\ar@{->>}[r]&
C^\adc_{n,0}}
\end{equation}
Here the arrow labeled with a dot $\cdot$ is defined by one of the multiplications given by the right $A(\symt{n})^\add$-module structure of
$C_{n,*}^\add$.
The natural projections $C_{n,*}\twoheadrightarrow C_{n,*}^\add\twoheadrightarrow C_{n,*}^\adc$ define an
$E_\infty$-quadratic pair algebra morphism $C\twoheadrightarrow C^\adc$, the unit of the reflection. The product
and the cup-one product in $C^\adc$ are determined by this fact.

\begin{lem}\label{sesi1}
Let $C$ be an $E_\infty$-quadratic pair algebra with trivial $k$-invariants such that $C_{n,0}^\add$ is a flat right $\sym{n}$-module for all
$n\geq 0$. Then the natural projection $C\twoheadrightarrow C^\adc$ is a quasi-isomorphism.
\end{lem}

\begin{proof}
Since $C_{n,*}$ has trivial $k$-invariant the natural projection $C_{n,*}\twoheadrightarrow C_{n,*}^\add$ is a
quasi-isomorphism by Proposition \ref{til}. Moreover, since $C_{n,0}^\add$ is flat as a right $\sym{n}$-module then the
lower horizontal arrow in the left of diagram (\ref{adc}) is injective. This implies that the upper one, labeled
$\cdot$, is also injective, hence the ``snake lemma'' applied to diagram (\ref{adc}) implies that
$C_{n,*}^\add\twoheadrightarrow C_{n,*}^\adc$ is also a quasi-isomorphism, and the lemma follows.
\end{proof}

The technical condition on $C_{n,0}^\add$ in the statement of Lemma \ref{sesi1} is satisfied, up to a blow-up, by
any $E_\infty$-quadratic pair algebra, as we show in the following lemma.

\begin{lem}\label{sesi2}
Let $C$ be a $0$-free $E_\infty$-quadratic pair algebra. There is a $0$-free quadratic pair algebra
$\widetilde{C}$ with $\widetilde{C}_{n,0}^\add$ a flat right $\sym{n}$-module for all $n\geq0$
and a natural projection $\widetilde{C}\twoheadrightarrow C$ which is a quasi-isomorphism.
\end{lem}

\begin{proof}
Since $C$ is $0$-free there are pointed sets $E_n$ such that $C_{n,(0)}=\Z_\nill[E_n]$, $n\geq 0$. 
The laws of a module over a quadratic pair algebra show that for any $e\in E_n$ and $g\in\sym{n}$ we have
$H(e\cdot[g])=0$. The unique elements of $\grupo{E_n}_\nill$ in the kernel of $H$ are the elements of the
basis $E_n$, compare the proof of \cite[Lemma 12]{qaI}. Therefore the action of $\sym{n}$ on $\grupo{E_n}_\nill$ is induced by
an action of $\sym{n}$ on the pointed set $E_n$.

The right $\sym{n}$-module $C_{n,0}^\add=\Z[E_n]$ needs not be flat. It would be flat if $E_n$ were a free pointed
$\sym{n}$-set. We are going to construct $\widetilde{C}$ by replacing $E_n$
by a free $\sym{n}$-set without changing $h_0$ and $h_1$. Let $$\varphi_n\colon E'_n=E_n\wedge(\sym{n}_+)\To E_n$$
be the pointed map induced by the action of $\sym{n}$ on $E_n$. This map is $\sym{n}$-equivariant with respect to
the action of $\sym{n}$ on $E'_n$ given by multiplication on the second coordinate of the smash product. 
Moreover, $E_*'$ is a graded monoid with multiplication given coordinatewise by the product in $E_*$ and the
cross product homomorphisms in (\ref{cp}). This multiplication is equivariant with respect to the cross product
homomorphisms. Let $\varphi\colon E'\r E$ be the graded map given by $\varphi_n$.

Since $E_n'$ is a free $\sym{n}$-set for all $n\geq 0$ the $E_\infty$-quadratic pair algebra 
$\widetilde{C}=\varphi^*C$ in Lemma \ref{kqt} and the natural projection
$\varphi_\#\colon\widetilde{C}\twoheadrightarrow C$ satisfy the properties stated in the lemma. 
\end{proof}

Now we are ready to prove Theorem \ref{ca1}.

\begin{proof}[Proof of Theorem \ref{ca1}]
The $E_\infty$-quadratic pair algebra $\pi_{*,*}Q$ in Theorem \ref{ca2} is $0$-free, therefore by Lemmas
\ref{sesi1} and \ref{sesi2} the functorial $E_\infty$-pair algebra
\begin{eqnarray*}
\pi^\adc_{*,*}Q&=&(\widetilde{\pi_{*,*}Q})^\adc
\end{eqnarray*}
proves Theorem \ref{ca1}. The natural zig-zag of quasi-isomorphisms between the pair algebras 
$\pi^\adc_{*,*}Q$ and $\pi_{*,*}^\add Q$ is obtained by applying the additivization functor (\ref{ref2}) to the
diagram
$$(\widetilde{\pi_{*,*}Q})^\adc\twoheadleftarrow\widetilde{\pi_{*,*}Q}\twoheadrightarrow \pi_{*,*}Q,$$
i.e. it is 
\begin{equation}\label{zz}
\pi^\adc_{*,*}Q\twoheadleftarrow(\widetilde{\pi_{*,*}Q})^\add\twoheadrightarrow \pi_{*,*}^\add Q.
\end{equation}
For this we recall that we have defined $\pi_{*,*}^\add Q$ as $(\pi_{*,*} Q)^\add$ in the proof of Theorem \ref{a1}.
\end{proof}

\begin{rem}\label{kiki}
In order to deduce Theorem \ref{ae1} from Theorem \ref{ae2} one extends the results of this section to algebras over an $E_\infty$-quadratic pair algebra $C$. More precisely, the natural projection $C\twoheadrightarrow C^{adc}$ induces an inclusion
$$\left(C^{adc}\text{-algebras}\right)\subset \left(C\text{-algebras}\right).$$
This inclusion admits a reflection
$$(-)^{ada}\colon\left(C\text{-algebras}\right)\To \left(C^{adc}\text{-algebras}\right)$$
sending a $C$-algebra $B$ to the $C^{adc}$-algebra $B^{ada}$ obtained from the additivization $B^{add}$ by killing the action of the symmetric groups as in (\ref{adc}). Lemmas \ref{sesi1} and \ref{sesi2} admit a generalized version for $C$-algebras. By using this one can check as in the proof of Theorem \ref{ca1} that the functor $\pi_{*,*}^{ada}$  in Theorem \ref{ae1} can be obtained taking first $\pi_{*,*}$ in Theorem \ref{ae2} and then applying $(-)^{ada}$ to a blow-up of the resulting $E_\infty$-quadratic pair algebra.
\end{rem}

\bibliographystyle{amsalpha}
\bibliography{Fernando}
\end{document}